\documentclass[11pt,reqno]{amsart}
\usepackage{amsmath,amssymb,amsthm,amsfonts,mathrsfs,color,cite} %times,
\usepackage{amsmath,amsfonts,mathtools}
\usepackage{amsthm}
\usepackage{graphicx}
\usepackage{color}
\usepackage{multicol}
\usepackage[margin=1in]{geometry}

\newtheorem{thm}{Theorem}[section]
\newtheorem{lemma}[thm]{Lemma}
\newtheorem{theorem}[thm]{Theorem}

\theoremstyle{definition}

\newtheorem{remark}{Remark}[section]
\newtheorem{notation}{Notation}[section]

\numberwithin{equation}{section}

\begin{document}

\title[Cauchy Problem of a System of Balance Laws]{Large Time Behavior and Diffusion Limit for a System of Balance Laws From Chemotaxis in Multi-dimensions}
\author{Tong Li}
\address[T. Li]{Department of Mathematics, University of Iowa, Iowa City, IA 52246, USA.}
\email{tong-li@uiowa.edu}
\author{Dehua Wang}
\address[D. Wang]{Department of Mathematics, University of Pittsburgh, Pittsburgh, PA 15260, USA.}
\email{dwang@math.pitt.edu}
\author{Fang Wang}
\address[F. Wang]{School of Mathematics and Statistics, Changsha University of Science and Technology, Changsha 410114, Hunan Province,   China.}
\email{wangfang1209@csust.edu.cn}
\author{Zhi-An Wang}
\address[Z.A. Wang]{Department of Applied Mathematics, Hong Kong Polytechnic University, Kowloon, Hong Kong.}
\email{mawza@polyu.edu.hk}
\author{Kun Zhao}
\address[K. Zhao]{Department of Mathematics, Tulane University, New Orleans, LA 70118, USA.}
\email{kzhao@tulane.edu}


\begin{abstract}
We consider the Cauchy problem for a system of balance laws derived from a chemotaxis model with singular sensitivity in multiple space dimensions. Utilizing energy methods, we first prove the global well-posedness of classical solutions to the Cauchy problem when only the energy of the first order spatial derivatives of the initial data is sufficiently small, and the solutions are shown to converge to the prescribed constant equilibrium states as time goes to infinity. Then we prove that the solutions of the fully dissipative model converge to those of the corresponding partially dissipative model when the chemical diffusion coefficient tends to zero.

\end{abstract}

\maketitle

\section{Introduction}

In this paper, we consider the system of balance laws:
\begin{equation}\label{hpbl}
\left\{
\begin{aligned}
&\partial_tp-\nabla\cdot(p\mathbf{v})=\Delta p,\\
&\partial_t\mathbf{v}-\nabla \left(p-\varepsilon|\mathbf{v}|^2\right)=\varepsilon\Delta\mathbf{v},
\end{aligned}
\right.\qquad \mathbf{x}\in\mathbb{R}^n,\ t>0,
\end{equation}
where $p(\mathbf{x},t)\in\mathbb{R}$ and $\mathbf{v}(\mathbf{x},t)\in\mathbb{R}^n$ are unknown functions, $n=2,3$, and $\varepsilon\ge0$ is a constant. The purpose of this paper is to study the {qualitative behavior, such as global well-posedness, long time behavior, and zero diffusion limit (as $\varepsilon\to0$), of classical solutions to the Cauchy problem} of \eqref{hpbl} in multiple space dimensions.

\subsection{Background} System \eqref{hpbl} {can be derived from the following chemotaxis model with logarithmic sensitivity}:
\begin{equation}\label{OS}
\left\{
\begin{aligned}
\partial_t u&=D\Delta u-\chi\nabla\cdot\left(u\nabla \log(c)\right),\\
\partial_t c&=\varepsilon\Delta c-\mu\, u\,c- \sigma c,
\end{aligned}
\right.\qquad \mathbf{x}\in\mathbb{R}^n,\ t>0,
\end{equation}
which was originally proposed in \cite{LS,OS} to describe the movement of chemotactic populations, such as myxobacteria, that deposit little- or non-diffusive chemical signals that modify the local environment for succeeding passages. System \eqref{OS} with $\sigma=0$ also appeared as a sub-model in \cite{LSN} to understand the {underlying mechanism of} tumor angiogenesis, and in particular the role of protease inhibitors in stopping angiogenesis.

System \eqref{OS} belongs to a family of nonlinear reaction-diffusion models, which are nowadays called the Keller-Segel type chemotaxis models, taking the canonical form:
\begin{equation}\label{KS}
\left\{
\begin{aligned}
\partial_t u&=D\Delta u-\chi\nabla\cdot\left(u\nabla\Phi(c)\right),\\
\partial_t c&=\varepsilon\Delta c+f(u,c).
\end{aligned}
\right.
\end{equation}
In the pioneering work \cite{KS71b}, {inspired by the existence of traveling bands in the chemotactic movement of {\it E. Coli} produced by Adler \cite{Adler}, Keller and Segel successfully reproduced such an experimental result through developing their original model by taking $\Phi(c)=\log(c)$, $f(u,c)=-\mu\, u\, c^m\,(0\le m<1)$ with $\chi,\mu>0$}. Since then the Keller-Segel model had provided a cornerstone for chemotaxis research, and its {capability of describing fundamental phenomena in chemotactic movement}, such as aggregation and uniform distribution (leveling out), inspired much of the later works investigating such a fundamental process in biochemistry. Numerous variations of the original Keller-Segel model have been developed to account for specific biological processes/environments involving chemotactic movements. Biologically, such kind of models describe the movement of certain biological organisms in response to the chemical signals that they release in the local environment for succeeding passages, while both entities are naturally diffusing and reacting (growing, dying, degrading, {\it et al}). Because of the biological background and analytical difficulties stemming from nonlinear advection (chemotaxis), the mathematical study of \eqref{KS} also attracted considerable attention from the community of nonlinear partial differential equations in recent decades. We refer the reader to the review papers \cite{Bellomo-review,HP09,DHreview1,Wang-DCDSB} and the references therein for more information.

System \eqref{OS} is a special case of \eqref{KS} when $\Phi(c)=\log(c)$ and $f(u,c)=-\mu\, u\, c-\sigma c$. The unknown functions and system parameters appearing in \eqref{OS} are interpreted as follows:
%\begin{itemize}
%\item $u(x,t)$: density of cellular population at position $x$ and time $t$,
%\item $c(x,t)$: concentration of chemical signal at position $x$ and time $t$,
%\item $D>0$: {diffusion coefficient of cellular density},
%\item $\chi\neq0$: coefficient of chemotactic sensitivity,
%\item $\varepsilon\ge0$: {diffusion coefficient of chemical concentration},
%\item $\mu\neq0$: {rate of density-dependent production/degradation of chemical signal},
%\item $\sigma>0$: {rate of natural degradation of chemical signal}.
%\end{itemize}
 $u(x,t)$ denotes the density of cellular population at position $x$ and time $t$,
 $c(x,t)$   the concentration of chemical signal at position $x$ and time $t$,
 $D>0$   the {diffusion coefficient of cellular density},
 $\chi\neq0$    the  coefficient of chemotactic sensitivity,
  $\varepsilon\ge0$   the {diffusion coefficient of chemical concentration},
 $\mu\neq0$   the {rate of density-dependent production/degradation of chemical signal}, and
 $\sigma>0$ denotes the  {rate of natural degradation of chemical signal}.
 We would like to emphasize that one of the most important parameters in \eqref{OS} {is $\chi$, whose sign dictates whether the chemotaxis is attractive ($\chi>0$) or repulsive ($\chi<0$)}, with $|\chi|$ measuring the strength of chemotactic response. It has been commonly acknowledged that the introduction of the nonlinear advection term in \eqref{KS} is the major contribution of the Keller-Segel type model, which captures the intrinsic features elucidating the underlying mechanisms of chemotactic movements.

Another important characteristic of \eqref{OS} is the logarithmic (singular) sensitivity function in the first equation, which entails that the chemotactic response of cellular population to chemical signal follows the Weber-Fechner's law, {a fundamental hypothesis in psychophysics, stating that subjective sensation is proportional to the logarithm of the stimulus intensity, and has played important roles in the modeling of biological processes} (cf. \cite{AL,BM,DLK,KJTW}). Indeed, the same sensitivity function was incorporated in the original Keller-Segel model \cite{KS71b}, whose significance was exemplified through demonstrating the model's capability of producing traveling wave solutions corroborating the experimental result reported in \cite{Adler}.

On the other hand, despite its importance in biological modeling, the possible singularity emanating from the logarithmic sensitivity function {brings about significant challenges} to the qualitative analysis of \eqref{OS}, from both the analytical and numerical perspectives. Soon after the model was initiated, it was discovered that (cf. \cite{LS}) the technical barrier raised by the singular sensitivity function can be removed by taking the Cole-Hopf transformation:
$$
\mathbf{V}=\nabla_{\mathbf{x}}\left(\log\left(e^{\sigma t}c(\mathbf{x},t)\right)\right),
$$
after which \eqref{OS} becomes a system of balance laws (also denoting $P\equiv u$):
\begin{equation}\label{TOS0}
\left\{
\begin{aligned}
&\partial_t P+\nabla\cdot\left(\chi P\mathbf{V}\right)=D\Delta P,\\
&\partial_t \mathbf{V}+\nabla\left(\mu P-\varepsilon|\mathbf{V}|^2\right)=\varepsilon\nabla(\nabla\cdot \mathbf{V}).
\end{aligned}
\right.
\end{equation}
Since $\Delta \mathbf{V}=\nabla(\nabla\cdot\mathbf{V})-\nabla\times(\nabla\times\mathbf{V})$ and $\mathbf{V}$ is a gradient field, for classical solutions, the system \eqref{TOS0} is equivalent to the following system of equations:
\begin{equation}\label{TOS}
\left\{
\begin{aligned}
&\partial_t P+\nabla\cdot\left(\chi P\mathbf{V}\right)=D\Delta P,\\
&\partial_t \mathbf{V}+\nabla\left(\mu P-\varepsilon|\mathbf{V}|^2\right)=\varepsilon\Delta \mathbf{V}.
\end{aligned}
\right.
\end{equation}
From the point of view of rigorous analysis we observe that the sign of the product of $\chi$ and $\mu$ plays an indispensable role in the qualitative study of the model. In fact, by applying the following re-scalings:
$$
t\rightarrow \frac{|\chi\mu|}{D}\,t,\quad\quad \mathbf{x}\rightarrow \frac{\sqrt{|\chi\mu|}}{D}\,\mathbf{x},\quad\quad \mathbf{V}\rightarrow -\mathrm{sign}(\chi)\sqrt{\frac{|\chi|}{|\mu|}}\,\mathbf{V},
$$
to the transformed system \eqref{TOS}, we obtain a clean version of the model:
\begin{equation}\label{RTOS}
\left\{
\begin{aligned}
&\partial_t P-\nabla\cdot\left(P\mathbf{V}\right)=\Delta P,\\
&\partial_t \mathbf{V}-\nabla\left(\mathrm{sign}(\chi\mu)P-\frac{\varepsilon}{\chi}|\mathbf{V}|^2\right)=\frac{\varepsilon}{D}\Delta \mathbf{V}.
\end{aligned}
\right.
\end{equation}
In the one-dimensional case, by a direction calculation, we can show that the characteristics (eigenvalues) associated with the flux on the left hand sides of the equations in \eqref{RTOS} are
$$
\lambda_\pm=\dfrac{\left(\dfrac{2\varepsilon}{\chi}-1\right)V\pm\sqrt{\left(\dfrac{2\varepsilon}{\chi}-1\right)^2V^2+4\,\mathrm{sign}(\chi\mu)P}}{2}.
$$
{Hence, the principle part of \eqref{RTOS} is hyperbolic when $\chi\mu>0$ in biologically relevant regimes where $P(\mathbf{x},t)>0$; while the system may change type when $\chi\mu<0$ (cf.\,\cite{LS})}. We refer the readers to \cite{LLW2016} for a recent study of the mixed type case, where the oscillatory traveling wave solutions are investigated. In this paper, we shall consider the case when $\chi\mu>0$, since otherwise
the possible change of type of the system may bring intractable difficulties to the underlying analysis in this paper.

Formally, {when $\chi>0$, the first equation of \eqref{OS} shows that the concentration gradient of the chemical signal drives the cellular population in the opposite direction of diffusion, indicating that the cellular population may aggregate as time evolves. On the other hand, the (exponentially) rapid degradation (due to $\mu>0$)} in the second equation of \eqref{OS} illustrates that the force driving the cellular population to aggregate is diminishing as time goes on. Hence, one may expect that the system will enter an equilibrium state in the long time run, due to the balance between cellular aggregation and chemical degradation. Similarly, when $\chi<0$ and $\mu<0$, because of the interaction between chemotactic repulsion and chemical production, the system is also expected to reach a steady state as time goes on. Collectively, when $\chi\mu>0$, finite time singularities are not anticipated to develop in the system \eqref{OS}, and the synergy of diffusion, chemotactic attraction/repulsion, and chemical degradation/production makes the dynamics of the model an intriguing problem to pursue.

In this paper, we aim to understand the dynamics of the {chemotaxis model with  logarithmic sensitivity, \eqref{OS},} through studying the qualitative behavior of solutions to the transformed system \eqref{RTOS} for fixed values of $\chi,\mu$ and $D$ when $\chi\mu>0$. Hence for brevity, we assume $\chi=\mu=D=1$ throughout the paper, {leading to} the following system of equations:
\begin{equation}\label{RTOS1}
\left\{
\begin{aligned}
&\partial_t P-\nabla\cdot\left(P\mathbf{V}\right)=\Delta P,\\
&\partial_t \mathbf{V}-\nabla\left(P-\varepsilon|\mathbf{V}|^2\right)=\varepsilon\Delta \mathbf{V},
\end{aligned}
\right.
\end{equation}
which is formally identical to the system \eqref{hpbl}. However, it should be mentioned that the system of balance laws \eqref{hpbl} is more general than the system \eqref{RTOS1}, due to the solution component $\mathbf{V}$ in the latter one is a gradient field (hence curl free). In this paper, we consider the general model \eqref{hpbl} and specify the condition under which the model automatically generates curl free solutions $\mathbf{V}$, from which one can recover the solutions to the original chemotaxis model \eqref{OS}.

\subsection{Literature Review and Motivations}

{To put things into perspective, now we would like to point out the existing results that are related to this work. When the spatial dimension is one, the following results for \eqref{hpbl} are available in the literature}:
%\begin{itemize}
%\item explicit and numerical solutions on finite intervals \cite{LS},
%\item shock wave formation for the Riemann problem on {$\mathbb{R}$} \cite{WH-shock},
%\item global well-posedness and long-time behavior of small data solutions on finite intervals \cite{ZhangZhu},
%\item local stability of traveling wave solutions on {$\mathbb{R}$} \cite{JLW,LiWang09,LLW,LW092,LWJDE,Li-Wang2011},
%\item global well-posedness of large data solutions on {$\mathbb{R}$} \cite{GXZZ},
%\item global well-posedness of large data solutions on finite intervals \cite{FFH},
%\item long-time behavior and chemical diffusion limit of large-amplitude classical solutions on finite intervals  \cite{LZ2015,LPZ2011,PWZZ,Tao-Wang2,WZ},
%\item long-time behavior, chemical diffusion limit and spatial analyticity of large data solutions on {$\mathbb{R}$} \cite{MWZ,LPZ2015},
%\item boundary layer formation and characterization of large data solutions on finite intervals \cite{HWZ2016,HLWW,LZ2015,PWZZ}.
%\end{itemize}
global well-posedness and large-time behavior with small or large data in \cite{ZhangZhu,GXZZ,FFH, LZ2015,LPZ2011,PWZZ,Tao-Wang2,WZ, MWZ,LPZ2015},
local stability of traveling wave solutions on {$\mathbb{R}$} in \cite{Choi3,JLW,LiWang09,LLW,LW092,LWJDE,Li-Wang2011,PW2018},
boundary layer formation and characterization of large data solutions  in \cite{HWZ2016,HLWW,LZ2015,PWZZ},
shock wave formation in \cite{WH-shock},   explicit and numerical solutions in \cite{LS}, and so on.
In particular, the results reported in \cite{LPZ2015,LZ2015,LPZ2011,MWZ,PWZZ,Tao-Wang2,WZ} indicate that when $\chi\mu>0$, no matter how strong the chemotactic sensitivity is and how large the energy of initial data is, the cellular population always distributes uniformly over space as time evolves.

We note that one of the main ingredients of the proofs constructed in \cite{LPZ2015,LZ2015,LPZ2011,MWZ,PWZZ,Tao-Wang2,WZ} is the implementation of the free energy (weak Lyapunov functional) associated with \eqref{hpbl}:
\begin{equation}\label{entropy1}
\frac{d}{dt}\left(\int E(p,\bar{p})dx+\|v\|_{L^2}^2\right)+\int\frac{(p_x)^2}{p}dx+\varepsilon\|v_x\|_{L^2}^2=0,
\end{equation}
where {$\bar{p}>0$} is a constant equilibrium state and the ``entropy expansion'' is defined by
$$
E(p,\bar{p})=[p\ln (p)-p]-[\bar{p}\ln(\bar{p})-\bar{p}]-\ln(\bar{p})(p-\bar{p}).
$$
The entropy type estimate \eqref{entropy1} lays down a foundation for the subsequent energy estimates that lead to the global well-posedness of large data classical solutions to the one-dimensional version of \eqref{hpbl} and the global stability of constant equilibrium states (uniform distribution of cellular population) associated with the model.

On the other hand, when the space dimension is greater than one, \eqref{entropy1} takes a different form:
\begin{equation}\label{entropy2}
\frac{d}{dt}\left(\int E(p,\bar{p})d\mathbf{x}+\|\mathbf{v}\|_{L^2}^2\right)+\int\frac{|\nabla p|^2}{p}d\mathbf{x}+\varepsilon\|\nabla\mathbf{v}\|_{L^2}^2=\varepsilon\int|\mathbf{v}|^2\nabla\cdot\mathbf{v}\ d\mathbf{x},
\end{equation}
where the integral on the right hand side does not vanish, and is not sign-preserving either. Moreover, by a direct calculation, we can show that \eqref{hpbl} is invariant under the scaling
$$
\left(p,\mathbf{v}\right)\rightarrow
(p^\xi,\mathbf{v}^\xi):=\left(\xi^2p\left(\xi \mathbf{x}, \xi^2t\right),\xi \mathbf{v}\left(\xi \mathbf{x}, \xi^2t\right)\right).
$$
Under the scaling, when the initial data are perturbed around the zero ground state, it holds that
$$
\|p_0^\xi\|^2_{L^2}=\xi^{4-n}\|p_0\|_{L^2}^2\ \ \ \ \ \mathrm{and}\ \ \ \ \ \|\mathbf{v}_0^\xi\|_{L^2}^2=\xi^{2-n}\|\mathbf{v}_0\|_{L^2}^2,
$$
which reveals that norm-inflation (especially for the $\mathbf{v}$-component) is not possible when $n\ge2$.

The aforementioned (unfavorable) features of the multi-dimensional version of \eqref{hpbl} brings substantial difficulties to the rigorous analysis of some of the fundamental properties of the model, such as global well-posedness of large data classical solutions. Unlike the one-dimensional case, to the authors' knowledge, only some scattered results have been obtained in the multi-dimensional case, such as,
the local well-posedness and blowup criteria of large data classical solutions to the Cauchy problem in $\mathbb{R}^2$ and $\mathbb{R}^3$ in \cite{FZ,LLZ},
  global well-posedness and long-time behavior of (full spectrum) small data classical solutions to the Cauchy problem in $\mathbb{R}^2$ and $\mathbb{R}^3$ in \cite{Hao,LLZ},
  global well-posedness of classical solutions to the Cauchy problem in $\mathbb{R}^3$ when only {$\|p_0-\bar{p}\|_{L^2}+\|\mathbf{v}_0\|_{H^1}$} is small and long-time behavior when {$\|p_0-\bar{p}\|_{H^2}+\|\mathbf{v}_0\|_{H^1}$} is small in \cite{DL},
  global well-posedness and long-time behavior of classical solutions to the Cauchy problem in $\mathbb{R}^3$ when only {$\|(p_0-\bar{p},\mathbf{v}_0)\|_{L^2}$} is small in \cite{PWZ},
 global well-posedness and long-time behavior of strong solutions to the Cauchy problem in $\mathbb{R}^2$ and $\mathbb{R}^3$ when only { $\|(p_0-\bar{p},\mathbf{v}_0)\|_{H^1}$} is small in \cite{WXY},
 global well-posedness and long-time behavior of classical solutions to the Cauchy problem in $\mathbb{R}^3$ when only { $\int E(p,\bar{p})d\mathbf{x}+\|\mathbf{v}\|_{L^2}^2$} (cf. \eqref{entropy2}) is small and in $\mathbb{R}^2$ when only $\|(p_0-\bar{p},\mathbf{v}_0)\|_{L^2}$ is small in \cite{WWZ},
 global well-posedness and long-time behavior of classical solutions on bounded domains in $\mathbb{R}^2$ and $\mathbb{R}^3$ when only { $\int E(p,\bar{p})d\mathbf{x}+\|\mathbf{v}\|_{L^2}^2$} is small in \cite{RWWZZ},
 global existence of intermediate weak solutions on bounded domains in $\mathbb{R}^2$ and $\mathbb{R}^3$ with Neumann boundary conditions in \cite{LS2015},
 global existence of generalized (weak) solutions on bounded domains in $\mathbb{R}^2$ with Neumann boundary conditions in \cite{Wink2} followed with a work addressing the eventual smoothness of solutions in \cite{Wink3},
  convergence of boundary layer solutions in the half plane in \cite{HWJMPA}, stability of planar traveling waves in a two-dimensional infinite cylindrical domain in \cite{Choi1, Choi2}, and so on.

By closely examining the results listed above we found that although they provide useful information for the basic understanding of the model in multi-dimensional spaces, none of them gives a positive answer to the question of global well-posedness of classical solutions when the initial data carry potentially large $L^2$-norm of the zeroth frequency of the perturbations around prescribed constant ground states. The current work is primarily motivated by such a fact.

The second fact that motivates this work is due to the observation that in certain biological environment involving chemotactic movement, the chemical signals deposited by the organism that modify the local environment for succeeding passages are little- or non-diffusive (cf. \cite{OS}), which is modeled by the smallness of the diffusion coefficient $\varepsilon$ in \eqref{hpbl}. Such a feature is also observed in the process of tumor angiogenesis as the interaction between vascular endothelial growth factor (VEGF), modeled by the function $c$ in \eqref{OS}, and vascular endothelial cells (VEC), modeled by the function $u$ in \eqref{OS}, is much more significant that the diffusion of VEGF \cite{LSN}. Hence, it is interesting and desirable to know whether the chemically diffusive (realistic) model can be approximated by the non-diffusive (ideal) one when $\varepsilon$ is small. Equivalently, the question is ``Does the solution of the slightly diffusive model converge (in certain topology) to the solution of the non-diffusive model, as $\varepsilon\to0$?" Such a topic has been investigated in some of the works mentioned above, in particular, in \cite{HWZ2016,LZ2015,MWZ,PWZZ,WZ} for the one-dimensional case, and in \cite{PWZ,RWWZZ,WWZ,WXY} for the multi-dimensional case. Nevertheless, the question of vanishing chemical diffusion coefficient limit of solutions to \eqref{hpbl} in multi-dimensional spaces with potentially large $L^2$-norm of the zeroth frequency of the perturbations around constant equilibrium states is widely open.

\subsection{Statement of Results}

Motivated by the facts mentioned above, we devote this paper to the study of the global well-posedness, long-time behavior and zero chemical diffusion limit of solutions to the Cauchy problem of \eqref{hpbl} in multiple space dimensions with initial data carrying potentially large $L^2$-norm of the zeroth frequency of the perturbations around prescribed constant ground states. The point of study of this paper is the following Cauchy (initial value) problem:
\begin{equation}\label{ivp}
\left\{
\begin{aligned}
&\partial_tp-\nabla\cdot(p\mathbf{v})-\bar{p}\nabla\cdot\mathbf{v}=\Delta p,\\
&\partial_t\mathbf{v}-\nabla \left(p-\varepsilon|\mathbf{v}|^2\right)=\varepsilon\Delta\mathbf{v},\\
&(p_0,\mathbf{v}_0)\in H^3(\mathbb{R}^n),\quad\quad p_0+\bar{p}>0,\quad\quad \nabla\times\mathbf{v}_0=\mathbf{0},
\end{aligned}
\right.
\end{equation}
where $\bar{p}>0$ is a constant, and $(p,\mathbf{v})$ denotes the perturbation of the original solution to \eqref{hpbl} around the constant state $(\bar{p},\mathbf{0})$. The main results of this paper are summarized in the following theorems.

\begin{notation}
Throughout the rest part of this paper, we use {\normalsize$\|\cdot\|$} to denote {\normalsize$\|\cdot\|_{L^2}$}.
\end{notation}

\begin{theorem}\label{thm1}
Let $n=3$, and consider the Cauchy problem \eqref{ivp}. Define
\begin{equation}\label{energy2}
\begin{aligned}
\kappa&\equiv 2\,(1+1/\bar{p})\left(\|\nabla p_0\|^2+\bar{p}\,\|\nabla\cdot\mathbf{v}_0\|^2\right),\\[2mm]
N_1&\equiv (1+1/\bar{p})\left(\|p_0\|^2+\bar{p}\,\|\mathbf{v}_0\|^2\right)+1.
\end{aligned}
\end{equation}
If
\begin{equation}\label{smallness3d}
N_1\kappa \le \min\left\{\left(\frac{1}{16c_1c_2}\right)^4,\ \frac{1}{8c_1c_2},\ \frac{1}{54\,c_3^4},\ 1 \right\},
\end{equation}
where $c_1$ and $c_2$ are the generic constants appearing in the Gagliardo-Nirenberg interpolation inequalities (cf. \eqref{gn1}--\eqref{gn2}), then there exists a unique solution to the Cauchy problem \eqref{ivp} for any $\varepsilon\ge0$, such that it holds that
\begin{equation}\label{ub}
\|(p,\mathbf{v})(t)\|^2_{H^3}+\int_0^t \left(\|(\nabla p,\sqrt{\varepsilon}\,\nabla\cdot\mathbf{v})(\tau)\|^2_{H^3}+\|\nabla\cdot\mathbf{v}(\tau)\|^2_{H^2}\right)d\tau \le C,\quad \forall\,t>0,
\end{equation}
where the positive constant $C$ is independent of $t$ and remains bounded as $\varepsilon\to0$. In addition, the following decay estimate holds:
\begin{equation}\label{de}
\lim_{t\to\infty}\left(\|(p,\mathbf{v})(t)\|^2_{W^{1,\infty}}+\|(\nabla p,\nabla\cdot\mathbf{v})(t)\|_{H^2}^2\right) = 0,
\end{equation}
for any {\normalsize$\varepsilon\ge0$}. Furthermore, let {\normalsize$(p^\varepsilon,\mathbf{v}^\varepsilon)$} and {\normalsize$(p^0,\mathbf{v}^0)$} be the solutions to the Cauchy problem with {\normalsize$\varepsilon>0$} and {\normalsize$\varepsilon=0$}, respectively, for the same initial data. Then for any {\normalsize$t>0$}, it holds that
\begin{equation}\label{dl}
\begin{aligned}
&\|(p^\varepsilon-p^0,\mathbf{v}^\varepsilon-\mathbf{v}^0)(t)\|^2_{H^1}\le D_1\,e^{t}\,\varepsilon^2,\\
&\|(\Delta p^\varepsilon-\Delta p^0,\Delta\mathbf{v}^\varepsilon-\Delta\mathbf{v}^0)(t)\|^2 \le D_2\,e^{t}(1+\varepsilon)\,\varepsilon,
\end{aligned}
\end{equation}
where the positive constants {\normalsize$D_1$} and {\normalsize$D_2$} are independent of $t$ and remain bounded as $\varepsilon\to 0$.
\end{theorem}

\begin{theorem}\label{thm2}
Let $n=2$, and consider the Cauchy problem \eqref{ivp}. Let
\begin{equation}\label{energy1a}
\begin{aligned}
H_1&\equiv \frac14\|p_0\|^2 + \left\|\frac{p_0}{2}-\frac{2\varepsilon+1}{6}p_0^2 \right\|^2 + \left(\frac{2(2\varepsilon+1)^2}{9}+\frac{8\varepsilon^2+8\varepsilon^3}{12\bar{p}}\right)\|p_0\|_{L^4}^4+\frac{\bar{p}}{4}\|\mathbf{v}_0\|^2+\\[2mm]
&\quad\ \ \frac{\bar{p}}{4}\left\|\left(1-\frac{2\varepsilon}{\bar{p}}p_0\right)\mathbf{v}_0 \right\|^2+\frac{\varepsilon^2}{\bar{p}}\|p_0\mathbf{v}_0\|^2+\frac{1}{12}\left(\frac{4\,\varepsilon+4\,\varepsilon^2}{\bar{p}}+9\,\varepsilon+32\,\varepsilon^3\right)\|\mathbf{v}_0\|_{L^4}^4,\\[2mm]
H_2&\equiv \frac{\bar{p}^2}{2}\,\|\nabla p_0\|^2+\bar{p}^3\|\nabla\cdot\mathbf{v}_0\|^2+\frac12\,\|\bar{p}\nabla p_0-p_0\nabla p_0\|^2+\frac12\,\|p_0\nabla p_0\|^2,
\end{aligned}
\end{equation}
and define
\begin{equation}\label{energy1b}
\begin{aligned}
M_1&\equiv 4(1+1/\bar{p}) H_1+1,\\[2mm]
\delta&\equiv \frac{2(2\bar{p}+1)}{\bar{p}^3}\exp\left\{\frac{4}{\bar{p}}\left(2916\,d_1^8\bar{p}^2M_1+3\,d_1^4\bar{p}^2+27\,d_1^8M_1\bar{p}^2+2916\,d_1^8d_5^4(\bar{p})^{-1}\right)H_1 \right\}H_2,
\end{aligned}
\end{equation}
where $d_1$ and $d_5$ are the generic constants appearing in the Gagliardo-Nirenberg interpolation inequalities (cf. \eqref{gn2d1}, \eqref{gn2d5}). If
\begin{equation}\label{smallness2d}
\max\{M_1\delta,\ M_1\delta^2,\ \delta\} \le \varrho,
\end{equation}
for some %absolute
positive constant $\varrho$ which is sufficiently small such that \eqref{B2} and \eqref{D2} are fulfilled, then there exists a unique solution to the Cauchy problem \eqref{ivp} for any $\varepsilon\ge0$, such that the solution obeys similar estimates as \eqref{ub}, \eqref{de} and \eqref{dl}.
\end{theorem}

\begin{remark}
In view of \eqref{energy2} and \eqref{energy1a}--\eqref{energy1b} we see that the smallness assumptions, \eqref{smallness3d} and \eqref{smallness2d}, can be realized by taking the $L^2$-norm of the first order spatial derivatives of the initial perturbations to be sufficiently small, while the $L^2$-norm of the initial perturbations can be potentially large. We provide explicit examples in the Appendix, which fulfill such requirements.
\end{remark}

\begin{remark}
%In the theorems above, we assumed $\mathbf{v}_0$ to be curl free, due to the consideration that the system of balance laws \eqref{hpbl} reduces to the transformed system \eqref{RTOS1} when $\mathbf{v}$ is curl free, and the latter one originates from the chemotaxis model, \eqref{OS}, through the transformation $\mathbf{V}=\nabla_{\mathbf{x}}\left(\log\left(e^{\sigma t}c(\mathbf{x},t)\right)\right)$.
In Theorems \ref{thm1} and \ref{thm2}, we assumed that $\mathbf{v}_0$ is curl free, which is natural since  the system of balance laws \eqref{hpbl} reduces to the transformed system \eqref{RTOS1} when $\mathbf{v}$ is curl free, and the latter one originates from the chemotaxis model, \eqref{OS}, through the transformation $\mathbf{V}=\nabla_{\mathbf{x}}\left(\log\left(e^{\sigma t}c(\mathbf{x},t)\right)\right)$.
Since the second equation in \eqref{hpbl}  automatically generates curl free solutions when $\mathbf{v}_0$ is curl free, due to $\nabla\times\mathbf{v}$ satisfies $(\nabla\times\mathbf{v})_t=\varepsilon\Delta(\nabla\times\mathbf{v})$, one can recover the solution to the original chemotaxis model \eqref{OS} from that to \eqref{hpbl} under the initial curl free condition.
\end{remark}

\subsection{Difficulties and Idea of Proof}

We prove Theorems \ref{thm1}--\ref{thm2} by developing {$L^p$}-based energy methods. Since we only assume the smallness of a fraction of the total Sobolev norm of the initial data, the major technical difficulty consists in closing the energy estimate for each individual frequency of the solution, without combining low and high frequencies (as is usually done in the case of classical solutions with small total Sobolev norm). Because the energy of the zeroth frequency part of the perturbation is allowed to be potentially large, the estimate for the zeroth frequency part is challenging, due to the lack of the Poincar\'{e} inequality in the whole space case. Moreover, because the Gagliardo-Nirenberg interpolation inequalities generate less powers of high frequencies of a function in {$\mathbb{R}^2$} than in {$\mathbb{R}^3$}, the proof of the two-dimensional case is considerably more involved than the three-dimensional case. We overcome the difficulties by terminating low frequencies through creating higher order nonlinearities, taking full advantage of the dissipation mechanisms and the smallness assumption on the first frequency, and utilizing various Gagliardo-Nirenberg interpolation inequalities. {In addition, since we also aim to establish the zero chemical diffusion limit of the chemically diffusive solution, it is vital to obtain the uniform \emph{$\varepsilon$-independent} energy estimates of the solution to \eqref{ivp} when $\varepsilon>0$. We reach the goal by deriving a linear and inhomogeneous damping equation for the spatial divergence of $\mathbf{v}$ and taking advantage of the structures of the equations, and again using various Gagliardo-Nirenberg interpolation inequalities.}

The rest of the paper is organized as follows. In Section 2 we give a complete proof of Theorem \ref{thm1}. Since the proof of Theorem \ref{thm2} is much more involved than that of Theorem \ref{thm1}, we present the main steps of the proof of Theorem \ref{thm2} in Section 3, while leave some tedious calculations in the Appendix.

\section{Proof of Theorem \ref{thm1}}

In this section we shall prove the Theorem \ref{thm1}. For the reader's convenience, we re-state the Cauchy problem:
\begin{equation}\label{e49h}
\left\{
\begin{aligned}
&\partial_tp-\nabla\cdot(p\mathbf{v})-\bar{p}\nabla\cdot\mathbf{v}=\Delta p,\\
&\partial_t\mathbf{v}-\nabla \left(p-\varepsilon|\mathbf{v}|^2\right)=\varepsilon\Delta\mathbf{v},\\
&(p_0,\mathbf{v}_0)\in H^3(\mathbb{R}^3),\quad\quad p_0+\bar{p}>0,\quad\quad \nabla\times\mathbf{v}_0=\mathbf{0},
\end{aligned}
\right.
\end{equation}
where $(p,\mathbf{v})$ denotes the perturbation of the original solution to \eqref{hpbl} around the constant state $(\bar{p},\mathbf{0})$.

First of all, we note that the local well-posedness of classical solutions to the Cauchy problem \eqref{e49h} can be established by applying Kawashima's theory on a general system of balance laws \cite{Kawashima}. Moreover, it follows from the maximum principle (cf. \cite{Friedman}) that the local solution satisfies $p+\bar{p}>0$ within its lifespan. We collect the results in the following:

\begin{lemma}[Local Well-posedness]\label{LWP}
Consider the Cauchy problem \eqref{e49h}. For any $\varepsilon\ge 0$, there exists a unique local-in-time solution such that $p+\bar{p}>0$ and
$$
\begin{aligned}
&(p-\bar{p},\mathbf{v})\in L^\infty([0,T_0];H^3(\mathbb{R}^3))\cap L^2([0,T_0];H^4(\mathbb{R}^3)),\quad &\textit{when}&\ \varepsilon>0;\\
&\qquad\begin{aligned}
p-\bar{p}&\in L^\infty([0,T_0];H^3(\mathbb{R}^3))\cap L^2([0,T_0];H^4(\mathbb{R}^3)),\\ \mathbf{v}&\in L^\infty([0,T_0];H^3(\mathbb{R}^3)),
\end{aligned}\quad &\textit{when}&\ \varepsilon=0,
\end{aligned}
$$
for some $T_0\in(0,\infty)$.
\end{lemma}

%Next, we shall derive {\it a priori} estimates for the local solution, under the conditions of Theorem \ref{thm1}, in order to extend it to be a global one.
We now establish {\it a priori} estimates for the local solution, in order to obtain a global solution in the  Theorem \ref{thm1}.
First of all, let us recall the following Gagliardo-Nirenberg interpolation inequalities:
\begin{align}
\|f\|_{L^3} &\le c_1\|\nabla f\|^{\frac12}\|f\|^{\frac12}, &\forall&\ f\in H^1(\mathbb{R}^3),\label{gn1}\\
\|f\|_{L^6} &\le c_2\|\nabla f\|, &\forall&\ f\in H^1(\mathbb{R}^3),\label{gn2}\\
\|f\|_{L^\infty}&\le  c_3\,\|\Delta f\|^\frac12\|\nabla f\|^\frac12, &\forall&\ f\in H^2(\mathbb{R}^3).\label{gn3}
\end{align}
Secondly, in addition to $\kappa$ and $N_1$ defined in \eqref{energy2}, let us define
\begin{equation}\label{N}
\begin{aligned}
N_2&\equiv (1+1/\bar{p})\exp\left\{\frac{3c_4}{2\bar{p}}\left(\|p_0\|_{H^1}^2+\bar{p}\,\|\mathbf{v}_0\|_{H^1}^2\right)\right\}\left(\|\Delta p_0\|^2+\bar{p}\,\|\Delta \mathbf{v}_0\|^2\right)+1,\\
\hat{N}&\equiv\frac{3c_4}{2(\bar{p}+1)}\left(\|p_0\|_{H^1}^2+\bar{p}\,\|\mathbf{v}_0\|_{H^1}^2\right)(N_2-1)+\left(\|\Delta p_0\|^2+\bar{p}\,\|\Delta \mathbf{v}_0\|^2\right),\\
N_3&\equiv (1+1/\bar{p})\exp\left\{\frac{3c_{5}}{2\bar{p}}\left(\|p_0\|_{H^1}^2+\bar{p}\,\|\mathbf{v}_0\|_{H^1}^2\right)\right\}\Big(\|\nabla\Delta p_0\|^2+\bar{p}\,\|\Delta (\nabla\cdot\mathbf{v}_0)\|^2+\\
&\qquad\qquad\qquad\qquad \qquad\qquad \qquad\qquad\qquad\qquad\quad   c_{5}(1+N_2)\hat{N}\Big)+1,
\end{aligned}
\end{equation}
where $c_4$ and $c_5$ are defined below by \eqref{c4} and \eqref{c5}, respectively. Then we observe that
\begin{equation}\label{ismall}
\begin{aligned}
&\|p_0\|^2+\|\mathbf{v}_0\|^2< N_1-1,\\
&\|\nabla p_0\|^2+\|\nabla\cdot\mathbf{v}_0\|^2< \frac{\kappa}{2},\\
&\|\Delta p_0\|^2+\|\Delta \mathbf{v}_0\|^2< N_2-1,\\
&\|\nabla\Delta p_0\|^2+\|\Delta \nabla\cdot\mathbf{v}_0\|^2 < N_3-1.
\end{aligned}
\end{equation}
Then it follows from Lemma \ref{LWP} that there exists $T_1\in (0,T_0]$, such that
\begin{equation}\label{n3xx}
\begin{aligned}
&\sup_{0\le t\le T_1}\left(\|p(t)\|^2+\|\mathbf{v}(t)\|^2\right)\le N_1,\\
&\sup_{0\le t\le T_1}\left(\|\nabla p(t)\|^2+\|\nabla\cdot\mathbf{v}(t)\|^2\right)\le \kappa,\\
&\sup_{0\le t\le T_1}\left(\|\Delta p(t)\|^2+\|\Delta \mathbf{v}(t)\|^2\right)\le N_2,\\
&\sup_{0\le t\le T_1}\left(\|\nabla\Delta p(t)\|^2+\|\Delta \nabla\cdot\mathbf{v}(t)\|^2\right)\le N_3.
\end{aligned}
\end{equation}
Next, we derive {\it a priori} estimates for the local solution within the time interval $[0,T_1]$. It was mentioned that under the initial curl free condition the system of equations, \eqref{e49h}, automatically produces curl free solutions. Hence, it suffices to estimate the divergence of {\normalsize$\mathbf{v}$}, i.e. {\normalsize$\nabla\cdot\mathbf{v}$}, in order to control the spatial derivatives of {\normalsize$\mathbf{v}$}. We first establish the estimate of the $L^2$-norm of the zeroth frequency part of the solution.

\subsection{$L^2$-Estimate}\label{3dl2}
By testing the equations in \eqref{e49h} with the targeting functions, and using \eqref{gn1}--\eqref{gn2}, we can show that
{\normalsize
\begin{equation}\label{n7h}
\begin{aligned}
&\ \frac12\frac{d}{dt}\left(\|p\|^2+\bar{p}\,\|\mathbf{v}\|^2\right)+\|\nabla p\|^2+\varepsilon\,\bar{p}\,\|\nabla \cdot\mathbf{v}\|^2\\
=&\ -\int_{\mathbb{R}^3}p(\mathbf{v}\cdot\nabla p)d\mathbf{x}+\varepsilon\,\bar{p}\int_{\mathbb{R}^3}|\mathbf{v}|^2\nabla\cdot\mathbf{v}\ d\mathbf{x}\\
\le&\ \|p\|_{L^6}\|\mathbf{v}\|_{L^3}\|\nabla p\|+\varepsilon\,\bar{p}\,\|\mathbf{v}\|_{L^3}\|\mathbf{v}\|_{L^6}\|\nabla\cdot\mathbf{v}\|\\
\le&\ c_1c_2\left( \|\nabla p\|\|\nabla\cdot\mathbf{v}\|^{\frac12}\|\mathbf{v}\|^{\frac12}\|\nabla p\|+\varepsilon\,\bar{p}\, \|\nabla\cdot\mathbf{v}\|^{\frac12}\|\mathbf{v}\|^{\frac12}\|\nabla\cdot \mathbf{v}\|^2\right)\\
\le&\ c_1c_2\,(N_1\kappa)^{\frac14}\left(\|\nabla p\|^2+\varepsilon\,\bar{p}\,\|\nabla\cdot\mathbf{v}\|^2\right),
\end{aligned}
\end{equation}}where we applied \eqref{n3xx} in the last inequality. Hence, when
\begin{equation}\label{kappa1}
N_1\kappa\le \left(\frac{1}{2c_1c_2}\right)^4,
\end{equation}
it holds that
{\normalsize
\begin{equation}\label{n8h}
\begin{aligned}
\frac{d}{dt}\left(\|p\|^2+\bar{p}\,\|\mathbf{v}\|^2\right)+\|\nabla p\|^2+\varepsilon\,\bar{p}\,\|\nabla \cdot\mathbf{v}\|^2\le 0,
\end{aligned}
\end{equation}}which yields
{\normalsize
\begin{equation}\label{n9m}
\|p(t)\|^2+\bar{p}\,\|\mathbf{v}(t)\|^2+\int_0^t\left(\|\nabla p(\tau)\|^2+\varepsilon\,\bar{p}\,\|\nabla \cdot\mathbf{v}(\tau)\|^2\right)d\tau\le \|p_0\|^2+\bar{p}\,\|\mathbf{v}_0\|^2.
\end{equation}}Therefore, in view of \eqref{energy2}, we see that
{\normalsize\begin{equation}\label{assump1}
\|p(t)\|^2+\|\mathbf{v}(t)\|^2 \le N_1-1.
\end{equation}}This completes the proof for the $L^2$-estimate. \hfill$\square$

Next, we make estimates on  the first order spatial derivatives of the solution.

\subsection{$H^1$-Estimate}\label{3dh1}
Testing the equations in \eqref{e49h} by the {\normalsize$-\Delta$} of the targeting functions, and using \eqref{gn3}, we can show that
{\normalsize
\begin{align}
&\frac12\frac{d}{dt}\left(\|\nabla p\|^2+\bar{p}\,\|\nabla\cdot \mathbf{v}\|^2\right)+\|\Delta p\|^2+\varepsilon\,\bar{p}\,\|\Delta\mathbf{v}\|^2\nonumber\\
=&\ -\int_{\mathbb{R}^3}\nabla\cdot(p\mathbf{v})\Delta p\ d\mathbf{x}+\varepsilon\,\bar{p}\int_{\mathbb{R}^3}\nabla(|\mathbf{v}|^2)\Delta\mathbf{v}\ d\mathbf{x}\nonumber\\
\le&\ \|p\|_{L^\infty}\|\nabla\cdot\mathbf{v}\|\|\Delta p\|+\|\nabla p\|_{L^6}\|\mathbf{v}\|_{L^3}\|\Delta p\|+2\,\varepsilon\,\bar{p}\,\|\mathbf{v}\|_{L^3}\|\nabla\mathbf{v}\|_{L^6}\|\Delta\mathbf{v}\| \nonumber\\
\le&\ c_3\,\|\nabla p\|^{\frac12}\|\Delta p\|^{\frac12}\|\nabla\cdot\mathbf{v}\|\|\Delta p\|+c_1c_2\,\|\mathbf{v}\|^{\frac12}\|\nabla\cdot\mathbf{v}\|^{\frac12}\|\Delta p\|^2+c_1c_2\,\varepsilon\,\bar{p}\,\|\nabla\cdot\mathbf{v}\|^{\frac12}\|\mathbf{v}\|^{\frac12}\|\Delta \mathbf{v}\|^2 \nonumber\\
\le&\ \left(\frac14+c_1c_2\,\|\mathbf{v}\|^{\frac12}\|\nabla\cdot\mathbf{v}\|^{\frac12}\right)\|\Delta p\|^2+\frac{27}{4}c_3^4\,\|\nabla p\|^2\|\nabla\cdot\mathbf{v}\|^4+c_1c_2\,\varepsilon\,\bar{p}\,\|\nabla\cdot\mathbf{v}\|^{\frac12}\|\mathbf{v}\|^{\frac12}\|\Delta \mathbf{v}\|^2\nonumber\\
\le&\ \left(\frac14+c_1c_2\,(N_1\kappa)^{\frac14}\right)\|\Delta p\|^2+\frac{27}{4}c_3^4\,\kappa^2\|\nabla p\|^2+c_1c_2\,\varepsilon\,\bar{p}\,(N_1\kappa)^{\frac14}\|\Delta \mathbf{v}\|^2,\label{n10h}
\end{align}}where we applied Young's inequality. Hence, when
\begin{equation}\label{kappa2}
N_1\kappa\le \left(\frac{1}{4c_1c_2}\right)^4,
\end{equation}
it holds that
{\normalsize
\begin{equation}\label{n11m}
\begin{aligned}
\frac{d}{dt}\left(\|\nabla p\|^2+\bar{p}\,\|\nabla\cdot \mathbf{v}\|^2\right)+\|\Delta p\|^2+\varepsilon\,\bar{p}\,\|\Delta\mathbf{v}\|^2\le \frac{27}{2}c_3^4\,\kappa^2\|\nabla p\|^2.
\end{aligned}
\end{equation}}Integrating \eqref{n11m} with respect to time and using \eqref{n9m}, we obtain
{\normalsize
\begin{equation}\label{n13m}
\begin{aligned}
&\ \|\nabla p(t)\|^2+\bar{p}\,\|\nabla\cdot\mathbf{v}(t)\|^2+\int_0^t\left(\|\Delta p(\tau)\|^2+\varepsilon\,\bar{p}\,\|\Delta\mathbf{v}(\tau)\|^2\right)d\tau \\
\le&\ \|\nabla p_0\|^2+\bar{p}\,\|\nabla\cdot\mathbf{v}_0\|^2+\frac{27}{2}c_3^4\,\kappa^2\left(\|p_0\|^2+\bar{p}\,\|\mathbf{v}_0\|^2\right),
\end{aligned}
\end{equation}}which implies
{\normalsize\begin{equation}\label{h1}
\|\nabla p(t)\|^2+\|\nabla\cdot\mathbf{v}(t)\|^2\le (1+1/\bar{p})\left[\|\nabla p_0\|^2+\bar{p}\,\|\nabla\cdot\mathbf{v}_0\|^2+\frac{27}{2}c_3^4\,\kappa^2\left(\|p_0\|^2+\bar{p}\,\|\mathbf{v}_0\|^2\right)\right].
\end{equation}}In view of \eqref{energy2}, we see that
$$
\begin{aligned}
\frac{27}{2}c_3^4\,\kappa^2\left(\|p_0\|^2+\bar{p}\,\|\mathbf{v}_0\|^2\right)&=54\,c_3^4(1+1/\bar{p})^2\left[\|\nabla p_0\|^2+\bar{p}\,\|\nabla\cdot\mathbf{v}_0\|^2\right]^2\left(\|p_0\|^2+\bar{p}\,\|\mathbf{v}_0\|^2\right)\\
&\le 27\,c_3^4\,\kappa\,N_1\left[\|\nabla p_0\|^2+\bar{p}\,\|\nabla\cdot\mathbf{v}_0\|^2\right]\\
&\le \frac12\left[\|\nabla p_0\|^2+\bar{p}\,\|\nabla\cdot\mathbf{v}_0\|^2\right],
\end{aligned}
$$
provided that
\begin{equation}\label{kappa2a}
N_1\kappa \le \frac{1}{54\ c_3^4}.
\end{equation}
Hence, we update \eqref{h1} as
{\normalsize\begin{equation}\label{h2}
\|\nabla p(t)\|^2+\|\nabla\cdot\mathbf{v}(t)\|^2\le \frac{3}{2}\,(1+1/\bar{p})\left[\|\nabla p_0\|^2+\bar{p}\,\|\nabla\cdot\mathbf{v}_0\|^2\right]=\frac34\kappa.
\end{equation}}In addition, we deduce from \eqref{n13m} that
{\normalsize
\begin{equation}\label{n13h}
\begin{aligned}
\int_0^t \left(\|\Delta p(\tau)\|^2+\varepsilon\,\bar{p}\,\|\Delta\mathbf{v}(\tau)\|^2\right)d\tau \le \frac32\left(\|\nabla p_0\|^2+\bar{p}\,\|\nabla\cdot\mathbf{v}_0\|^2\right),
\end{aligned}
\end{equation}}which will be utilized in the subsequent section. This completes the proof for the $H^1$-estimate. \hfill$\square$

\subsection{$H^2$-Estimate}\label{3dh2}
By computing the second order {\normalsize$L^2$} inner products, we can show that
{\normalsize
\begin{equation}\label{n14h}
\begin{aligned}
&\frac12\frac{d}{dt}\left(\|\Delta p\|^2+\bar{p}\,\|\Delta \mathbf{v}\|^2\right)+\|\nabla\Delta p\|^2+\varepsilon\,\bar{p}\,\|\Delta(\nabla\cdot\mathbf{v})\|^2\\
=&-\int_{\mathbb{R}^3}\nabla\left(\nabla\cdot(p\mathbf{v})\right)\cdot\nabla(\Delta p)\ d\mathbf{x}+\varepsilon\,\bar{p}\int_{\mathbb{R}^3}\Delta(|\mathbf{v}|^2)\Delta(\nabla\cdot\mathbf{v})\ d\mathbf{x}.
\end{aligned}
\end{equation}}For the first term on the RHS of \eqref{n14h}, by using H\"{o}lder, Gagliardo-Nirenberg and Young inequalities, we can show that
{\normalsize
\begin{equation*}
\begin{aligned}
&\ \left|-\int_{\mathbb{R}^3}\nabla\left(\nabla\cdot(p\mathbf{v})\right)\cdot\nabla(\Delta p)\ d\mathbf{x}\right|\\
\le&\ \left(\|p\|_{L^\infty}\|\Delta\mathbf{v}\|+\|\nabla p\|_{L^3}\|\nabla\mathbf{v}\|_{L^6}+\|\Delta p\|_{L^6}\|\mathbf{v}\|_{L^3}\right)\|\nabla\Delta p\|\\
\le&\ \left(c_3\|\nabla p\|^{\frac12}\|\Delta p\|^{\frac12}\|\Delta\mathbf{v}\|+c_1c_2\|\nabla p\|^{\frac12}\|\Delta p\|^{\frac12}\|\Delta \mathbf{v}\|+c_1c_2\|\nabla\Delta p\|\|\mathbf{v}\|^{\frac12}\|\nabla\cdot\mathbf{v}\|^{\frac12}\right)\|\nabla\Delta p\|\\
\le&\ \left(\frac14+c_1c_2\,(N_1\kappa)^{\frac14}\right)\|\nabla\Delta p\|^2 + \left(c_3+c_1c_2\right)^2\,\|\nabla p\|\|\Delta p\|\|\Delta\mathbf{v}\|^2\\
\le&\ \left(\frac14+c_1c_2\,(N_1\kappa)^{\frac14}\right)\|\nabla\Delta p\|^2 + \frac{\left(c_3+c_1c_2\right)^2}{2}\left(\|\nabla p\|^2+\|\Delta p\|^2\right)\|\Delta\mathbf{v}\|^2.
\end{aligned}
\end{equation*}}In a similar fashion, we can show that
{\normalsize
\begin{equation*}
\begin{aligned}
&\ \left|\varepsilon\,\bar{p}\int_{\mathbb{R}^3}\Delta(|\mathbf{v}|^2)\Delta(\nabla\cdot\mathbf{v})\ d\mathbf{x}\right|\\
\le&\ 2\,\varepsilon\,\bar{p}\,\left(\|\nabla\mathbf{v}\|_{L^3}\|\nabla\mathbf{v}\|_{L^6}+\|\mathbf{v}\|_{L^3}\|\nabla^2\mathbf{v}\|_{L^6}\right)\|\Delta(\nabla\cdot\mathbf{v})\|\\
\le&\ 2c_1c_2\,\varepsilon\,\bar{p}\, \left(\|\nabla\cdot\mathbf{v}\|^{\frac12}\|\Delta\mathbf{v}\|^{\frac32}+\|\mathbf{v}\|^{\frac12}\|\nabla\cdot\mathbf{v}\|^{\frac12}\|\Delta(\nabla\cdot\mathbf{v})\|\right)\|\Delta(\nabla\cdot\mathbf{v})\|\\
\le&\ \varepsilon\,\bar{p}\left(\frac14+2c_1c_2\,(N_1\kappa)^{\frac14}\right)\|\Delta(\nabla\cdot\mathbf{v})\|^2+ 4c_1^2c_2^2\,\varepsilon\,\bar{p}\,\|\nabla\cdot\mathbf{v}\|\|\Delta\mathbf{v}\|\|\Delta\mathbf{v}\|^2\\
\le&\ \varepsilon\,\bar{p}\left(\frac14+2c_1c_2\,(N_1\kappa)^{\frac14}\right)\|\Delta(\nabla\cdot\mathbf{v})\|^2+2c_1^2c_2^2\left(\varepsilon\,\bar{p}\,\|\nabla\cdot\mathbf{v}\|^2+\varepsilon\,\bar{p}\,\|\Delta\mathbf{v}\|^2\right)\|\Delta\mathbf{v}\|^2.
\end{aligned}
\end{equation*}}Hence, when
\begin{equation}\label{kappa3}
N_1\kappa\le \left(\frac{1}{8c_1c_2}\right)^4,
\end{equation}
it holds that
{\normalsize
\begin{equation}\label{n16h}
\begin{aligned}
&\ \frac{d}{dt}\left(\|\Delta p\|^2+\bar{p}\,\|\Delta \mathbf{v}\|^2\right)+\|\nabla\Delta p\|^2+\varepsilon\,\bar{p}\,\|\Delta(\nabla\cdot\mathbf{v})\|^2 \\
\le&\ \frac{c_4}{\bar{p}}\,\left(\|\nabla p\|^2+\|\Delta p\|^2+\varepsilon\,\bar{p}\,\|\nabla\cdot\mathbf{v}\|^2+\varepsilon\,\bar{p}\,\|\Delta\mathbf{v}\|^2\right)\left(\|\Delta p\|^2+\bar{p}\,\|\Delta\mathbf{v}\|^2\right),
\end{aligned}
\end{equation}}where
\begin{equation}\label{c4}
c_4=\max\left\{ (c_3+c_1c_2)^2, 4c_1^2c_2^2\right\}.
\end{equation}
Applying the Gronwall inequality to \eqref{n16h} and using \eqref{n9m} and \eqref{n13h}, we have
\begin{equation}\label{n17h}
\begin{aligned}
&\ \|\Delta p(t)\|^2+\bar{p}\,\|\Delta \mathbf{v}(t)\|^2 \\
\le &\ \exp\left\{\frac{c_4}{\bar{p}}\int_0^t\left(\|\nabla p\|^2+\|\Delta p\|^2+\varepsilon\,\bar{p}\,\|\nabla\cdot\mathbf{v}\|^2+\varepsilon\,\bar{p}\,\|\Delta\mathbf{v}\|^2\right)d\tau\right\}\left(\|\Delta p_0\|^2+\bar{p}\,\|\Delta \mathbf{v}_0\|^2\right)\\
\le &\ \exp\left\{\frac{3c_{4}}{2\bar{p}}\left(\|p_0\|_{H^1}^2+\bar{p}\,\|\mathbf{v}_0\|_{H^1}^2\right)\right\}\left(\|\Delta p_0\|^2+\bar{p}\,\|\Delta \mathbf{v}_0\|^2\right)\\
=&\ \dfrac{\bar{p}}{\bar{p}+1}(N_2-1),
\end{aligned}
\end{equation}
where we used the definition in \eqref{N}. From \eqref{n17h} we deduce that
\begin{equation}\label{h3}
\begin{aligned}
\|\Delta p(t)\|^2+\|\Delta \mathbf{v}(t)\|^2 \le N_2-1.
\end{aligned}
\end{equation}
In addition, by plugging \eqref{n17h} into \eqref{n16h}, then integrating the result with respect to $t$, we can show that
{\normalsize
\begin{equation}\label{n18h}
\begin{aligned}
&\ \int_0^t \left(\|\nabla\Delta p(\tau)\|^2+\varepsilon\,\bar{p}\,\|\Delta(\nabla\cdot\mathbf{v})(\tau)\|^2\right)d\tau \\
\le &\ \frac{3c_4}{2(\bar{p}+1)}\left(\|p_0\|_{H^1}^2+\bar{p}\,\|\mathbf{v}_0\|_{H^1}^2\right)(N_2-1)+\left(\|\Delta p_0\|^2+\bar{p}\,\|\Delta \mathbf{v}_0\|^2\right)=\hat{N}.
\end{aligned}
\end{equation}}We note that the constant $\hat{N}$ depends only on $\bar{p}$, the initial data and Gagliardo-Nirenberg constants. This completes the proof for the $H^2$-estimate. \hfill$\square$

\subsection{$H^3$-Estimate}\label{3dh3}
For the third order estimate, we can show that
\begin{equation}\label{H31}
\begin{aligned}
&\ \frac12\frac{d}{dt}\left(\|\nabla\Delta p\|^2+\bar{p}\,\|\Delta (\nabla\cdot\mathbf{v})\|^2\right)+\|\Delta^2 p\|^2+\varepsilon\,\bar{p}\,\|\Delta^2\mathbf{v}\|^2\\
\le &\ \left(\|p\|_{L^\infty}\|\Delta(\nabla\cdot\mathbf{v})\|+3\|\nabla p\|_{L^3}\|\Delta\mathbf{v}\|_{L^6}+3\|\Delta p\|_{L^6}\|\nabla\cdot\mathbf{v}\|_{L^3}+\|\nabla\Delta p\|\|\mathbf{v}\|_{L^\infty}\right)\|\Delta^2 p\|\\
&\quad +\varepsilon\,\bar{p}\left(2\|\Delta(\nabla\cdot\mathbf{v})\| \|\mathbf{v}\|_{L^\infty}+6\|\nabla\cdot\mathbf{v}\|_{L^3}\|\Delta\mathbf{v}\|_{L^6}\right)\|\Delta^2\mathbf{v}\|,
\end{aligned}
\end{equation}
where the terms on the right hand side are estimated by applying Gagliardo-Nirenberg inequalities and Sobolev embeddings as follows:
\begin{equation*}
\begin{aligned}
\bullet\quad \|p\|_{L^\infty}\|\Delta(\nabla\cdot\mathbf{v})\| \cdot \|\Delta^2 p\| &\le c_3\|\nabla p\|^{\frac12}\|\Delta p\|^{\frac12} \|\Delta(\nabla\cdot\mathbf{v})\| \cdot \|\Delta^2 p\| \\
&\le 2c_3^2\|\nabla p\|\cdot \|\Delta p\|\cdot \|\Delta(\nabla\cdot\mathbf{v})\|^2+\frac18\|\Delta^2 p\|^2\\
&\le c_3^2 (\|\nabla p\|^2+ \|\Delta p\|^2)\|\Delta(\nabla\cdot\mathbf{v})\|^2+\frac18\|\Delta^2 p\|^2,
\end{aligned}
\end{equation*}
\begin{equation*}
\begin{aligned}
\bullet\quad 3\|\nabla p\|_{L^3}\|\Delta\mathbf{v}\|_{L^6}\|\Delta^2 p\| &\le 3c_1c_2\|\nabla p\|^{\frac12}\|\Delta p\|^{\frac12} \|\Delta(\nabla\cdot\mathbf{v})\|\cdot \|\Delta^2 p\| \\
&\le 18c_1^2c_2^2\|\nabla p\|\cdot \|\Delta p\| \cdot\|\Delta(\nabla\cdot\mathbf{v})\|^2+\frac18\|\Delta^2 p\|^2\\
&\le 9c_1^2c_2^2(\|\nabla p\|^2+ \|\Delta p\|^2)\|\Delta(\nabla\cdot\mathbf{v})\|^2+\frac18\|\Delta^2 p\|^2,
\end{aligned}
\end{equation*}
\begin{equation*}
\begin{aligned}
\bullet\quad 3\|\Delta p\|_{L^6}\|\nabla\cdot\mathbf{v}\|_{L^3}\|\Delta^2 p\| &\le 3c_1c_2\|\nabla \Delta p\|\cdot \|\nabla\cdot\mathbf{v}\|^{\frac12} \|\Delta\mathbf{v}\|^{\frac12} \|\Delta^2 p\| \\
&\le 18c_1^2c_2^2\|\nabla\cdot\mathbf{v}\|\cdot \|\Delta\mathbf{v}\|\cdot \|\nabla \Delta p\|^2+\frac18\|\Delta^2 p\|^2,\\
&\le 9c_1^2c_2^2(\|\nabla\cdot\mathbf{v}\|^2+\|\Delta\mathbf{v}\|^2)\|\nabla \Delta p\|^2+\frac18\|\Delta^2 p\|^2,
\end{aligned}
\end{equation*}
\begin{equation*}
\begin{aligned}
\bullet\quad \|\nabla\Delta p\|\cdot \|\mathbf{v}\|_{L^\infty}\|\Delta^2 p\| &\le c_3\|\nabla \Delta p\|\cdot \|\nabla\cdot\mathbf{v}\|^{\frac12} \|\Delta\mathbf{v}\|^{\frac12} \|\Delta^2 p\| \\
&\le 2c_3^2\|\nabla\cdot\mathbf{v}\|\cdot \|\Delta\mathbf{v}\|\cdot \|\nabla \Delta p\|^2+\frac18\|\Delta^2 p\|^2,\\
&\le c_3^2(\|\nabla\cdot\mathbf{v}\|^2+\|\Delta\mathbf{v}\|^2)\|\nabla \Delta p\|^2+\frac18\|\Delta^2 p\|^2,
\end{aligned}
\end{equation*}
\begin{equation*}
\begin{aligned}
\bullet\quad 2\varepsilon\,\bar{p}\|\Delta(\nabla\cdot\mathbf{v})\| \cdot \|\mathbf{v}\|_{L^\infty}\|\Delta^2\mathbf{v}\|& \le 2c_3\varepsilon\,\bar{p}\|\nabla\cdot\mathbf{v}\|^{\frac12} \|\Delta\mathbf{v}\|^{\frac12} \|\Delta(\nabla\cdot\mathbf{v})\|\cdot\|\Delta^2\mathbf{v}\|\\
&\le 4c_3^2\varepsilon\,\bar{p}\|\nabla\cdot\mathbf{v}\|\cdot \|\Delta\mathbf{v}\|\cdot \|\Delta(\nabla\cdot\mathbf{v})\|^2+\frac14\varepsilon\,\bar{p}\|\Delta^2\mathbf{v}\|^2\\
&\le 2c_3^2\varepsilon\,\bar{p}(\|\nabla\cdot\mathbf{v}\|^2+ \|\Delta\mathbf{v}\|^2) \|\Delta(\nabla\cdot\mathbf{v})\|^2+\frac14\varepsilon\,\bar{p}\|\Delta^2\mathbf{v}\|^2,
\end{aligned}
\end{equation*}
\begin{equation*}
\begin{aligned}
\bullet\quad 6\varepsilon\,\bar{p}\|\nabla\cdot\mathbf{v}\|_{L^3}\|\Delta\mathbf{v}\|_{L^6}\|\Delta^2\mathbf{v}\|& \le 6c_1c_2\varepsilon\,\bar{p}\|\nabla\cdot\mathbf{v}\|^{\frac12} \|\Delta\mathbf{v}\|^{\frac12} \|\Delta(\nabla\cdot\mathbf{v})\|\cdot\|\Delta^2\mathbf{v}\|\\
&\le 36c_1^2c_2^2\varepsilon\,\bar{p}\|\nabla\cdot\mathbf{v}\|\cdot \|\Delta\mathbf{v}\|\cdot \|\Delta(\nabla\cdot\mathbf{v})\|^2+\frac14\varepsilon\,\bar{p}\|\Delta^2\mathbf{v}\|^2\\
&\le 18c_1^2c_2^2\varepsilon\,\bar{p}(\|\nabla\cdot\mathbf{v}\|^2+ \|\Delta\mathbf{v}\|^2) \|\Delta(\nabla\cdot\mathbf{v})\|^2+\frac14\varepsilon\,\bar{p}\|\Delta^2\mathbf{v}\|^2.
\end{aligned}
\end{equation*}
By substituting the above estimates into \eqref{H31}, we get
\begin{equation}\label{H32}
\begin{aligned}
&\ \frac{d}{dt}\left(\|\nabla\Delta p\|^2+\bar{p}\,\|\Delta (\nabla\cdot\mathbf{v})\|^2\right)+\|\Delta^2 p\|^2+\varepsilon\,\bar{p}\,\|\Delta^2\mathbf{v}\|^2\\
\le&\ c_5(\|\nabla p\|^2+ \|\Delta p\|^2)\|\Delta(\nabla\cdot\mathbf{v})\|^2+2c_5(\|\nabla\cdot\mathbf{v}\|^2+ \|\Delta\mathbf{v}\|^2)\times\\
&\qquad\qquad\qquad\qquad\qquad\qquad\qquad\qquad\qquad\qquad \big(\|\nabla\Delta p\|^2+\varepsilon\,\bar{p}\,\|\Delta(\nabla\cdot\mathbf{v})\|^2\big)\\
\le&\ \frac{c_5}{\bar{p}}(\|\nabla p\|^2+ \|\Delta p\|^2)\big(\|\nabla\Delta p\|^2+\bar{p}\,\|\Delta(\nabla\cdot\mathbf{v})\|^2\big)+c_5(\kappa+N_2)\times\\
&\qquad\qquad\qquad\qquad\qquad\qquad\qquad\qquad\qquad\qquad \big(\|\nabla\Delta p\|^2+\varepsilon\,\bar{p}\,\|\Delta(\nabla\cdot\mathbf{v})\|^2\big),
\end{aligned}
\end{equation}
where
\begin{equation}\label{c5}
c_5=9c_1^2c_2^2+c_3^2.
\end{equation}
By applying the Gronwall inequality to \eqref{H32} and using \eqref{n9m}, \eqref{n13h} and \eqref{n18h}, we have
\begin{equation}\label{H33}
\begin{aligned}
&\ \|\nabla\Delta p(t)\|^2+\bar{p}\,\|\Delta (\nabla\cdot\mathbf{v})(t)\|^2\\
\le &\ \exp\left\{\frac{c_{5}}{\bar{p}}\int_0^t (\|\nabla p(\tau)\|^2+\|\Delta p(\tau)\|^2)d\tau\right\}\bigg(\|\nabla\Delta p_0\|^2+\bar{p}\,\|\Delta (\nabla\cdot\mathbf{v}_0)\|^2+\\
&\qquad\qquad c_{5}(\kappa+N_2)\int_0^t\big(\|\nabla\Delta p(\tau)\|^2+\varepsilon\,\bar{p}\,\|\Delta(\nabla\cdot\mathbf{v})(\tau)\|^2\big)d\tau\bigg)\\
\le &\ \exp\left\{\frac{3c_{5}}{2\bar{p}}\left(\|p_0\|_{H^1}^2+\bar{p}\,\|\mathbf{v}_0\|_{H^1}^2\right)\right\}\Big(\|\nabla\Delta p_0\|^2+\bar{p}\,\|\Delta (\nabla\cdot\mathbf{v}_0)\|^2+c_{5}(\kappa+N_2)\hat{N}\Big).
\end{aligned}
\end{equation}
In view of \eqref{H33} and \eqref{N}, we see that
\begin{align}
&\ \|\nabla\Delta p(t)\|^2+\|\Delta (\nabla\cdot\mathbf{v})(t)\|^2\nonumber\\
\le &\ (1+1/\bar{p})\exp\left\{\frac{3c_{5}}{2\bar{p}}\left(\|p_0\|_{H^1}^2+\bar{p}\,\|\mathbf{v}_0\|_{H^1}^2\right)\right\}\Big(\|\nabla\Delta p_0\|^2+\bar{p}\,\|\Delta (\nabla\cdot\mathbf{v}_0)\|^2+c_{5}(\kappa+N_2)\hat{N}\Big)\nonumber\\
\le&\ (1+1/\bar{p})\exp\left\{\frac{3c_{5}}{2\bar{p}}\left(\|p_0\|_{H^1}^2+\bar{p}\,\|\mathbf{v}_0\|_{H^1}^2\right)\right\}\Big(\|\nabla\Delta p_0\|^2+\bar{p}\,\|\Delta (\nabla\cdot\mathbf{v}_0)\|^2+c_{5}(1+N_2)\hat{N}\Big)\nonumber\\
=&\ N_3-1,\label{h4}
\end{align}
provided that $\kappa\le 1$. This completes the proof for the $H^3$-estimate. \hfill$\square$

\subsection{Further Estimate for $\mathbf{v}$}\label{3dv}

From previous sections we see that the temporal integral of the spatial derivatives of {\normalsize$\mathbf{v}$} are inversely proportional to {\normalsize$\varepsilon$} (cf. \eqref{n9m}, \eqref{n13h}, \eqref{n18h}). In this section, we improve the temporal integrability of the spatial derivatives of {\normalsize$\mathbf{v}$} to be proportional to {\normalsize$\varepsilon$}, which will be used later for proving the zero chemical diffusion limit result. For this purpose, we take the divergence of the second equation of \eqref{e49h}, then combine the result with the first equation to get
{\normalsize
\begin{equation}\label{nn1h}
\partial_t (\nabla\cdot\mathbf{v})+\bar{p}\,\nabla\cdot\mathbf{v}=\varepsilon\Delta(\nabla\cdot\mathbf{v})+\partial_t p-\varepsilon\Delta(|\mathbf{v}|^2)-\nabla\cdot(p\mathbf{v}).
\end{equation}}Taking the {\normalsize$L^2$} inner product of \eqref{nn1h} with {\normalsize$\nabla\cdot\mathbf{v}$}, we have
{\normalsize
\begin{equation}\label{nn2h}
\begin{aligned}
&\ \frac12\frac{d}{dt}\|\nabla\cdot\mathbf{v}\|^2+\bar{p}\,\|\nabla\cdot\mathbf{v}\|^2+\varepsilon\,\|\Delta\mathbf{v}\|^2\\
=&\ \int_{\mathbb{R}^3}(\partial_t p)(\nabla\cdot\mathbf{v})d\mathbf{x}-\varepsilon\int_{\mathbb{R}^3}\Delta(|\mathbf{v}|^2)(\nabla\cdot\mathbf{v})d\mathbf{x}-\int_{\mathbb{R}^3}\left(\nabla\cdot(p\mathbf{v})\right)(\nabla\cdot\mathbf{v})d\mathbf{x}.
\end{aligned}
\end{equation}}We note that
{\normalsize
\begin{equation*}\label{nn3h}
\begin{aligned}
\int_{\mathbb{R}^3}(\partial_t p)(\nabla\cdot\mathbf{v})d\mathbf{x} &=\frac{d}{dt}\int_{\mathbb{R}^3}p(\nabla\cdot\mathbf{v})d\mathbf{x}-\int_{\mathbb{R}^3}p(\partial_t\nabla\cdot\mathbf{v})d\mathbf{x}\\
%&=\frac{d}{dt}\int_{\mathbb{R}^3}p(\nabla\cdot\mathbf{v})d\mathbf{x}-\int_{\mathbb{R}^3}p(\Delta p)d\mathbf{x}-\int_{\mathbb{R}^3}p\left(\varepsilon\Delta(\nabla\cdot\mathbf{v})-\varepsilon\Delta(|\mathbf{v}|^2)\right)d\mathbf{x}\\
&=\frac{d}{dt}\int_{\mathbb{R}^3}p(\nabla\cdot\mathbf{v})d\mathbf{x}+\|\nabla p\|^2-\int_{\mathbb{R}^3}p\left(\varepsilon\Delta(\nabla\cdot\mathbf{v})-\varepsilon\Delta(|\mathbf{v}|^2)\right)d\mathbf{x},
\end{aligned}
\end{equation*}}where we used the second equation of \eqref{e49h}. Then we update \eqref{nn2h} as
{\normalsize
\begin{equation}\label{nn4h}
\begin{aligned}
&\ \frac{d}{dt}\left(\frac12\|\nabla\cdot\mathbf{v}\|^2-\int_{\mathbb{R}^3}p(\nabla\cdot\mathbf{v})d\mathbf{x}\right)+\bar{p}\,\|\nabla\cdot\mathbf{v}\|^2+\varepsilon\,\|\Delta\mathbf{v}\|^2\\
%=&\ \|\nabla p\|^2-\varepsilon\int_{\mathbb{R}^3}\Delta(|\mathbf{v}|^2)(\nabla\cdot\mathbf{v})d\mathbf{x}-\int_{\mathbb{R}^3}\left(\nabla\cdot(p\mathbf{v})\right)(\nabla\cdot\mathbf{v})d\mathbf{x}-\\
%&\ \int_{\mathbb{R}^3}p\left(\varepsilon\Delta(\nabla\cdot\mathbf{v})-\varepsilon\Delta(|\mathbf{v}|^2)\right)d\mathbf{x}\\
=&\ \|\nabla p\|^2+\varepsilon\int_{\mathbb{R}^3}\nabla(|\mathbf{v}|^2)\cdot(\Delta\mathbf{v})d\mathbf{x}-\int_{\mathbb{R}^3}\left(\nabla\cdot(p\mathbf{v})\right)(\nabla\cdot\mathbf{v})d\mathbf{x}+\\
&\ \int_{\mathbb{R}^3}\nabla p\cdot\left(\varepsilon\nabla(\nabla\cdot\mathbf{v})-\varepsilon\nabla(|\mathbf{v}|^2)\right)d\mathbf{x}.
\end{aligned}
\end{equation}}For the second term on the RHS of \eqref{nn4h}, according to \eqref{gn1} and \eqref{gn2}, we have
{\normalsize\begin{equation*}\label{nn4ah}
\begin{aligned}
\left|\varepsilon\int_{\mathbb{R}^3}\nabla(|\mathbf{v}|^2)\cdot(\Delta\mathbf{v})d\mathbf{x} \right| &\le 2\varepsilon\|\mathbf{v}\|_{L^3}\|\nabla\cdot\mathbf{v}\|_{L^6}\|\Delta\mathbf{v}\|\\
&\le 2c_1c_2\,\varepsilon\,\|\mathbf{v}\|^{\frac12}\|\nabla\cdot\mathbf{v}\|^{\frac12}\|\Delta \mathbf{v}\|^2\\
&\le 2c_1c_2\,\varepsilon\,(N_1\kappa)^{\frac14}\|\Delta \mathbf{v}\|^2.
\end{aligned}
\end{equation*}}For the third term on the RHS of \eqref{nn4h}, by using similar arguments as in \eqref{n10h}, we can show that
{\normalsize
\begin{equation*}\label{nn5h}
\begin{aligned}
\left| -\int_{\mathbb{R}^3}\left(\nabla\cdot(p\mathbf{v})\right)(\nabla\cdot\mathbf{v})d\mathbf{x}\right| &\le \left(\|\nabla p\|_{L^3}\|\mathbf{v}\|_{L^6}+ \|p\|_{L^\infty}\|\|\nabla\cdot\mathbf{v}\|\|\right)\|\nabla\cdot\mathbf{v}\|\\
&\le (c_1c_2+c_3)\|\nabla p\|^{\frac12}\|\Delta p\|^{\frac12}\|\nabla\cdot\mathbf{v}\|^2\\
&\le \frac{(c_1c_2+c_3)^2}{2\bar{p}}\|\nabla p\|\|\Delta p\|\|\nabla\cdot\mathbf{v}\|^2+\frac{\bar{p}}{2}\,\|\nabla\cdot\mathbf{v}\|^2\\
&\le  \frac{(c_1c_2+c_3)^2}{\bar{p}}\left(\|\nabla p\|^2+\|\Delta p\|^2\right)\kappa+\frac{\bar{p}}{2}\,\|\nabla\cdot\mathbf{v}\|^2\\
&\le  \frac{(c_1c_2+c_3)^2}{\bar{p}}\left(\|\nabla p\|^2+\|\Delta p\|^2\right)+\frac{\bar{p}}{2}\,\|\nabla\cdot\mathbf{v}\|^2,
\end{aligned}
\end{equation*}}provided that $\kappa\le 1$. For the fourth term on the RHS of \eqref{nn4h}, we can show that
{\normalsize
\begin{equation*}\label{nn5ah}
\begin{aligned}
\left|\int_{\mathbb{R}^3}\nabla p\cdot\left(\varepsilon\nabla(\nabla\cdot\mathbf{v})-\varepsilon\nabla(|\mathbf{v}|^2)\right)d\mathbf{x}\right| &\le \varepsilon\,\|\nabla p\|\|\nabla(\nabla\cdot\mathbf{v})\|+2\,\varepsilon\,\|\nabla p\|\|\mathbf{v}\|_{L^3}\|\nabla\mathbf{v}\|_{L^6}\\
&\le 2\,\varepsilon\,\|\nabla p\|^2+\frac{\varepsilon}{4}\,\|\Delta\mathbf{v}\|^2+\varepsilon\,\|\mathbf{v}\|_{L^3}^2\|\nabla\mathbf{v}\|_{L^6}^2 \\
&\le 2\,\varepsilon\,\|\nabla p\|^2+\frac{\varepsilon}{4}\,\|\Delta\mathbf{v}\|^2+c_1c_2\,\varepsilon\,\|\mathbf{v}\|\|\nabla\cdot\mathbf{v}\|\|\Delta\mathbf{v}\|^2\\
&\le 2\,\varepsilon\,\|\nabla p\|^2+\frac{\varepsilon}{4}\,\|\Delta\mathbf{v}\|^2+c_1c_2\,\varepsilon\,(N_1\kappa)\|\Delta\mathbf{v}\|^2.
\end{aligned}
\end{equation*}}Hence, when
\begin{equation}\label{kappa4}
N_1\kappa \le \min\left\{\left(\frac{1}{16c_1c_2}\right)^4,\frac{1}{8c_1c_2} \right\},
\end{equation}
we update \eqref{nn4h} as
{\normalsize
\begin{equation}\label{nn6h}
\begin{aligned}
&\ \frac{d}{dt}\left(\frac12\|\nabla\cdot\mathbf{v}\|^2-\int_{\mathbb{R}^3}p(\nabla\cdot\mathbf{v})d\mathbf{x}\right)+\frac{\bar{p}}{2}\,\|\nabla\cdot\mathbf{v}\|^2 +\frac{\varepsilon}{2}\,\|\Delta\mathbf{v}\|^2\\
\le&\ (1+2\varepsilon)\|\nabla p\|^2+\frac{(c_1c_2+c_3)^2}{\bar{p}}\left(\|\nabla p\|^2+\|\Delta p\|^2\right).
\end{aligned}
\end{equation}}By multiplying \eqref{n8h} by 2, then adding the result to \eqref{nn6h}, we have
{\normalsize
\begin{equation}\label{nn7h}
\begin{aligned}
&\ \frac{d}{dt}E(t)+\frac{\bar{p}}{2}\,\|\nabla\cdot\mathbf{v}\|^2+\frac{\varepsilon}{2}\,\|\Delta\mathbf{v}\|^2+\|\nabla p\|^2 +2\,\varepsilon\,\bar{p}\,\|\nabla\cdot\mathbf{v}\|^2 \\
\le &\ \left(2\varepsilon+\frac{(c_1c_2+c_3)^2}{\bar{p}}\right)\|\nabla p\|^2+\frac{(c_1c_2+c_3)^2}{\bar{p}}\|\Delta p\|^2,
\end{aligned}
\end{equation}}where
{\normalsize$$
\begin{aligned}
E(t)&=\frac12\,\|\nabla\cdot\mathbf{v}\|^2-\int_{\mathbb{R}^3}p(\nabla\cdot\mathbf{v})d\mathbf{x}+2\,\|p\|^2+2\,\bar{p}\,\|\mathbf{v}\|^2\\
&=\frac14\,\|\nabla\cdot\mathbf{v}\|^2+\left\|\frac12\nabla\cdot\mathbf{v}-p\right\|^2+\|p\|^2+2\,\bar{p}\,\|\mathbf{v}\|^2.
\end{aligned}
$$}Integrating \eqref{nn7h} with respect to time and using \eqref{n9m} and \eqref{n13h} then yield, in particular,
{\normalsize\begin{equation}\label{nn8h}
\begin{aligned}
\frac{\bar{p}}{2}\int_0^t\|\nabla\cdot\mathbf{v}(\tau)\|^2d\tau \le E(0)+&\left(2\varepsilon+\frac{(c_1c_2+c_3)^2}{\bar{p}}\right)\left(\|p_0\|^2+\bar{p}\,\|\mathbf{v}_0\|^2\right) + \\
&\quad \frac{3(c_1c_2+c_3)^2}{2\bar{p}}\left(\|\nabla p_0\|^2+\bar{p}\,\|\nabla\cdot\mathbf{v}_0\|^2\right),
\end{aligned}
\end{equation}}where the constant on the RHS is independent of {\normalsize$t$} and remains bounded as {\normalsize$\varepsilon\to0$}. In a completely similar fashion, by working on the non-homogeneous damped equation \eqref{nn1h} and using the $H^1$-, $H^2$- and $H^3$-estimates established in the previous sections, we can show that
{\normalsize
\begin{equation*}
\int_0^t \left(\|\nabla(\nabla\cdot\mathbf{v})(\tau)\|^2+\|\Delta(\nabla\cdot\mathbf{v})(\tau)\|^2\right)d\tau\le C(c_1,c_2,c_3,\bar{p},p_0,\mathbf{v}_0,\varepsilon).
\end{equation*}}We omit the technical details to simplify the presentation. \hfill$\square$

{It is worth mentioning that all the constants appearing in the energy estimates in Sections 2.1--2.5 remain bounded as {\normalsize$\varepsilon\to0$}, which allows us to establish the global well-posedness of solutions to \eqref{ivp} when $\varepsilon=0$. Indeed, it can be readily checked that by repeating the arguments in Sections 2.1--2.5, one can establish similar energy estimates for the solution to \eqref{ivp} when $\varepsilon=0$, and in particular the solution satisfies \eqref{ub} with $\varepsilon=0$. More importantly, the energy estimates derived in Sections 2.1--2.5 allow us to take the zero chemical diffusion limit of the solution in a later section.} Furthermore, in view of \eqref{kappa1}, \eqref{kappa2}, \eqref{kappa2a}, \eqref{kappa3}, \eqref{kappa4} and the derivation of \eqref{h4}, we see that all of the energy estimates derived in Sections \ref{3dl2}--\ref{3dh3} are valid when
\begin{equation}\label{kappa5}
N_1\kappa \le \min\left\{\left(\frac{1}{16c_1c_2}\right)^4,\ \frac{1}{8c_1c_2},\ \frac{1}{54\,c_3^4},\ 1 \right\}.
\end{equation}
Note that \eqref{kappa5} is the same smallness condition as stated in Theorem \ref{thm1}.

\subsection{Global Well-posedness} We now prove the global well-posedness of classical solutions to \eqref{e49h} in  Theorem \ref{thm1}. %For this purpose,
First,  from \eqref{assump1}, \eqref{h2}, \eqref{h3} and \eqref{h4}, we have
\begin{equation}\label{ismall1}
\begin{aligned}
&\|p(t)\|^2+\|\mathbf{v}(t)\|^2 \le N_1-1,\\
&\|\nabla p(t)\|^2+\|\nabla\cdot\mathbf{v}(t)\|^2 \le \frac34\kappa,\\
&\|\Delta p(t)\|^2+\|\Delta\mathbf{v}(t)\|^2 \le N_2-1,\\
&\|\nabla\Delta p(t)\|^2+\|\Delta (\nabla\cdot\mathbf{v})(t)\|^2 \le N_3-1,
\end{aligned}\qquad\forall\ t\in[0,T_1],
\end{equation}
which, in particular, imply that
\begin{equation}\label{ismall2}
\begin{aligned}
&\|p(T_1)\|^2+\|\mathbf{v}(T_1)\|^2 \le N_1-1,\\
&\|\nabla p(T_1)\|^2+\|\nabla\cdot\mathbf{v}(T_1)\|^2 \le \frac34\kappa,\\
&\|\Delta p(T_1)\|^2+\|\Delta\mathbf{v}(T_1)\|^2 \le N_2-1,\\
&\|\nabla\Delta p(T_1)\|^2+\|\Delta (\nabla\cdot\mathbf{v})(T_1)\|^2 \le N_3-1.
\end{aligned}
\end{equation}
From the local existence in Lemma \ref{LWP},  for some $\hat{T}\in(0,\infty)$   there exists a unique classical solution to \eqref{e49h} on the time interval $[T_1,T_1+\hat{T}]$ satisfying
\begin{equation}\label{ismall3}
\begin{aligned}
&\|p(t)\|^2+\|\mathbf{v}(t)\|^2 \le N_1,\\
&\|\nabla p(t)\|^2+\|\nabla\cdot\mathbf{v}(t)\|^2 \le \kappa,\\
&\|\Delta p(t)\|^2+\|\Delta\mathbf{v}(t)\|^2 \le N_2,\\
&\|\nabla\Delta p(t)\|^2+\|\Delta (\nabla\cdot\mathbf{v})(t)\|^2 \le N_3,
\end{aligned}\qquad\forall\ t\in[T_1,T_1+\hat{T}].
\end{equation}
In view of \eqref{n3xx} and \eqref{ismall3}, we see that
\begin{equation}
\begin{aligned}
&\|p(t)\|^2+\|\mathbf{v}(t)\|^2 \le N_1,\\
&\|\nabla p(t)\|^2+\|\nabla\cdot\mathbf{v}(t)\|^2 \le \kappa,\\
&\|\Delta p(t)\|^2+\|\Delta\mathbf{v}(t)\|^2 \le N_2,\\
&\|\nabla\Delta p(t)\|^2+\|\Delta (\nabla\cdot\mathbf{v})(t)\|^2 \le N_3.
\end{aligned}\qquad\forall\ t\in[0,T_1+\hat{T}].
\end{equation}
%By repeating the arguments in Sections \ref{3dl2}--\ref{3dh3}, we know that the solution satisfies the same estimates as in \eqref{ismall1} for $\forall\ t\in[0,T_1+\hat{T}]$, and in particular, the estimates in \eqref{ismall2} are valid when $T_1$ is replaced by $T_1+\hat{T}$. Hence, starting from $T_1+\hat{T}$, by applying the local well-posedness again, we know that the solution can be extended by $\hat{T}$ again, such that it satisfies the same uniform energy estimates within the time interval $[0,T_1+2\hat{T}]$. By repeating the procedure, we know that the solution indeed exists globally in time. Moreover, by repeating the arguments in Section \ref{3dv}, we know that the solution satisfies the same energy estimates derived therein. This completes the proof of the global well-posedness result recorded in Theorem \ref{thm1}.
From the standard procedure as in \cite{RWWZZ}, we conclude that the solution indeed exists globally in time.
This completes the proof of the global well-posedness result recorded in Theorem \ref{thm1}.
\hfill$\square$

\subsection{Long-time Behavior}

In this section, we derive the long-time behavior of the solution obtained above. We shall combine the energy estimates obtained in the previous subsections with the fact that any function $f(t)$, belonging to $W^{1,1}(0, \infty)$, converges to zero as $t\to\infty$, to establish the decay estimate stated in Theorem \ref{thm1}. For brevity, we only present the proof for the decay of the first order spatial derivatives of the solution, in order to illustrate the main idea. The proof for the second and third order derivatives is in a similar fashion and we omit the technical details.

%\begin{lemma}\label{decay}
%Let $f\in W^{1,1}(0, \infty)$ be a nonnegative function. Then $f(t)\to 0$ as $t\to \infty$.
%\end{lemma}

%\begin{proof}
%Since $f\in W^{1,1}(0, \infty)$, $f$ can be identified with a uniformly continuous function a.e., so we may assume that $f$ is uniformly continuous. We prove the lemma by the definition of limit. For any $\epsilon>0$, the uniform continuity of $f$ implies that there exists $\delta>0$ such that, for any $0\le s<t<\infty$,
%$$
%|f(t) -f(s)| <\frac{\epsilon}{2} \quad\mbox{whenever}\quad |t-s|<\delta.
%$$
%Since $f\in L^1(0, \infty)$, there exists $T>0$ such that
%$$
%\int_T^\infty f(t) dt < \frac{\epsilon \delta}{2}.
%$$
%Especially, for any $k=0,1,2,\cdots$,
%$$
%\int_{T+ k\delta}^{T+ (k+1)\delta} f(t) dt < \frac{\epsilon \delta}{2},
%$$
%which implies that there exists $t_k \in (T+ k\delta, T+ (k+1)\delta)$ such that
%$$
%f(t_{k}) < \frac{\epsilon}{2}.
%$$
%Then, for any $t>T$, there exists $k_0$ satisfying $t\in (T+ k_0 \delta, T+ (k_0 +1)\delta]$ and thus
%$$
%f(t) \le f(t_{k_0}) + |f(t) -f(t_{k_0})| \le \frac{\epsilon}{2} + \frac{\epsilon}{2} =\epsilon.
%$$
%This concludes the proof of the lemma.
%\end{proof}

We note that according to \eqref{n9m} and \eqref{nn8h}, it holds that
{\normalsize
\begin{equation}\label{ltb1}
\|\nabla p(t)\|^2+\|\nabla\cdot \mathbf{v}(t)\|^2 \in L^1(0,\infty),
\end{equation}}which is valid for any {\normalsize$\varepsilon\ge0$}.

Next, by using similar arguments as in deriving \eqref{n10h}, we can show that
$$
\begin{aligned}
\frac12\left|\frac{d}{dt}\left(\|\nabla p\|^2+\bar{p}\,\|\nabla\cdot \mathbf{v}\|^2\right)\right| \le&\ \left(\frac14+c_1c_2\,(N_1\kappa)^{\frac14}\right)\|\Delta p\|^2+\frac{27}{4}c_3^4\,\kappa^2\|\nabla p\|^2+ \\
&\quad c_1c_2\,\varepsilon\,\bar{p}\,(N_1\kappa)^{\frac14}\|\Delta \mathbf{v}\|^2+\|\Delta p\|^2+\varepsilon\,\bar{p}\,\|\Delta\mathbf{v}\|^2,
\end{aligned}
$$
which, together with \eqref{kappa1} and \eqref{kappa2}, implies
{\normalsize\begin{equation}\label{ltb2}
\begin{aligned}
&\ \left|\frac{d}{dt}\left(\|\nabla p\|^2+\bar{p}\,\|\nabla\cdot \mathbf{v}\|^2\right)\right| \le 3\|\Delta p\|^2+\frac{27}{2}c_3^4\,\kappa^2\|\nabla p\|^2+3\varepsilon\,\bar{p}\,\|\Delta\mathbf{v}\|^2.
\end{aligned}
\end{equation}}By integrating \eqref{ltb2} with respect to {\normalsize$t$} and applying \eqref{n9m} and \eqref{n13h}, we have
{\normalsize
\begin{equation*}\label{ltb3}
\frac{d}{dt}\left(\|(\nabla p)(t)\|^2+\bar{p}\,\|(\nabla\cdot \mathbf{v})(t)\|^2\right) \in L^1(0,\infty).
\end{equation*}}By combining \eqref{ltb1} and \eqref{ltb2}, we conclude that
{\normalsize
\begin{equation*}\label{ltb4}
\|(\nabla p)(t)\|^2+\bar{p}\,\|(\nabla\cdot \mathbf{v})(t)\|^2 \in W^{1,1}(0,\infty),
\end{equation*}}which implies
{\normalsize$$
\lim_{t\to\infty} \left(\|(\nabla p)(t)\|^2+\bar{p}\,\|(\nabla\cdot \mathbf{v})(t)\|^2\right)=0.
$$}Furthermore, by working with \eqref{n14h} and \eqref{H31} we can obtain the similar result for the second and third order spatial derivatives of the solution. Next, we turn to the zero chemical diffusion limit of the solution and identify the convergence rate in terms of $\varepsilon$, in order to complete the proof for Theorem \ref{thm1}.

\subsection{Diffusion Limit}

Now we study the zero diffusion limit of the solution and quantify the convergence rate in terms of $\varepsilon$. Let {\normalsize$(p^\varepsilon,\mathbf{v}^\varepsilon)$} and {\normalsize$(p^0,\mathbf{v}^0)$} be the solutions to \eqref{e49h} with {\normalsize$\varepsilon>0$} and {\normalsize$\varepsilon=0$}, respectively, for the same initial data, and set {\normalsize$\tilde{p}=p^\varepsilon-p^0$} and {\normalsize$\tilde{\mathbf{v}}=\mathbf{v}^\varepsilon-\mathbf{v}^0$}. Then {\normalsize$(\tilde{p},\tilde{\mathbf{v}})$} satisfies
{\normalsize\begin{equation}\label{dl1}
\left\{
\begin{aligned}
&\partial_t\tilde{p}-\nabla\cdot\tilde{\mathbf{v}}=\Delta\tilde{p}+\nabla\cdot(\tilde{p}\mathbf{v}^\varepsilon+p^0\tilde{\mathbf{v}}),\\
&\partial_t\tilde{\mathbf{v}}-\nabla\tilde{p}=\varepsilon\Delta\mathbf{v}^\varepsilon-\varepsilon\nabla\left(|\mathbf{v}^\varepsilon|^2\right);\\
&(\tilde{p}_0,\tilde{\mathbf{v}}_0)=(0,\mathbf{0}),
\end{aligned}
\right.
\end{equation}}where for simplicity, we took {\normalsize$\bar{p}=1$}. We begin with the zeroth frequency estimate.

{\bf Step 1.} By taking the {\normalsize$L^2$} inner products, we find
{\normalsize\begin{equation}\label{dl2}
\begin{aligned}
&\ \frac12\frac{d}{dt}\left(\|\tilde{p}\|^2+\|\tilde{\mathbf{v}}\|^2\right)+\|\nabla\tilde{p}\|^2\\
=&\ -\int_{\mathbb{R}^3} (\tilde{p}\mathbf{v}^\varepsilon+p^0\tilde{\mathbf{v}})\cdot\nabla\tilde{p}\,d\mathbf{x}+\int_{\mathbb{R}^3}\left[\varepsilon\Delta\mathbf{v}^\varepsilon-\varepsilon\nabla\left(|\mathbf{v}^\varepsilon|^2\right)\right]\cdot\tilde{\mathbf{v}}\,d\mathbf{x}.
\end{aligned}
\end{equation}}For the first term on the RHS of \eqref{dl2}, by applying \eqref{gn3}, we have
{\normalsize\begin{equation}\label{dl3}
\begin{aligned}
&\ \left|-\int_{\mathbb{R}^3} (\tilde{p}\mathbf{v}^\varepsilon+p^0\tilde{\mathbf{v}})\cdot\nabla\tilde{p}\,d\mathbf{x}\right| \\
\le &\ \frac12\,\|\nabla\tilde{p}\|^2+\|\mathbf{v}^\varepsilon\|_{L^\infty}^2\|\tilde{p}\|^2+\|p^0\|_{L^\infty}^2\|\tilde{\mathbf{v}}\|^2\\
\le &\ \frac12\,\|\nabla\tilde{p}\|^2+c_3^2\left(\|\nabla\cdot\mathbf{v}^\varepsilon\|\cdot \|\Delta\mathbf{v}^\varepsilon\| \cdot \|\tilde{p}\|^2+\|\nabla p^0\|\cdot\|\Delta p^0\| \cdot \|\tilde{\mathbf{v}}\|^2\right)\\
\le &\ \frac12\,\|\nabla\tilde{p}\|^2+\frac{c_3^2}{2}\left[(\|\nabla\cdot\mathbf{v}^\varepsilon\|^2+ \|\Delta\mathbf{v}^\varepsilon\|^2)\|\tilde{p}\|^2+(\|\nabla p^0\|^2+\|\Delta p^0\|^2)\|\tilde{\mathbf{v}}\|^2\right].
\end{aligned}
\end{equation}}For the second term on the RHS of \eqref{dl2}, again by applying \eqref{gn3}, we have
{\normalsize\begin{equation}\label{dl4}
\begin{aligned}
&\ \left|\int_{\mathbb{R}^3}\left[\varepsilon\Delta\mathbf{v}^\varepsilon-\varepsilon\nabla\left(|\mathbf{v}^\varepsilon|^2\right)\right]\cdot\tilde{\mathbf{v}}\,d\mathbf{x}\right| \\
\le &\ \frac12\,\|\tilde{\mathbf{v}}\|^2+\varepsilon^2\,\|\Delta\mathbf{v}^\varepsilon\|^2+4\,\varepsilon^2\,\|\mathbf{v}^\varepsilon\|_{L^\infty}^2\|\nabla\cdot\mathbf{v}^\varepsilon\|^2\\
\le &\ \frac12\,\|\tilde{\mathbf{v}}\|^2+\varepsilon^2\,\|\Delta\mathbf{v}^\varepsilon\|^2+4c_3^2\,\varepsilon^2\,\|\nabla\cdot\mathbf{v}^\varepsilon\|\cdot \|\Delta\mathbf{v}^\varepsilon\| \cdot \|\nabla\cdot\mathbf{v}^\varepsilon\|^2\\
\le &\ \frac12\,\|\tilde{\mathbf{v}}\|^2+\varepsilon^2\,\|\Delta\mathbf{v}^\varepsilon\|^2+2c_3^2\sqrt{3\kappa(N_2-1)}\,\varepsilon^2\|\nabla\cdot\mathbf{v}^\varepsilon\|^2.
\end{aligned}
\end{equation}}where we applied \eqref{h2} and \eqref{h3}. By substituting \eqref{dl3} and \eqref{dl4} into \eqref{dl2} then multiplying through by 2, we have
{\normalsize\begin{equation}\label{dl5}
\begin{aligned}
&\ \frac{d}{dt}\left(\|\tilde{p}\|^2+\|\tilde{\mathbf{v}}\|^2\right)+\|\nabla\tilde{p}\|^2\\[2mm]
\le &\ \left[ c_3^2\left(\|\nabla\cdot\mathbf{v}^\varepsilon\|^2+ \|\Delta\mathbf{v}^\varepsilon\|^2+\|\nabla p^0\|^2+\|\Delta p^0\|^2\right)+1\right] \left(\|\tilde{p}\|^2+\|\tilde{\mathbf{v}}\|^2\right)+\\[2mm]
&\ 2\varepsilon^2\,\|\Delta\mathbf{v}^\varepsilon\|^2+4c_3^2\sqrt{3\kappa(N_2-1)}\,\varepsilon^2\|\nabla\cdot\mathbf{v}^\varepsilon\|^2.
\end{aligned}
\end{equation}}
By applying the Gronwall inequality to \eqref{dl5}, we have
{\normalsize\begin{equation*}\label{dl6}
\begin{aligned}
&\ \|\tilde{p}(t)\|^2+\|\tilde{\mathbf{v}}(t)\|^2 \\[2mm]
\le &\ \exp\left\{c_3^2\int_0^t \left(\|\nabla\cdot\mathbf{v}^\varepsilon(\tau)\|^2+ \|\Delta\mathbf{v}^\varepsilon(\tau)\|^2+\|\nabla p^0(\tau)\|^2+\|\Delta p^0(\tau)\|^2\right) d\tau+t\right\} \times\\[2mm]
&\ \left(2\int_0^t \|\Delta\mathbf{v}^\varepsilon(\tau)\|^2 d\tau+ 4c_3^2\sqrt{3\kappa(N_2-1)}\int_0^t\|\nabla\cdot\mathbf{v}^\varepsilon(\tau)\|^2d\tau\right)\varepsilon^2.
\end{aligned}
\end{equation*}}We note that according to the energy estimates recorded in Theorem \ref{thm1}, see \eqref{ub}, the integrals on the RHS of the above inequality are bounded by some constants that are independent of $t$ and remain bounded as $\varepsilon\to0$. We rewrite the above estimate by using short notations as
\begin{equation}\label{dl6}
\begin{aligned}
\|\tilde{p}(t)\|^2+\|\tilde{\mathbf{v}}(t)\|^2 \le e^{C_1+t}\, C_2\,\varepsilon^2.
\end{aligned}
\end{equation}
 Next, we consider the convergence of the first order derivatives of the perturbation.

{\bf Step 2.} By taking the {\normalsize$L^2$} inner products of the first two equations in \eqref{dl1} with the {\normalsize$-\Delta$} of the targeting functions, we deduce
{\normalsize\begin{equation}\label{dl7}
\begin{aligned}
&\ \frac12\frac{d}{dt}\left(\|\nabla\tilde{p}\|^2+\|\nabla\cdot\tilde{\mathbf{v}}\|^2\right)+\|\Delta\tilde{p}\|^2\\
=&\ -\int_{\mathbb{R}^3} \left[\nabla\cdot(\tilde{p}\mathbf{v}^\varepsilon+p^0\tilde{\mathbf{v}})\right]\Delta\tilde{p}\,d\mathbf{x}+\varepsilon\int_{\mathbb{R}^3} (\Delta\nabla\cdot\mathbf{v}^\varepsilon)(\nabla\cdot\tilde{\mathbf{v}})d\mathbf{x}-\\
&\ \varepsilon\int_{\mathbb{R}^3} \Delta\left(|\mathbf{v}^\varepsilon|^2\right) (\nabla\cdot\tilde{\mathbf{v}})d\mathbf{x}.
\end{aligned}
\end{equation}}For the first term on the RHS of \eqref{dl7}, we have
{\normalsize\begin{equation}\label{dl8}
\begin{aligned}
\left|-\int_{\mathbb{R}^3} \left[\nabla\cdot(\tilde{p}\mathbf{v}^\varepsilon+p^0\tilde{\mathbf{v}})\right]\Delta\tilde{p}\,d\mathbf{x}\right|\le \frac12\,\|\Delta\tilde{p}\|^2+\frac12\,\|\nabla\cdot(\tilde{p}\mathbf{v}^\varepsilon+p^0\tilde{\mathbf{v}})\|^2,
\end{aligned}
\end{equation}}where the second term on the RHS can be estimated as
{\normalsize\begin{equation*}\label{dl9}
\begin{aligned}
&\quad\frac12\,\|\nabla\cdot(\tilde{p}\mathbf{v}^\varepsilon+p^0\tilde{\mathbf{v}})\|^2\\[2mm]
&\le 2\left(\|\nabla\tilde{p}\cdot\mathbf{v}^\varepsilon\|^2+\|\tilde{p}(\nabla\cdot\mathbf{v}^\varepsilon)\|^2+\|\nabla p^0\cdot\tilde{\mathbf{v}}\|^2+\|p^0(\nabla\cdot\tilde{\mathbf{v}})\|^2\right)\\[2mm]
&\le 2\left(\|\nabla\tilde{p}\|^2\|\mathbf{v}^\varepsilon\|_{L^\infty}^2+\|\tilde{p}\|_{L^6}^2\|\nabla\cdot\mathbf{v}^\varepsilon\|_{L^3}^2+\|\nabla p^0\|_{L^3}^2\|\tilde{\mathbf{v}}\|_{L^6}^2+\|p^0\|_{L^\infty}^2\|\nabla\cdot\tilde{\mathbf{v}}\|^2\right)\\[2mm]
&\le 2\big(c_3^2\|\nabla\tilde{p}\|^2\|\nabla\cdot\mathbf{v}^\varepsilon\|\cdot\|\Delta \mathbf{v}^\varepsilon\|+c_1^2c_2^2\|\nabla\tilde{p}\|^2\|\nabla\cdot\mathbf{v}^\varepsilon\|\cdot\|\Delta \mathbf{v}^\varepsilon\| +\\[2mm]
&\qquad\qquad  c_1^2c_2^2\|\nabla p^0\|\cdot\|\Delta p^0\|\cdot\|\nabla\cdot\tilde{\mathbf{v}}\|^2+c_3^2\|\nabla p^0\| \cdot \|\Delta p^0\|\cdot \|\nabla\cdot\tilde{\mathbf{v}}\|^2\big)\\[2mm]
&\le (c_3^2+c_1^2c_2^2)\left(\|\nabla\cdot\mathbf{v}^\varepsilon\|^2+\|\Delta \mathbf{v}^\varepsilon\|^2 + \|\nabla p^0\|^2+\|\Delta p^0\|^2 \right)\left(\|\nabla\tilde{p}\|^2+\|\nabla\cdot\tilde{\mathbf{v}}\|^2\right).
\end{aligned}
\end{equation*}}So we can update \eqref{dl8} as
{\normalsize\begin{equation}\label{dl10}
\begin{aligned}
&\ \left|-\int_{\mathbb{R}^3} \left[\nabla\cdot(\tilde{p}\mathbf{v}^\varepsilon+p^0\tilde{\mathbf{v}})\right]\Delta\tilde{p}\,d\mathbf{x}\right| \\
\le &\ \frac12\,\|\Delta\tilde{p}\|^2+(c_3^2+c_1^2c_2^2)\left(\|\nabla\cdot\mathbf{v}^\varepsilon\|^2+\|\Delta \mathbf{v}^\varepsilon\|^2 + \|\nabla p^0\|^2+\|\Delta p^0\|^2 \right)\left(\|\nabla\tilde{p}\|^2+\|\nabla\cdot\tilde{\mathbf{v}}\|^2\right).
\end{aligned}
\end{equation}}For the second and third terms on the RHS of \eqref{dl7}, in a similar fashion, we can show that
{\normalsize\begin{equation}\label{dl11}
\begin{aligned}
&\ \left|\varepsilon\int_{\mathbb{R}^3} (\Delta\nabla\cdot\mathbf{v}^\varepsilon)(\nabla\cdot\tilde{\mathbf{v}})d\mathbf{x}-\varepsilon\int_{\mathbb{R}^3} \Delta\left(|\mathbf{v}^\varepsilon|^2\right) (\nabla\cdot\tilde{\mathbf{v}})d\mathbf{x}\right|\\
\le &\ \frac12\,\|\nabla\cdot\tilde{\mathbf{v}}\|^2+\varepsilon^2\,\|\Delta\nabla\cdot\mathbf{v}^\varepsilon\|^2+\varepsilon^2\left\|\Delta\left(|\mathbf{v}^\varepsilon|^2\right) \right\|^2\\
\le &\ \frac12\,\|\nabla\cdot\tilde{\mathbf{v}}\|^2+\varepsilon^2\,\|\Delta\nabla\cdot\mathbf{v}^\varepsilon\|^2+8\,\varepsilon^2\,\left(\|\Delta\mathbf{v}^\varepsilon\|^2\|\mathbf{v}^\varepsilon\|_{L^\infty}^2+\|\nabla\cdot\mathbf{v}^\varepsilon\|_{L^3}^2\|\nabla\mathbf{v}^\varepsilon\|_{L^6}^2\right)\\
\le &\ \frac12\,\|\nabla\cdot\tilde{\mathbf{v}}\|^2+\varepsilon^2\,\|\Delta\nabla\cdot\mathbf{v}^\varepsilon\|^2+8\left(c_3^2+c_1^2c_2^2\right)\,\varepsilon^2\,\|\nabla\cdot \mathbf{v}^\varepsilon\| \cdot\|\Delta\mathbf{v}^\varepsilon\| \cdot \|\Delta\mathbf{v}^\varepsilon\|^2\\
\le &\ \frac12\,\|\nabla\cdot\tilde{\mathbf{v}}\|^2+\varepsilon^2\|\Delta\nabla\cdot\mathbf{v}^\varepsilon\|^2+4\left(c_3^2+c_1^2c_2^2\right)\sqrt{3\kappa (N_2-1)}\,\|\Delta\mathbf{v}^\varepsilon\|^2\,\varepsilon^2.
\end{aligned}
\end{equation}}By feeding \eqref{dl10} and \eqref{dl11} into \eqref{dl7}, we find
{\normalsize\begin{equation}\label{dl12}
\begin{aligned}
&\ \frac{d}{dt}\left(\|\nabla\tilde{p}\|^2+\|\nabla\cdot\tilde{\mathbf{v}}\|^2\right)+\|\Delta\tilde{p}\|^2 \\[2mm]
\le &\ \left[2(c_3^2+c_1^2c_2^2)\left(\|\nabla\cdot\mathbf{v}^\varepsilon\|^2+\|\Delta \mathbf{v}^\varepsilon\|^2 + \|\nabla p^0\|^2+\|\Delta p^0\|^2\right)+1\right]\left(\|\nabla\tilde{p}\|^2+\|\nabla\cdot\tilde{\mathbf{v}}\|^2\right)+\\[2mm]
&\quad 2\varepsilon^2\|\Delta\nabla\cdot\mathbf{v}^\varepsilon\|^2+8\left(c_3^2+c_1^2c_2^2\right)\sqrt{3\kappa (N_2-1)}\,\|\Delta\mathbf{v}^\varepsilon\|^2\,\varepsilon^2.
\end{aligned}
\end{equation}}By applying the Gronwall inequality to \eqref{dl12}, we deduce
{\normalsize$$
\begin{aligned}
&\ \|\nabla\tilde{p}(t)\|^2+\|\nabla\cdot\tilde{\mathbf{v}}(t)\|^2\\[2mm]
\le &\ \exp\left\{2(c_3^2+c_1^2c_2^2)\int_0^t \left(\|\nabla\cdot\mathbf{v}^\varepsilon(\tau)\|^2+ \|\Delta\mathbf{v}^\varepsilon(\tau)\|^2+\|\nabla p^0(\tau)\|^2+\|\Delta p^0(\tau)\|^2\right) d\tau+t\right\} \times\\[2mm]
&\ \left(2\int_0^t \|\Delta \nabla\cdot\mathbf{v}^\varepsilon(\tau)\|^2 d\tau+ 8\left(c_3^2+c_1^2c_2^2\right)\sqrt{3\kappa (N_2-1)}\int_0^t\|\Delta\mathbf{v}^\varepsilon(\tau)\|^2d\tau\right)\varepsilon^2.
\end{aligned}
$$}We note again that the integrals on the RHS of the above inequality are bounded by some constants that are independent of $t$ and remain bounded as $\varepsilon\to0$. We rewrite the above estimate as
\begin{equation}\label{dl13}
\begin{aligned}
\|\nabla\tilde{p}(t)\|^2+\|\nabla\cdot\tilde{\mathbf{v}}(t)\|^2 \le e^{C_3+t}\, C_4\,\varepsilon^2.
\end{aligned}
\end{equation}

{\bf Step 3.} The convergence of the second order derivatives of the solution can be established in a similar fashion, except that the term $\varepsilon^2\|\Delta^2\mathbf{v}^\varepsilon(t)\|^2$ will appear in the energy estimates. According to \eqref{ub}, the temporal integral of such term, i.e.,
$$
\int_0^t\varepsilon^2\|\Delta^2\mathbf{v}^\varepsilon(\tau)\|^2d\tau,
$$
has the order of $O(\varepsilon)$, instead of $O(\varepsilon^2)$, as $\varepsilon\to 0$. The final energy estimate takes the form
$$
\|\Delta\tilde{p}(t)\|^2+\|\Delta\tilde{\mathbf{v}}(t)\|^2 \le e^{C_5+t}\, C_6\,(1+\varepsilon)\varepsilon.
$$
We omit the technical details to simplify the presentation. This completes the proof for Theorem \ref{thm1}. \hfill$\square$

%%%%%%%%%%%%%%%%%%%%%%%%%%%%%%%%%%%%%
%%%%%%%%%%%%%%%%%%%%%%%%%%%%%%%%%%%%%

\section{Proof of Theorem \ref{thm2}}

In this section we prove Theorem \ref{thm2}. Comparing with the 3D case, the proof in this section is much lengthier, due to (as was mentioned in the Introduction) the Gagliardo-Nirenberg interpolation inequalities generate less powers of high frequencies of a function in {$\mathbb{R}^2$} than in {$\mathbb{R}^3$}. Such a deficiency has a substantial impact on the energy estimates for all the individual frequencies of the solution (cf. \eqref{n7h}, \eqref{n10h}, {\it et al}), especially when the solution carries potentially large $L^2$-norm of the zeroth frequency of the perturbation. We overcome the difficulty by creating higher order nonlinearities, through coupling and cancellation, to compensate the deficiency.

Again, for the reader's convenience, we re-state the Cauchy problem:
\begin{equation}\label{n1}
\left\{
\begin{aligned}
&\partial_tp-\nabla\cdot(p\mathbf{v})-\bar{p}\nabla\cdot\mathbf{v}=\Delta p,\\
&\partial_t\mathbf{v}-\nabla \left(p-\varepsilon|\mathbf{v}|^2\right)=\varepsilon\Delta\mathbf{v},\\
&(p_0,\mathbf{v}_0)\in H^3(\mathbb{R}^2),\quad p_0+\bar{p}>0,\quad \nabla\times\mathbf{v}_0=\mathbf{0},
\end{aligned}
\right.
\end{equation}
where $(p,\mathbf{v})$ denotes the perturbation of the original solution to \eqref{hpbl} around the constant state $(\bar{p},\mathbf{0})$.
%\begin{notation}
%In this section, we use {\normalsize$d_i\,(i=1,2,...)$} to denote a generic constant which is independent of the unknown functions, {\normalsize$\varepsilon$}, {\normalsize$t$} and initial data, while {\normalsize$D_i\,(i=1,2,...)$} denotes a generic constant which is independent of the unknown functions, {\normalsize$\varepsilon$} and {\normalsize$t$}, but depends on the initial data. The values of the constants may vary line by line according to the context.
%\end{notation}
As in the three-dimensional case, the local well-posedness of classical solutions to \eqref{n1} follows as an application of Kawashima's theory \cite{Kawashima}, and the positivity of the function $(p+\bar{p})(\mathbf{x},t)$ follows from the maximum principle (cf. \cite{Friedman}).

Now we recall several Gagliardo-Nirenberg interpolation inequalities in 2D:
\begin{align}
&\|f\|_{L^4}\le d_1\,\|\nabla f\|^{\frac12}\|f\|^{\frac12}, \label{gn2d1}\\
&\|f\|_{L^6}\le d_2\,\|\nabla f\|^{\frac23}\|f\|^{\frac13}, \label{gn2d2}\\
&\|f\|_{L^8}\le d_3\,\|\nabla f\|^{\frac34}\|f\|^{\frac14}, \label{gn2d3}\\
&\|f\|_{L^{12}}\le d_4\,\|\nabla f\|^{\frac56}\|f\|^{\frac16}, \label{gn2d4}\\
&\|f\|_{L^\infty}\le d_5\,\|f\|^{\frac12}\|\Delta f\|^{\frac12}.\label{gn2d5}
\end{align}
In what follows, we assume that for a local existence time {\normalsize$T>0$} the following hold:
{\normalsize
\begin{equation}\label{n2}
\begin{aligned}
&\sup_{0\le t\le T}\left(\|p(t)\|^2_{L^2}+\|\mathbf{v}(t)\|^2_{L^2}\right)\le M_1,\\
&\sup_{0\le t\le T}\left(\|\nabla p(t)\|_{L^2}^2+\|\nabla\cdot\mathbf{v}(t)\|_{L^2}^2\right)\le \delta,
\end{aligned}
\end{equation}}where {\normalsize$M_1$} and $\delta$ are defined in Theorem \ref{thm2}. We being with the estimate of the $L^2$-norm of the zeroth frequency part of the solution.

\subsection{$L^2$-estimate}
By testing the equations in \eqref{n1} with the targeting functions, we have
{\normalsize
\begin{equation}\label{n3}
\begin{aligned}
&\ \frac12\frac{d}{dt}\left(\|p\|^2+\bar{p}\,\|\mathbf{v}\|^2\right)+\|\nabla p\|^2+\varepsilon\,\bar{p}\,\|\nabla\cdot\mathbf{v}\|^2\\
=&\ -\int_{\mathbb{R}^2}p\mathbf{v}\cdot\nabla p\ d\mathbf{x}+\varepsilon\,\bar{p}\int_{\mathbb{R}^2}|\mathbf{v}|^2\nabla\cdot\mathbf{v}\ d\mathbf{x}.
\end{aligned}
\end{equation}}We note that the two terms on the RHS of \eqref{n3} can not be estimated as in the three-dimensional case, as the interpolation inequalities in 2D does not generate enough powers of the higher frequencies of a function. Next, we eliminate those two terms through performing higher order energy estimates. During such a process, the higher order nonlinearities can be controlled by using the smallness of the {\normalsize$L^2$} norm of the first order derivatives of the solution. We divide the subsequent proof into six steps.

{\bf Step 1.} Taking the {\normalsize$L^2$} inner product of the first equation in \eqref{n1} with {\normalsize$-\varepsilon\,|\mathbf{v}|^2$}, we have
{\normalsize\begin{equation}\label{n3a}
\begin{aligned}
&\ -\int_{\mathbb{R}^2}\varepsilon\,|\mathbf{v}|^2\partial_t p\ d\mathbf{x}\\
=&\ -\varepsilon\int_{\mathbb{R}^2}|\mathbf{v}|^2\nabla\cdot(p\mathbf{v})\ d\mathbf{x}-\varepsilon\,\,\bar{p}\int_{\mathbb{R}^2}|\mathbf{v}|^2\nabla\cdot\mathbf{v}\ d\mathbf{x}-\varepsilon\int_{\mathbb{R}^2}|\mathbf{v}|^2\Delta p\ d\mathbf{x}.
\end{aligned}
\end{equation}}Taking the {\normalsize$L^2$} inner product of the second equation in \eqref{n1} with {\normalsize$-2\,\varepsilon\, p\,\mathbf{v}$}, we have
{\normalsize
\begin{equation}\label{n3b}
\begin{aligned}
&\ -\int_{\mathbb{R}^2}\varepsilon\, p\,\partial_t(|\mathbf{v}|^2)\ d\mathbf{x} \\
=&\ -2\,\varepsilon\int_{\mathbb{R}^2} p\mathbf{v}\cdot\nabla p\ d\mathbf{x}-2\,\varepsilon^2\int_{\mathbb{R}^2}p\mathbf{v}\cdot\Delta\mathbf{v}\ d\mathbf{x}+2\,\varepsilon^2\int_{\mathbb{R}^2}p\mathbf{v}\cdot\nabla(|\mathbf{v}|^2)\ d\mathbf{x}.
\end{aligned}
\end{equation}}Adding \eqref{n3a} and \eqref{n3b}, we find
{\normalsize
\begin{equation}\label{n3c}
\begin{aligned}
&\ -\frac{d}{dt}\left(\varepsilon\int_{\mathbb{R}^2}p|\mathbf{v}|^2d\mathbf{x}\right)\\
=&\ -\varepsilon\int_{\mathbb{R}^2}|\mathbf{v}|^2\nabla\cdot(p\mathbf{v})\ d\mathbf{x}-\varepsilon\,\bar{p}\int_{\mathbb{R}^2}|\mathbf{v}|^2\nabla\cdot\mathbf{v}\ d\mathbf{x}-\varepsilon\int_{\mathbb{R}^2}|\mathbf{v}|^2\Delta p\ d\mathbf{x}-\\
&2\,\varepsilon\int_{\mathbb{R}^2} p\mathbf{v}\cdot\nabla p\ d\mathbf{x}-2\,\varepsilon^2\int_{\mathbb{R}^2}p\mathbf{v}\cdot\Delta\mathbf{v}\ d\mathbf{x}+2\,\varepsilon^2\int_{\mathbb{R}^2}p\mathbf{v}\cdot\nabla(|\mathbf{v}|^2)\ d\mathbf{x}.
\end{aligned}
\end{equation}}We observe that the second term on the RHS of \eqref{n3c} terminates the second term on the RHS of \eqref{n3} upon addition. Adding \eqref{n3c} to \eqref{n3}, we have
{\normalsize
\begin{equation}\label{n3d}
\begin{aligned}
&\ \frac{d}{dt}\left(\frac12\,\|p\|^2+\frac{\bar{p}}{2}\,\|\mathbf{v}\|^2-\varepsilon\int_{\mathbb{R}^2}p|\mathbf{v}|^2d\mathbf{x}\right)+\|\nabla p\|^2+\varepsilon\,\bar{p}\|\nabla\cdot\mathbf{v}\|^2\\
=&\ -(2\,\varepsilon+1)\int_{\mathbb{R}^2}p\mathbf{v}\cdot\nabla p\ d\mathbf{x}+\varepsilon\,(2\,\varepsilon+1)\int_{\mathbb{R}^2}p\mathbf{v}\cdot\nabla(|\mathbf{v}|^2) \ d\mathbf{x}-\\
&\ \varepsilon\int_{\mathbb{R}^2}|\mathbf{v}|^2\Delta p\ d\mathbf{x}-2\,\varepsilon^2\int_{\mathbb{R}^2}p\mathbf{v}\cdot\Delta\mathbf{v}\ d\mathbf{x}.
\end{aligned}
\end{equation}}Taking the {\normalsize$L^2$} inner product of the first equation in \eqref{n1} with {\normalsize$-(2\,\varepsilon+1)\,p^2$}, we have
{\normalsize
\begin{equation}\label{n4}
\begin{aligned}
&\ -\frac{d}{dt}\left(\frac{2\,\varepsilon+1}{6}\int_{\mathbb{R}^2}p^3d\mathbf{x}\right)-(2\,\varepsilon+1)\int_{\mathbb{R}^2}p|\nabla p|^2d\mathbf{x}\\
=&\ (2\,\varepsilon+1)\int_{\mathbb{R}^2}p\mathbf{v}\cdot\nabla p\ d\mathbf{x}+(2\,\varepsilon+1)\,\bar{p}\int_{\mathbb{R}^2}p^2\mathbf{v}\cdot\nabla p\ d\mathbf{x},
\end{aligned}
\end{equation}}where the first term on the RHS terminates the first term on the RHS of \eqref{n3d} upon addition. Adding \eqref{n4} to \eqref{n3d}, we find
{\normalsize
\begin{equation}\label{n3da}
\begin{aligned}
&\ \frac{d}{dt}\left(\frac12\,\|p\|^2+\frac{\bar{p}}{2}\,\|\mathbf{v}\|^2-\varepsilon\int_{\mathbb{R}^2}p|\mathbf{v}|^2d\mathbf{x}-\frac{2\,\varepsilon+1}{6}\int_{\mathbb{R}^2}p^3d\mathbf{x}\right)+\\
&\quad\quad\quad\|\nabla p\|^2+\varepsilon\,\bar{p}\,\|\nabla\cdot\mathbf{v}\|^2-(2\,\varepsilon+1)\int_{\mathbb{R}^2}p|\nabla p|^2d\mathbf{x}\\
=&\ (2\,\varepsilon+1)\,\bar{p}\int_{\mathbb{R}^2}p^2\mathbf{v}\cdot\nabla p\ d\mathbf{x}+\varepsilon\,(2\,\varepsilon+1)\int_{\mathbb{R}^2}p\mathbf{v}\cdot\nabla(|\mathbf{v}|^2) \ d\mathbf{x}-\\
&\ \varepsilon\int_{\mathbb{R}^2}|\mathbf{v}|^2\Delta p\ d\mathbf{x}-2\,\varepsilon^2\int_{\mathbb{R}^2}p\mathbf{v}\cdot\Delta\mathbf{v}\ d\mathbf{x}.
\end{aligned}
\end{equation}}We note that the expression inside the parenthesis on the LHS of \eqref{n3da} is not necessarily positive. Hence, we need to supply terms in order to gain the positivity of the quantity.

{\bf Step 2.} First, for any positive constant {\normalsize$k_1$}, taking the {\normalsize$L^2$} inner product of the first equation in \eqref{n1} with {\normalsize$k_1p^3$}, we have
{\normalsize
\begin{equation}\label{n5}
\begin{aligned}
&\ \frac{d}{dt}\left(k_1\int_{\mathbb{R}^2}p^4d\mathbf{x}\right)+12\,k_1\int_{\mathbb{R}^2}p^2|\nabla p|^2d\mathbf{x}\\
=&\ -12\,k_1\,\bar{p}\int_{\mathbb{R}^2}p^2\mathbf{v}\cdot\nabla p\ d\mathbf{x}-12\,k_1\int_{\mathbb{R}^2}p^3\mathbf{v}\cdot\nabla p\ d\mathbf{x}.
\end{aligned}
\end{equation}}Adding \eqref{n5} to \eqref{n3da}, we have
{\normalsize
\begin{equation}\label{n6}
\begin{aligned}
&\ \frac{d}{dt}E_1(t)+D_1(t) \\
=&\ (2\,\varepsilon+1-12\,k_1)\,\bar{p}\int_{\mathbb{R}^2}p^2\mathbf{v}\cdot\nabla p\ d\mathbf{x}+\varepsilon\,(2\,\varepsilon+1)\int_{\mathbb{R}^2}p\mathbf{v}\cdot\nabla(|\mathbf{v}|^2) \ d\mathbf{x}-\\
&\varepsilon\int_{\mathbb{R}^2}|\mathbf{v}|^2\Delta p\ d\mathbf{x}-2\,\varepsilon^2\int_{\mathbb{R}^2}p\mathbf{v}\cdot\Delta\mathbf{v}\ d\mathbf{x}-12\,k_1\int_{\mathbb{R}^2}p^3\mathbf{v}\cdot\nabla p\ d\mathbf{x},
\end{aligned}
\end{equation}}where
{\normalsize$$
\begin{aligned}
E_1(t)&\equiv \frac12\,\|p\|^2+\frac{\bar{p}}{2}\,\|\mathbf{v}\|^2-\varepsilon\int_{\mathbb{R}^2}p|\mathbf{v}|^2d\mathbf{x}-\frac{2\,\varepsilon+1}{6}\int_{\mathbb{R}^2}p^3d\mathbf{x}+k_1\int_{\mathbb{R}^2}p^4d\mathbf{x},\\
D_1(t)&\equiv \|\nabla p\|^2+\varepsilon\,\bar{p}\,\|\nabla\cdot\mathbf{v}\|^2-(2\,\varepsilon+1)\int_{\mathbb{R}^2}p|\nabla p|^2d\mathbf{x}+12\,k_1\int_{\mathbb{R}^2}p^2|\nabla p|^2d\mathbf{x}.
\end{aligned}
$$}Second, by taking the {\normalsize$L^2$} inner product of the first equation in \eqref{n1} with {\normalsize$2\,\varepsilon^2|\mathbf{v}|^2p$}, we have
{\normalsize
\begin{equation}\label{n6a}
\begin{aligned}
&\ \int_{\mathbb{R}^2}\varepsilon^2|\mathbf{v}|^2\partial_t(p^2)\\
=&\ 2\,\varepsilon^2\int_{\mathbb{R}^2}|\mathbf{v}|^2p \nabla\cdot(p\mathbf{v})\ d\mathbf{x}+2\,\varepsilon^2\,\bar{p}\int_{\mathbb{R}^2}|\mathbf{v}|^2p \nabla\cdot\mathbf{v}\ d\mathbf{x}
+2\,\varepsilon^2\int_{\mathbb{R}^2}|\mathbf{v}|^2p \Delta p\ d\mathbf{x}.
\end{aligned}
\end{equation}}Taking the {\normalsize$L^2$} inner product of the second equation in \eqref{n1} with {\normalsize$2\,\varepsilon^2p^2\mathbf{v}$}, we have
{\normalsize
\begin{equation}\label{n6b}
\begin{aligned}
&\ \int_{\mathbb{R}^2}\varepsilon^2p^2\partial_t(|\mathbf{v}|^2)\\
=&\ 2\,\varepsilon^2\int_{\mathbb{R}^2}p^2\mathbf{v}\cdot \nabla p\ d\mathbf{x}+2\,\varepsilon^3\int_{\mathbb{R}^2}p^2\mathbf{v}\cdot\Delta\mathbf{v}\ d\mathbf{x}
-2\,\varepsilon^3\int_{\mathbb{R}^2}p^2\mathbf{v}\cdot\nabla(|\mathbf{v}|^2)\ d\mathbf{x}.
\end{aligned}
\end{equation}}Adding \eqref{n6a} and \eqref{n6b}, we find
{\normalsize
\begin{equation}\label{n6c}
\begin{aligned}
&\ \frac{d}{dt}\left(\int_{\mathbb{R}^2}\varepsilon^2p^2|\mathbf{v}|^2d\mathbf{x}\right)\\
=&\ 2\,\varepsilon^2\int_{\mathbb{R}^2}|\mathbf{v}|^2p \nabla\cdot(p\mathbf{v})\ d\mathbf{x}+2\,\varepsilon^2\,\bar{p}\int_{\mathbb{R}^2}|\mathbf{v}|^2p \nabla\cdot\mathbf{v}\ d\mathbf{x}
+2\,\varepsilon^2\int_{\mathbb{R}^2}|\mathbf{v}|^2p \Delta p\ d\mathbf{x}+\\
&\ 2\,\varepsilon^2\int_{\mathbb{R}^2}p^2\mathbf{v}\cdot \nabla p\ d\mathbf{x}+2\,\varepsilon^3\int_{\mathbb{R}^2}p^2\mathbf{v}\cdot\Delta\mathbf{v}\ d\mathbf{x}
-2\,\varepsilon^3\int_{\mathbb{R}^2}p^2\mathbf{v}\cdot\nabla(|\mathbf{v}|^2)\ d\mathbf{x}.
\end{aligned}
\end{equation}}We note that the third term on the RHS of \eqref{n6c}:
{\normalsize\begin{equation}\label{n6d}
2\,\varepsilon^2\int_{\mathbb{R}^2}|\mathbf{v}|^2p \Delta p\ d\mathbf{x}=-2\,\varepsilon^2\int_{\mathbb{R}^2}|\mathbf{v}|^2|\nabla p|^2 d\mathbf{x}-2\,\varepsilon^2\int_{\mathbb{R}^2}p\nabla(|\mathbf{v}|^2)\cdot\nabla p\ d\mathbf{x},
\end{equation}}and the fifth term:
{\normalsize\begin{equation}\label{n6e}
\begin{aligned}
&\ 2\,\varepsilon^3\int_{\mathbb{R}^2}p^2\mathbf{v}\cdot\Delta\mathbf{v}\ d\mathbf{x}\\
=&\ -2\,\varepsilon^3\int_{\mathbb{R}^2}p^2\left(|\nabla v_1|^2+|\nabla v_2|^2\right)d\mathbf{x}-4\,\varepsilon^3\int_{\mathbb{R}^2}\left(pv_1\nabla p\cdot\nabla v_1+pv_2\nabla p\cdot\nabla v_2\right)d\mathbf{x}.
\end{aligned}
\end{equation}}Plugging \eqref{n6d} and \eqref{n6e} into \eqref{n6c}, then multiplying the result by 2, we have
{\normalsize
\begin{equation}\label{n6f}
\begin{aligned}
&\ \frac{d}{dt}\left(\int_{\mathbb{R}^2}2\,\varepsilon^2\,p^2\,|\mathbf{v}|^2d\mathbf{x}\right)+4\,\varepsilon^2\int_{\mathbb{R}^2}|\mathbf{v}|^2|\nabla p|^2 d\mathbf{x}+4\,\varepsilon^3\int_{\mathbb{R}^2}p^2\left(|\nabla v_1|^2+|\nabla v_2|^2\right)d\mathbf{x}\\
=&\ 4\,\varepsilon^2\int_{\mathbb{R}^2}|\mathbf{v}|^2p \nabla\cdot(p\mathbf{v})\ d\mathbf{x}+4\varepsilon^2\,\bar{p}\int_{\mathbb{R}^2}|\mathbf{v}|^2p \nabla\cdot\mathbf{v}\ d\mathbf{x}-4\,\varepsilon^2\int_{\mathbb{R}^2}p\nabla(|\mathbf{v}|^2)\cdot\nabla p\ d\mathbf{x}+\\
&\ 4\,\varepsilon^2\int_{\mathbb{R}^2}p^2\mathbf{v}\cdot \nabla p\ d\mathbf{x}
-4\,\varepsilon^3\int_{\mathbb{R}^2}p^2\mathbf{v}\cdot\nabla(|\mathbf{v}|^2)\ d\mathbf{x}-\\
&\ 8\,\varepsilon^3\int_{\mathbb{R}^2}\left(pv_1\nabla p\cdot\nabla v_1+pv_2\nabla p\cdot\nabla v_2\right)d\mathbf{x}.
\end{aligned}
\end{equation}}Next, we carry out some preliminary energy estimates for both the RHS of \eqref{n6} and \eqref{n6f}.

{\bf Step 3.} For the third term on the RHS of \eqref{n6}, we have
{\normalsize
\begin{equation}\label{n6i}
\begin{aligned}
\left|\varepsilon\int_{\mathbb{R}^2}|\mathbf{v}|^2\Delta p\ d\mathbf{x}\right|&=\left|\varepsilon\int_{\mathbb{R}^2}\nabla\left(|\mathbf{v}|^2\right)\cdot\nabla p\ d\mathbf{x} \right|\\
&\le \frac{1}{4}\,\|\nabla p\|^2+8\,\varepsilon^2\int_{\mathbb{R}^2}\left((v_1)^2|\nabla v_1|^2+(v_2)^2|\nabla v_2|^2\right)d\mathbf{x}.
\end{aligned}
\end{equation}}For the fourth term on the RHS of \eqref{n6}, we have
{\normalsize
\begin{equation}\label{n6j}
\begin{aligned}
&\ \left|2\,\varepsilon^2\int_{\mathbb{R}^2}p\mathbf{v}\cdot\Delta\mathbf{v}\ d\mathbf{x}\right| \\
=&\ \left|2\,\varepsilon^2\int_{\mathbb{R}^2}p\left(|\nabla v_1|^2+|\nabla v_2|^2\right)d\mathbf{x}+4\,\varepsilon^2\int_{\mathbb{R}^2}\left(v_1\nabla p\cdot\nabla v_1+v_2\nabla p\cdot\nabla v_2\right)d\mathbf{x}\right|\\
\le&\ \frac{\varepsilon\,\bar{p}}{2}\,\|\nabla\cdot\mathbf{v}\|^2+\frac{2\,\varepsilon^3}{\bar{p}}\int_{\mathbb{R}^2}p^2\left(|\nabla v_1|^2+|\nabla v_2|^2\right)d\mathbf{x}+\frac{1}{4}\,\|\nabla p\|^2+\\
&\ 32\,\varepsilon^4\int_{\mathbb{R}^2}\left((v_1)^2|\nabla v_1|^2+(v_2)^2|\nabla v_2|^2\right)d\mathbf{x}.
\end{aligned}
\end{equation}}Similarly, for the third term on the RHS of \eqref{n6f}, we can show that
{\normalsize
\begin{equation}\label{n6k}
\begin{aligned}
&\ \left|4\,\varepsilon^2\int_{\mathbb{R}^2}p\nabla(|\mathbf{v}|^2)\cdot\nabla p\ d\mathbf{x}\right| \\
=&\ \left|8\,\varepsilon^2\int_{\mathbb{R}^2}\left(pv_1\nabla p\cdot\nabla v_1+pv_2\nabla p\cdot\nabla v_2\right)d\mathbf{x}\right|\\
\le&\ 8\,\varepsilon^2\int_{\mathbb{R}^2}p^2|\nabla p|^2d\mathbf{x}+4\,\varepsilon^2\int_{\mathbb{R}^2}\left((v_1)^2|\nabla v_1|^2+(v_2)^2|\nabla v_2|^2\right)d\mathbf{x}.
\end{aligned}
\end{equation}}In the same way, for the sixth term on the RHS of \eqref{n6f}, we have
{\normalsize
\begin{equation}\label{n6l}
\begin{aligned}
&\ \left|8\,\varepsilon^3\int_{\mathbb{R}^2}\left(pv_1\nabla p\cdot\nabla v_1+pv_2\nabla p\cdot\nabla v_2\right)d\mathbf{x}\right|\\
\le&\ 8\,\varepsilon^3\int_{\mathbb{R}^2}p^2|\nabla p|^2d\mathbf{x}+4\,\varepsilon^3\int_{\mathbb{R}^2}\left((v_1)^2|\nabla v_1|^2+(v_2)^2|\nabla v_2|^2\right)d\mathbf{x}.
\end{aligned}
\end{equation}}By feeding \eqref{n6i} and \eqref{n6j} into \eqref{n6}, we have
{\normalsize
\begin{equation}\label{n6m}
\begin{aligned}
&\ \frac{d}{dt}E_1(t)+D_1(t)\\
\le&\ (2\,\varepsilon+1-12\,k_1)\int_{\mathbb{R}^2}p^2\mathbf{v}\cdot\nabla p\ d\mathbf{x}+\varepsilon\,(2\,\varepsilon+1)\int_{\mathbb{R}^2}p\mathbf{v}\cdot\nabla(|\mathbf{v}|^2) \ d\mathbf{x}-\\
&\ 12\,k_1\int_{\mathbb{R}^2}p^3\mathbf{v}\cdot\nabla p\ d\mathbf{x}+\frac{\varepsilon\,\bar{p}}{2}\,\|\nabla\cdot\mathbf{v}\|^2+\frac{2\,\varepsilon^3}{\bar{p}}\int_{\mathbb{R}^2}p^2\left(|\nabla v_1|^2+|\nabla v_2|^2\right)d\mathbf{x}+\\
&\ \frac{1}{2}\,\|\nabla p\|^2+\left(8\,\varepsilon^2+32\,\varepsilon^4\right)\int_{\mathbb{R}^2}\left((v_1)^2|\nabla v_1|^2+(v_2)^2|\nabla v_2|^2\right)d\mathbf{x}.
\end{aligned}
\end{equation}}By feeding \eqref{n6k} and \eqref{n6l} into \eqref{n6f}, we have
{\normalsize
\begin{equation}\label{n6n}
\begin{aligned}
&\ \frac{d}{dt}\left(\int_{\mathbb{R}^2}2\,\varepsilon^2\,p^2\,|\mathbf{v}|^2d\mathbf{x}\right)+4\,\varepsilon^2\int_{\mathbb{R}^2}|\mathbf{v}|^2|\nabla p|^2 d\mathbf{x}+4\,\varepsilon^3\int_{\mathbb{R}^2}p^2\left(|\nabla v_1|^2+|\nabla v_2|^2\right)d\mathbf{x}\\
\le&\ 4\,\varepsilon^2\int_{\mathbb{R}^2}|\mathbf{v}|^2p \nabla\cdot(p\mathbf{v})\ d\mathbf{x}+4\,\varepsilon^2\int_{\mathbb{R}^2}|\mathbf{v}|^2p \nabla\cdot\mathbf{v}\ d\mathbf{x}+4\,\varepsilon^2\int_{\mathbb{R}^2}p^2\mathbf{v}\cdot \nabla p\ d\mathbf{x}-\\
&\ 4\,\varepsilon^3\int_{\mathbb{R}^2}p^2\mathbf{v}\cdot\nabla(|\mathbf{v}|^2)\ d\mathbf{x}+\left(8\,\varepsilon^2+8\,\varepsilon^3\right)\int_{\mathbb{R}^2}p^2|\nabla p|^2d\mathbf{x}+\\
&\ \left(4\,\varepsilon^2+4\,\varepsilon^3\right)\int_{\mathbb{R}^2}\left((v_1)^2|\nabla v_1|^2+(v_2)^2|\nabla v_2|^2\right)d\mathbf{x}.
\end{aligned}
\end{equation}}By adding \eqref{n6n}{\normalsize$\times \frac{1}{\bar{p}}$} to \eqref{n6m}, we find
{\normalsize
\begin{equation}\label{n6o}
\begin{aligned}
&\ \frac{d}{dt}E_2(t)+D_2(t) \\
\le &\ (2\,\varepsilon+1-12\,k_1)\int_{\mathbb{R}^2}p^2\mathbf{v}\cdot\nabla p\ d\mathbf{x}+\varepsilon\,(2\,\varepsilon+1)\int_{\mathbb{R}^2}p\mathbf{v}\cdot\nabla(|\mathbf{v}|^2) \ d\mathbf{x}-\\
&\ 12\,k_1\int_{\mathbb{R}^2}p^3\mathbf{v}\cdot\nabla p\ d\mathbf{x}+\frac{4\,\varepsilon^2}{\bar{p}}\int_{\mathbb{R}^2}|\mathbf{v}|^2p \nabla\cdot(p\mathbf{v})\ d\mathbf{x}+\frac{4\,\varepsilon^2}{\bar{p}}\int_{\mathbb{R}^2}|\mathbf{v}|^2p \nabla\cdot\mathbf{v}\ d\mathbf{x}+\\
&\ \frac{4\,\varepsilon^2}{\bar{p}}\int_{\mathbb{R}^2}p^2\mathbf{v}\cdot \nabla p\ d\mathbf{x}-\frac{4\,\varepsilon^3}{\bar{p}}\int_{\mathbb{R}^2}p^2\mathbf{v}\cdot\nabla(|\mathbf{v}|^2)\ d\mathbf{x}+\frac{\left(8\,\varepsilon^2+8\,\varepsilon^3\right)}{\bar{p}}\int_{\mathbb{R}^2}p^2|\nabla p|^2d\mathbf{x}+\\
&\ \left(\frac{\left(4\,\varepsilon^2+4\,\varepsilon^3\right)}{\bar{p}}+8\,\varepsilon^2+32\,\varepsilon^4\right)\int_{\mathbb{R}^2}\left((v_1)^2|\nabla v_1|^2+(v_2)^2|\nabla v_2|^2\right)d\mathbf{x},
\end{aligned}
\end{equation}}where
{\normalsize$$
\begin{aligned}
E_2(t)\equiv &\ \frac12\,\|p\|^2+\frac{\bar{p}}{2}\,\|\mathbf{v}\|^2-\varepsilon\int_{\mathbb{R}^2}p|\mathbf{v}|^2d\mathbf{x}-\frac{2\,\varepsilon+1}{6}\int_{\mathbb{R}^2}p^3d\mathbf{x}+\\
&\ k_1\int_{\mathbb{R}^2}p^4d\mathbf{x}+\frac{2\,\varepsilon^2}{\bar{p}}\int_{\mathbb{R}^2}p^2|\mathbf{v}|^2d\mathbf{x},\\
D_2(t)\equiv &\ \frac12\,\|\nabla p\|^2+\frac{\varepsilon\,\bar{p}}{2}\,\|\nabla\cdot\mathbf{v}\|^2-(2\,\varepsilon+1)\int_{\mathbb{R}^2}p|\nabla p|^2d\mathbf{x}+12\,k_1\int_{\mathbb{R}^2}p^2|\nabla p|^2d\mathbf{x}+\\
&\ \frac{4\,\varepsilon^2}{\bar{p}}\int_{\mathbb{R}^2}|\mathbf{v}|^2|\nabla p|^2 d\mathbf{x}+\frac{2\,\varepsilon^3}{\bar{p}}\int_{\mathbb{R}^2}p^2\left(|\nabla v_1|^2+|\nabla v_2|^2\right)d\mathbf{x}.
\end{aligned}
$$}We observe that the last term on the RHS of \eqref{n6o} is troublesome, due to the expression of {\normalsize$D_2(t)$} contains only the square of $\nabla \mathbf{v}$. Next, we cook up a quantity to dominate such a bad term.

{\bf Step 4.} We write the second equation in \eqref{n1} in the component form as
{\normalsize
\begin{equation}\label{n6g}
\begin{aligned}
&\partial_tv_1-\partial_x p=\varepsilon\Delta v_1-\varepsilon\partial_x(|\mathbf{v}|^2),\\
&\partial_tv_2-\partial_y p=\varepsilon\Delta v_2-\varepsilon\partial_y(|\mathbf{v}|^2).
\end{aligned}
\end{equation}}For any positive constant {\normalsize$k_2$}, taking the {\normalsize$L^2$} inner product of the first equation of \eqref{n6g} with {\normalsize$4\,k_2(v_1)^3$}, the second equation with {\normalsize$4\,k_2(v_2)^3$}, then adding the results, we can show that
{\normalsize
\begin{equation}\label{n6h}
\begin{aligned}
&\ \frac{d}{dt}\left(k_2\int_{\mathbb{R}^2}(v_1)^4+(v_2)^4d\mathbf{x}\right)+12\,k_2\,\varepsilon \int_{\mathbb{R}^2}\left((v_1)^2|\nabla v_1|^2+(v_2)^2|\nabla v_2|^2\right)d\mathbf{x}\\
=&\ 4\,k_2\,\int_{\mathbb{R}^2}\left( (v_1)^3\partial_x p+(v_2)^3\partial_y p\right)d\mathbf{x}
-4\,k_2\,\varepsilon \int_{\mathbb{R}^2}\left( (v_1)^3\partial_x (|\mathbf{v}|^2)+(v_2)^3\partial_y (|\mathbf{v}|^2)\right)d\mathbf{x}.
\end{aligned}
\end{equation}}By adding \eqref{n6h} to \eqref{n6o}, we have
{\normalsize
\begin{equation}\label{n6p}
\begin{aligned}
&\ \frac{d}{dt}E_3(t)+D_3(t) \\
\le &\ \left(2\,\varepsilon+1-12\,k_1+\frac{4\,\varepsilon^2}{\bar{p}}\right)\int_{\mathbb{R}^2}p^2\mathbf{v}\cdot\nabla p\ d\mathbf{x}+\varepsilon\,(2\,\varepsilon+1)\int_{\mathbb{R}^2}p\mathbf{v}\cdot\nabla(|\mathbf{v}|^2) \ d\mathbf{x}-\\
&\ 12k_1\int_{\mathbb{R}^2}p^3\mathbf{v}\cdot\nabla p\ d\mathbf{x}+\frac{4\,\varepsilon^2}{\bar{p}}\int_{\mathbb{R}^2}|\mathbf{v}|^2p \nabla\cdot(p\mathbf{v})\ d\mathbf{x}+\frac{4\,\varepsilon^2}{\bar{p}}\int_{\mathbb{R}^2}|\mathbf{v}|^2p \nabla\cdot\mathbf{v}\ d\mathbf{x}-\\
&\ \frac{4\,\varepsilon^3}{\bar{p}}\int_{\mathbb{R}^2}p^2\mathbf{v}\cdot\nabla(|\mathbf{v}|^2)\ d\mathbf{x}+4\,k_2\int_{\mathbb{R}^2}\left( (v_1)^3\partial_x p+(v_2)^3\partial_y p\right)d\mathbf{x}-\\
&\ 4\,k_2\,\varepsilon \int_{\mathbb{R}^2}\left( (v_1)^3\partial_x (|\mathbf{v}|^2)+(v_2)^3\partial_y (|\mathbf{v}|^2)\right)d\mathbf{x}\ \equiv\ \sum_{k=1}^{8}I_k(t),
\end{aligned}
\end{equation}}where
{\normalsize$$
\begin{aligned}
E_3(t)\equiv &\ \frac12\,\|p\|^2+\frac{\bar{p}}{2}\,\|\mathbf{v}\|^2-\varepsilon\int_{\mathbb{R}^2}p|\mathbf{v}|^2d\mathbf{x}-\frac{2\,\varepsilon+1}{6}\int_{\mathbb{R}^2}p^3d\mathbf{x}+k_1\int_{\mathbb{R}^2}p^4d\mathbf{x}+\\
&\ \frac{2\,\varepsilon^2}{\bar{p}}\int_{\mathbb{R}^2}p^2|\mathbf{v}|^2d\mathbf{x}+k_2\int_{\mathbb{R}^2}(v_1)^4+(v_2)^4d\mathbf{x},\\
D_3(t)\equiv &\ \frac12\,\|\nabla p\|^2+\frac{\varepsilon\,\bar{p}}{2}\,\|\nabla\cdot\mathbf{v}\|^2-(2\,\varepsilon+1)\int_{\mathbb{R}^2}p|\nabla p|^2d\mathbf{x}+\\
&\ \left(12k_1-\frac{(8\,\varepsilon^2+8\,\varepsilon^3)}{\bar{p}}\right)\int_{\mathbb{R}^2}p^2|\nabla p|^2d\mathbf{x}+\\
&\ \frac{4\,\varepsilon^2}{\bar{p}}\int_{\mathbb{R}^2}|\mathbf{v}|^2|\nabla p|^2 d\mathbf{x}+\frac{2\,\varepsilon^3}{\bar{p}}\int_{\mathbb{R}^2}p^2\left(|\nabla v_1|^2+|\nabla v_2|^2\right)d\mathbf{x}+\\
&\ \left(12\,k_2\,\varepsilon-\left[\frac{\left(4\,\varepsilon^2+4\,\varepsilon^3\right)}{\bar{p}}+8\,\varepsilon^2+32\,\varepsilon^4\right]
\right) \int_{\mathbb{R}^2}\left((v_1)^2|\nabla v_1|^2+(v_2)^2|\nabla v_2|^2\right)d\mathbf{x}.
\end{aligned}
$$}By choosing
{\normalsize\begin{equation}\label{k1}
12\,k_1=\frac{8\,\varepsilon^2+8\,\varepsilon^3}{\bar{p}}+3\,(2\,\varepsilon+1)^2,\ \ \ \ \ \ \
12\,k_2=\frac{4\,\varepsilon+4\,\varepsilon^2}{\bar{p}}+9\,\varepsilon+32\,\varepsilon^3,
\end{equation}}we have
{\normalsize\begin{equation}\label{e3d3}
\begin{aligned}
E_3(t)=&\ \frac14\|p\|^2 + \left\|\frac{p}{2}-\frac{2\varepsilon+1}{6}p^2 \right\|^2 + \left(\frac{2(2\varepsilon+1)^2}{9}+\frac{8\varepsilon^2+8\varepsilon^3}{12\bar{p}}\right)\|p\|_{L^4}^4+\\[2mm]
&\ \frac{\bar{p}}{4}\|\mathbf{v}\|^2+\frac{\bar{p}}{4}\left\|\left(1-\frac{2\varepsilon}{\bar{p}}p\right)\mathbf{v} \right\|^2+\frac{\varepsilon^2}{\bar{p}}\|p\mathbf{v}\|^2+\frac{1}{12}\left(\frac{4\,\varepsilon+4\,\varepsilon^2}{\bar{p}}+9\,\varepsilon+32\,\varepsilon^3\right)\|\mathbf{v}\|_{L^4}^4,\\[2mm]
D_3(t) =&\ \frac38\|\nabla p\|^2+\frac18\left\|(1-4\,(2\,\varepsilon+1)\,p)\nabla p\right\|^2+(2\,\varepsilon+1)^2\|p\nabla p\|^2+\\[2mm]
&\ \frac{\varepsilon\,\bar{p}}{2}\,\|\nabla\cdot\mathbf{v}\|^2+\frac{4\,\varepsilon^2}{\bar{p}}\||\mathbf{v}||\nabla p|\|^2 +\frac{2\,\varepsilon^3}{\bar{p}}\|p\nabla\mathbf{v}\|^2+\varepsilon^2\left(\|v_1\nabla v_1\|^2+\|v_2\nabla v_2\|^2\right).
\end{aligned}
\end{equation}}We note that $k_1=\frac14, k_2=0$ when $\varepsilon=0$, and $E_3(0)$ is the same as $H_1$ in the statement of Theorem \ref{thm2}. Next, we carry out energy estimates for the RHS of \eqref{n6p} by applying the Gagliardo-Nirenberg interpolation inequalities, \eqref{gn2d1}--\eqref{gn2d5}.

{\bf Step 5.} By using $\eqref{gn2d1}$, $\eqref{gn2d3}$ and \eqref{n2}, we can show that
{\normalsize\begin{equation}\label{n6aa}
\begin{aligned}
\left|I_1(t)\right|&=\left|\left(2\,\varepsilon+1-12\,k_1+\frac{4\,\varepsilon^2}{\bar{p}}\right)\int_{\mathbb{R}^2}p^2\mathbf{v}\cdot\nabla p\ d\mathbf{x}\right|\\
&\le \left(2\,\varepsilon+1+12\,k_1+\frac{4\,\varepsilon^2}{\bar{p}}\right)\|p\|_{L^8}^2\|\mathbf{v}\|_{L^4}\|\nabla p\|\\
&\le d_1d_3^2\left(2\,\varepsilon+1+12\,k_1+\frac{4\,\varepsilon^2}{\bar{p}}\right)\|\nabla p\|^{\frac32}\|p\|^{\frac12}\|\nabla\cdot\mathbf{v}\|^{\frac12}\|\mathbf{v}\|^{\frac12}\|\nabla p\|\\
&\le d_1d_3^2\left(2\,\varepsilon+1+12\,k_1+\frac{4\,\varepsilon^2}{\bar{p}}\right)(M_1\delta )^{\frac12}\|\nabla p\|^2.
\end{aligned}
\end{equation}}Similarly, we can show that
{\normalsize\begin{equation}\label{n6bb}
\begin{aligned}
\left|I_2(t)\right|&=\left|\varepsilon\,(2\,\varepsilon+1)\int_{\mathbb{R}^2}p\mathbf{v}\cdot\nabla(|\mathbf{v}|^2) \ d\mathbf{x}\right|\\
&\le 2\,(2\,\varepsilon+1)\,\varepsilon\,\|p\|_{L^4}\|\mathbf{v}\|_{L^8}^2\|\nabla\cdot\mathbf{v}\|\\
&\le 2\,d_1d_3^2\,(2\,\varepsilon+1)\,\varepsilon\,\|\nabla p\|^{\frac12}\|p\|^{\frac12}\|\nabla\cdot\mathbf{v}\|^{\frac32}\|\mathbf{v}\|^{\frac12}\|\nabla\cdot\mathbf{v}\|\\
&\le 2\,d_1d_3^2\,(2\,\varepsilon+1)(M_1\delta)^{\frac12}\,(\varepsilon\,\|\nabla\cdot\mathbf{v}\|^2).
\end{aligned}
\end{equation}}By using $\eqref{gn2d4}$, we can show that
{\normalsize\begin{equation}\label{n6cc}
\begin{aligned}
\left|I_3(t)\right|&=\left|12\,k_1\int_{\mathbb{R}^2}p^3\mathbf{v}\cdot\nabla p\ d\mathbf{x}\right|\\
&\le 12\,k_1\,\|p\|_{L^{12}}^3\|\mathbf{v}\|_{L^4}\|\nabla p\|\\
&\le 12\,d_1d_4^2\,k_1\,\|\nabla p\|^{\frac52}\|p\|^{\frac12}\|\nabla\cdot\mathbf{v}\|^{\frac12}\|\mathbf{v}\|^{\frac12}\|\nabla p\|\\
&\le 12\,d_1d_4^2\,k_1\,(M_1\delta^2)^{\frac12}\|\nabla p\|^2\\
&= d_1d_4^2\left(\frac{8\,\varepsilon^2+8\,\varepsilon^3}{\bar{p}}+3\,(2\,\varepsilon+1)^2\right)(M_1\delta^2)^{\frac12}\|\nabla p\|^2.
\end{aligned}
\end{equation}}In a similar fashion, we can show that
{\normalsize\begin{equation}\label{n6dd}
\begin{aligned}
\left|I_4(t)\right|&=\left|\frac{4\,\varepsilon^2}{\bar{p}}\int_{\mathbb{R}^2}|\mathbf{v}|^2p \nabla\cdot(p\mathbf{v})\ d\mathbf{x}\right|\\
&\le \frac{4\,\varepsilon^2}{\bar{p}}\left(\|\mathbf{v}\|_{L^8}^2\|p\|_{L^8}^2\|\nabla\cdot\mathbf{v}\|+\|\mathbf{v}\|_{L^{12}}^3\|p\|_{L^4}\|\nabla p\|\right)\\
&\le \frac{4\,\varepsilon^2}{\bar{p}}\Big(d_3^4\|\nabla\cdot\mathbf{v}\|^{\frac32}\|\mathbf{v}\|^{\frac12}\|\nabla p\|^{\frac32}\|p\|^{\frac12}\|\nabla\cdot\mathbf{v}\|+\\
&\qquad \qquad \qquad d_1d_4^3\|\nabla\cdot\mathbf{v}\|^{\frac52}\|\mathbf{v}\|^{\frac12}\|\nabla p\|^{\frac12}\|p\|^{\frac12}\|\nabla p\|\Big)\\
&\le \frac{4(d_3^4+d_1d_4^3)\,\varepsilon}{\bar{p}}\,(M_1\delta^2)^{\frac12}\,(\varepsilon\,\|\nabla\cdot\mathbf{v}\|^2).
\end{aligned}
\end{equation}}Similar to \eqref{n6aa}, we can show that
{\normalsize\begin{equation}\label{n6ee}
\begin{aligned}
\left|I_5(t)\right|&=\left|\frac{4\,\varepsilon^2}{\bar{p}}\int_{\mathbb{R}^2}|\mathbf{v}|^2p \nabla\cdot\mathbf{v}\ d\mathbf{x}\right|\\
&\le \frac{4\,d_1d_3^2\,\varepsilon}{\bar{p}}\,(M_1\delta)^{\frac12}\,(\varepsilon\,\|\nabla\cdot\mathbf{v}\|^2).
\end{aligned}
\end{equation}}Similar to \eqref{n6dd}, we can show that
{\normalsize\begin{equation}\label{n6ff}
\begin{aligned}
\left|I_6(t)\right|&=\left|\frac{4\,\varepsilon^3}{\bar{p}}\int_{\mathbb{R}^2}p^2\mathbf{v}\cdot\nabla(|\mathbf{v}|^2)\ d\mathbf{x}\right|\\
&\le \frac{8\,d_3^4\,\varepsilon^2}{\bar{p}}\,(M_1\delta^2)^{\frac12}\,(\varepsilon\,\|\nabla\cdot\mathbf{v}\|^2).
\end{aligned}
\end{equation}}By using $\eqref{gn2d2}$, we can show that
{\normalsize
\begin{equation}\label{n6gg}
\begin{aligned}
\left|I_7(t)\right|&=\left|4\,k_2\int_{\mathbb{R}^2}\left( (v_1)^3\partial_x p+(v_2)^3\partial_y p\right)d\mathbf{x}\right|\\
&\le 8\,k_2\,\|\mathbf{v}\|_{L^6}^3\|\nabla p\|\\
&\le 8\,d_2^3\,k_2\,\|\nabla\cdot\mathbf{v}\|^2\|\mathbf{v}\|\|\nabla p\|\\
&\le \frac23\,d_2^3\left(\frac{4+4\,\varepsilon}{\bar{p}}+9+32\,\varepsilon^2\right)(M_1\delta )^{\frac12}\,(\varepsilon\,\|\nabla\cdot\mathbf{v}\|^2).
\end{aligned}
\end{equation}}Similarly, we have
{\normalsize
\begin{equation}\label{n6hh}
\begin{aligned}
\left|I_8(t)\right|&=\left|4\,k_2\,\varepsilon \int_{\mathbb{R}^2}\left( (v_1)^3\partial_x (|\mathbf{v}|^2)+(v_2)^3\partial_y (|\mathbf{v}|^2)\right)d\mathbf{x}\right|\\
&\le 32\,k_2\,\varepsilon\,\|\mathbf{v}\|_{L^8}^4\|\nabla\cdot\mathbf{v}\|\\
&\le 32\,d_3^4\,k_2\,\varepsilon\,\|\nabla\cdot\mathbf{v}\|^3\|\mathbf{v}\|\|\nabla \cdot\mathbf{v}\|\\
&\le 32\,d_3^4\,k_2(M_1\delta^2)^{\frac12}\,(\varepsilon\,\|\nabla\cdot\mathbf{v}\|^2).
\end{aligned}
\end{equation}}This completes the energy estimation for the entire RHS of \eqref{n6p}.

{\bf Step 6.} By substituting \eqref{n6aa}--\eqref{n6hh} into \eqref{n6p}, we obtain
\begin{equation}\label{n6pa}
\frac{d}{dt}E_3(t)+D_3(t) \le (\mathcal{B}_1+\mathcal{B}_3)\|\nabla p\|^2+(\mathcal{B}_2+\mathcal{B}_4+\mathcal{B}_5+\mathcal{B}_6+\mathcal{B}_7+\mathcal{B}_8)(\varepsilon\,\|\nabla\cdot\mathbf{v}\|^2),
\end{equation}
where
\begin{equation}\label{B1}
\begin{aligned}
\mathcal{B}_1 & = d_1d_3^2\left(2\,\varepsilon+1+12\,k_1+\frac{4\,\varepsilon^2}{\bar{p}}\right)(M_1\delta )^{\frac12},  &\mathcal{B}&_2 = 2\,d_1d_3^2\,(2\,\varepsilon+1)(M_1\delta)^{\frac12},\\
\mathcal{B}_3 &= d_1d_4^2\left(\frac{8\,\varepsilon^2+8\,\varepsilon^3}{\bar{p}}+3\,(2\,\varepsilon+1)^2\right)(M_1\delta^2)^{\frac12},  &\mathcal{B}&_4 = \frac{4(d_3^4+d_1d_4^3)\,\varepsilon}{\bar{p}}\,(M_1\delta^2)^{\frac12},\\
\mathcal{B}_5 & = \frac{4\,d_1d_3^2\,\varepsilon}{\bar{p}}\,(M_1\delta)^{\frac12},  &\mathcal{B}&_6 = \frac{8\,d_3^4\,\varepsilon^2}{\bar{p}}\,(M_1\delta^2)^{\frac12},\\
\mathcal{B}_7 & = \frac23\,d_2^3\left(\frac{4+4\,\varepsilon}{\bar{p}}+9+32\,\varepsilon^2\right)(M_1\delta )^{\frac12},  &\mathcal{B}&_8 = 32\,d_3^4\,k_2(M_1\delta^2)^{\frac12}.
\end{aligned}
\end{equation}
Hence, when
\begin{equation}\label{B2}
\begin{aligned}
\mathcal{B}_i &\le \frac{1}{16},\quad i=1,3,\\[2mm]
\mathcal{B}_i &\le \frac{\bar{p}}{24},\quad i=2,4,5,6,7,8,
\end{aligned}
\end{equation}
we get from \eqref{n6pa} that
\begin{equation}\label{n6jj}
\frac{d}{dt}E_3(t)+\left(D_3(t)-\frac18\,\|\nabla p\|^2-\frac{\varepsilon\,\bar{p}}{4}\|\nabla\cdot\mathbf{v}\|^2\right)\le 0.
\end{equation}
In view of \eqref{e3d3} we see that
\begin{equation}\label{n6mm}
D_3(t)-\frac18\,\|\nabla p\|^2-\frac{\varepsilon\,\bar{p}}{4}\|\nabla\cdot\mathbf{v}\|^2\ge \frac14\,\|\nabla p\|^2+\frac{\varepsilon\,\bar{p}}{4}\,\|\nabla\cdot\mathbf{v}\|^2.
\end{equation}
By integrating \eqref{n6jj} with respect to time and using the definition of $E_3(t)$ (cf. \eqref{e3d3}), we find in particular that
{\normalsize
\begin{equation*}\label{n6kk}
\frac14\,\|p(t)\|^2+\frac{\bar{p}}{4}\,\|\mathbf{v}(t)\|^2+\int_0^t\left(\frac14\,\|\nabla p(\tau)\|^2+\frac{\varepsilon\,\bar{p}}{4}\,\|\nabla\cdot\mathbf{v}(\tau)\|^2\right)d\tau\le E_3(0).
\end{equation*}}This implies that
$$
\|p(t)\|^2+\|\mathbf{v}(t)\|^2 \le 4\,(1+1/\bar{p})\,E_3(0),
$$
and
\begin{equation}\label{M_1a}
\int_0^t\left(\|\nabla p(\tau)\|^2+\varepsilon\,\bar{p}\|\nabla\cdot\mathbf{v}(\tau)\|^2\right)d\tau\le 4\,E_3(0).
\end{equation}
In view of \eqref{energy1b} we see that
{\normalsize\begin{equation}\label{M_1}
\|p(t)\|^2+\|\mathbf{v}(t)\|^2 \le M_1-1.
\end{equation}}Next, we estimate the first order derivatives of the solution.

%%%%%%%%%%%%%%%%%%%%%%%%%%%%%%%%%%%%%%%%%%%%%%%%%%%%%%%%%%%%%%%%%%%%%%%%

\subsection{$H^1$-estimate}
We remark that the estimates of the first order derivatives of the solution bear the same level of difficulty as those for the zeroth frequency. Hence, we shall continue with the process of cancellation and coupling through higher order estimates.

{\bf Step 1.} By testing the equations in \eqref{n1} with the {\normalsize$-\Delta$} of the targeting functions, we can show that
{\normalsize
\begin{equation}\label{n10}
\begin{aligned}
&\ \frac{d}{dt}\left(\|\nabla p\|^2+\bar{p}\,\|\nabla\cdot \mathbf{v}\|^2\right)+2\,\|\Delta p\|^2+2\,\varepsilon\,\bar{p}\,\|\Delta\mathbf{v}\|^2\\
=&\ -\int_{\mathbb{R}^2}\left[2\,p(\nabla\cdot\mathbf{v})+2\,\nabla p\cdot\mathbf{v} \right]\Delta p\ d\mathbf{x}+2\,\varepsilon\,\bar{p}\int_{\mathbb{R}^2}\nabla(|\mathbf{v}|^2)\cdot(\Delta\mathbf{v})\ d\mathbf{x}.
\end{aligned}
\end{equation}}A direct calculation by using the first equation in \eqref{n1} shows that
{\normalsize
\begin{equation}\label{n13}
\begin{aligned}
&\ \frac12\frac{d}{dt}\int_{\mathbb{R}^2}p^2\Delta p\ d\mathbf{x}\\
=&\ \int_{\mathbb{R}^2}p\Delta p\left[\Delta p+\nabla\cdot(p\mathbf{v})+\bar{p}\nabla\cdot\mathbf{v}\right]d\mathbf{x}+\int_{\mathbb{R}^2}\frac{p^2}{2}\Delta\left[\Delta p+\nabla\cdot(p\mathbf{v}) +\bar{p} \nabla\cdot\mathbf{v}\right]d\mathbf{x}.
\end{aligned}
\end{equation}}We note that the second term on the RHS of \eqref{n13}
{\normalsize
\begin{equation*}\label{n14}
\int_{\mathbb{R}^2}\frac{p^2}{2}\Delta\left[\Delta p+\nabla\cdot(p\mathbf{v}) + \bar{p}\nabla\cdot\mathbf{v}\right]d\mathbf{x}=\int_{\mathbb{R}^2}\left(p\Delta p+|\nabla p|^2 \right)\left[\Delta p+\nabla\cdot(p\mathbf{v}) + \bar{p}\nabla\cdot\mathbf{v}\right]d\mathbf{x}.
\end{equation*}}So we update \eqref{n13}, after integrating by parts, as
{\normalsize
\begin{equation}\label{n15}
-\frac{d}{dt}\int_{\mathbb{R}^2}p|\nabla p|^2\ d\mathbf{x}=\int_{\mathbb{R}^2}\left(2\,p\Delta p+|\nabla p|^2\right)\left[\Delta p+\nabla\cdot(p\mathbf{v})+\bar{p}\nabla\cdot\mathbf{v}\right]d\mathbf{x}.
\end{equation}}By summing up \eqref{n10}{\normalsize$\times \bar{p}$} and \eqref{n15}, we have
{\normalsize
\begin{equation}\label{n16}
\begin{aligned}
&\frac{d}{dt}\left(\bar{p}\,\|\nabla p\|^2+\bar{p}^2\,\|\nabla\cdot\mathbf{v}\|^2-\int_{\mathbb{R}^2}p|\nabla p|^2\ d\mathbf{x}\right)+\\
&\qquad\qquad\qquad\qquad 2\,\bar{p}\,\|\Delta p\|^2-2\int_{\mathbb{R}^2}p(\Delta p)^2d\mathbf{x}+2\,\varepsilon\,\bar{p}^2\,\|\Delta\mathbf{v}\|^2=N_1(t),
\end{aligned}
\end{equation}}where
{\normalsize
\begin{equation}\label{n17}
\begin{aligned}
N_1(t)&=-2\,\bar{p}\int_{\mathbb{R}^2}(\nabla p\cdot\mathbf{v})\Delta p\ d\mathbf{x}+2\int_{\mathbb{R}^2} p\Delta p\nabla\cdot(p\mathbf{v})\ d\mathbf{x}+\int_{\mathbb{R}^2}|\nabla p|^2\Delta p\ d\mathbf{x}+\\
&\ \ \ \ \ \ \ \ \ \ \ \int_{\mathbb{R}^2}|\nabla p|^2\nabla\cdot(p\mathbf{v})\ d\mathbf{x}+\bar{p}\int_{\mathbb{R}^2}|\nabla p|^2\nabla\cdot\mathbf{v}\ d\mathbf{x}+2\,\varepsilon\,\bar{p}^2\int_{\mathbb{R}^2}\nabla(|\mathbf{v}|^2)\cdot(\Delta\mathbf{v})\ d\mathbf{x}.
\end{aligned}
\end{equation}}We note that the expression inside the parenthesis on the LHS of \eqref{n16} is not necessarily positive. Hence, we need to supply a positive term in order to gain the positivity of the quantity. For this purpose, a direct calculation by using the first equation in \eqref{n1} shows that
{\normalsize
\begin{equation}\label{n15a}
\begin{aligned}
&\ \frac{d}{dt}\int_{\mathbb{R}^2}p^2|\nabla p|^2\ d\mathbf{x}+2\int_{\mathbb{R}^2}p^2|\Delta p|^2d\mathbf{x}\\
=&\ -2\int_{\mathbb{R}^2}p^2\Delta p\left[\nabla\cdot(p\mathbf{v})+\bar{p}\nabla\cdot\mathbf{v}\right]d\mathbf{x}-2\int_{\mathbb{R}^2}p|\nabla p|^2\left[\Delta p+\nabla\cdot(p\mathbf{v})+\bar{p}\nabla\cdot\mathbf{v}\right]d\mathbf{x}.
\end{aligned}
\end{equation}}By adding \eqref{n15a} to \eqref{n16}{\normalsize$\times\bar{p}$}, we find
{\normalsize
\begin{equation}\label{n16a}
\begin{aligned}
&\ \frac{d}{dt}\left(\bar{p}^2\,\|\nabla p\|^2+\bar{p}^3\,\|\nabla\cdot\mathbf{v}\|^2-\bar{p}\int_{\mathbb{R}^2}p|\nabla p|^2\ d\mathbf{x}+\int_{\mathbb{R}^2}p^2|\nabla p|^2\ d\mathbf{x}\right)+\\
&\qquad\qquad 2\,\bar{p}^2\,\|\Delta p\|^2-2\,\bar{p}\int_{\mathbb{R}^2}p(\Delta p)^2d\mathbf{x}+2\,\varepsilon\,\bar{p}^3\,\|\Delta\mathbf{v}\|^2+2\int_{\mathbb{R}^2}p^2|\Delta p|^2d\mathbf{x}=N_2(t),
\end{aligned}
\end{equation}}where
{\normalsize$$
\begin{aligned}
N_2(t)=&\ \bar{p}\,N_1(t)-2\int_{\mathbb{R}^2}p^2\Delta p\nabla\cdot(p\mathbf{v})d\mathbf{x}-2\,\bar{p}\int_{\mathbb{R}^2}p^2\Delta p\nabla\cdot\mathbf{v}\ d\mathbf{x}-2\int_{\mathbb{R}^2}p|\nabla p|^2\Delta p\ d\mathbf{x}- \\
&\ 2\int_{\mathbb{R}^2}p|\nabla p|^2\nabla\cdot(p\mathbf{v})d\mathbf{x}-2\,\bar{p}\int_{\mathbb{R}^2}p|\nabla p|^2\nabla\cdot\mathbf{v}\ d\mathbf{x}\equiv \sum_{k=1}^{11}J_k(t).
\end{aligned}
$$}Next, we carry out energy estimates for {\normalsize$N_2(t)$}.

{\bf Step 2.} By using $\eqref{gn2d1}$ and Young's inequality, we can show that
{\normalsize
\begin{equation}\label{n18}
\begin{aligned}
|J_1(t)|&=2\,\bar{p}^2\left|\int_{\mathbb{R}^2}(\nabla p\cdot\mathbf{v})\Delta p\ d\mathbf{x}\right|\\
&\le 2\,\bar{p}^2\,\|\nabla p\|_{L^4}\|\mathbf{v}\|_{L^4}\|\Delta p\|\\
&\le 2\,d_1^2\,\bar{p}^2\,\|\nabla p\|^{\frac12}\|\Delta p\|^{\frac32}\|\mathbf{v}\|^{\frac12}\|\nabla\cdot\mathbf{v}\|^{\frac12}\\
&\le \frac{\bar{p}^2}{12}\,\|\Delta p\|^2+2916\,d_1^8\,\bar{p}^2\,\|\nabla p\|^2\|\mathbf{v}\|^2\|\nabla\cdot\mathbf{v}\|^2\\
&\le \frac{\bar{p}^2}{12}\,\|\Delta p\|^2+2916\,d_1^8\,\bar{p}^2\,M_1\,\|\nabla p\|^2\|\nabla\cdot\mathbf{v}\|^2.
\end{aligned}
\end{equation}}By using $\eqref{gn2d1}$ and $\eqref{gn2d5}$, we can show that
{\normalsize
\begin{equation}\label{n19}
\begin{aligned}
|J_2(t)|&=2\,\bar{p}\left|\int_{\mathbb{R}^2} p\Delta p\nabla\cdot(p\mathbf{v})\ d\mathbf{x}\right|\\
&=2\,\bar{p}\left|\int_{\mathbb{R}^2} p\Delta p\left(p\nabla\cdot\mathbf{v}+\mathbf{v}\cdot\nabla p\right)\ d\mathbf{x}\right|\\
&\le 2\,\bar{p}\left(\|\Delta p\|\|p\|_{L^\infty}^2\|\nabla\cdot \mathbf{v}\|+\|\Delta p\|\|p\|_{L^\infty}\|\mathbf{v}\|_{L^4}\|\nabla p\|_{L^4}\right)\\
&\le 2\,\bar{p}\left(d_5^2\,\|\Delta p\|^2\|p\|\|\nabla\cdot\mathbf{v}\|+d_1^2d_5\,\|\Delta p\|^2\|p\|^{\frac12}\|\mathbf{v}\|^{\frac12}\|\nabla\cdot\mathbf{v}\|^{\frac12}\|\nabla p\|^{\frac12}\right)\\
&\le 2(d_5^2+d_1^2d_5)\,\bar{p}\,(M_1\delta)^{\frac12}\|\Delta p\|^2.
\end{aligned}
\end{equation}}By using $\eqref{gn2d1}$, we can show that
{\normalsize
\begin{equation}\label{n20}
\begin{aligned}
|J_3(t)|&=\bar{p}\left|\int_{\mathbb{R}^2}|\nabla p|^2\Delta p\ d\mathbf{x}\right|\\
&\le \bar{p}\,\|\nabla p\|^2_{L^4}\|\Delta p\|\\
&\le d_1^2\,\bar{p}\,\|\nabla p\|\|\Delta p\|^2\\
&\le d_1^2\,\bar{p}\,\delta^{\frac12}\|\Delta p\|^2.
\end{aligned}
\end{equation}}By using similar arguments as in \eqref{n19}, we can show that
{\normalsize
\begin{equation}\label{n22}
\begin{aligned}
|J_4(t)|&=\bar{p}\left|\int_{\mathbb{R}^2}|\nabla p|^2\nabla\cdot(p\mathbf{v})\,d\mathbf{x}\right|\\
&=2\,\bar{p}\left|\int_{\mathbb{R}^2}\left(\nabla p \cdot \mathbb{H}(p)\right)\cdot (p\mathbf{v})\,d\mathbf{x} \right|\\
&\le 2\,\bar{p}\,\|\Delta p\|\|p\|_{L^\infty}\|\nabla p\|_{L^4}\|\mathbf{v}\|_{L^4}\\
&\le 2\,d_1^2d_5\,\bar{p}\,\|\Delta p\|^2\|p\|^{\frac12}\|\nabla p\|^{\frac12}\|\mathbf{v}\|^{\frac12}\|\nabla\cdot\mathbf{v}\|^{\frac12}\\
&\le 2\,d_1^2d_5\,\bar{p}\,(M_1\delta)^{\frac12}\|\Delta p\|^2,
\end{aligned}
\end{equation}}where {\normalsize$\mathbb{H}(p)$} denotes the Hessian matrix of {\normalsize$p$}. Similar to \eqref{n20}, we can show that
{\normalsize
\begin{equation}\label{n23}
\begin{aligned}
|J_5(t)|&=\bar{p}^2\left|\int_{\mathbb{R}^2}|\nabla p|^2\nabla\cdot\mathbf{v}\ d\mathbf{x}\right|\\
&\le \bar{p}^2\,\|\nabla p\|^2_{L^4}\|\nabla\cdot\mathbf{v}\|\\
&\le d_1^2\,\bar{p}^2\,\|\nabla p\|\|\Delta p\|\|\nabla\cdot\mathbf{v}\|\\
&\le \frac{\bar{p}^2}{12}\,\|\Delta p\|^2+3\,d_1^4\,\bar{p}^2\,\|\nabla p\|^2\|\nabla\cdot\mathbf{v}\|^2.
\end{aligned}
\end{equation}}Again, by using \eqref{gn2d1}, we can show that
{\normalsize
\begin{equation}\label{n23a}
\begin{aligned}
|J_6(t)|&=2\,\varepsilon\,\bar{p}^3\left|\int_{\mathbb{R}^2}\nabla(|\mathbf{v}|^2)\cdot(\Delta\mathbf{v})\ d\mathbf{x}\right|\\
&\le 4\,\varepsilon\,\bar{p}^3\,\|\mathbf{v}\|_{L^4}\|\nabla\mathbf{v}\|_{L^4}\|\Delta\mathbf{v}\|\\
&\le 4\,d_1^2\,\varepsilon\,\bar{p}^3\,\|\mathbf{v}\|^{\frac12}\|\nabla\cdot\mathbf{v}\|\|\Delta\mathbf{v}\|^{\frac32}\\
&\le \varepsilon\,\bar{p}^3\, \|\Delta\mathbf{v}\|^2+ 27\,d_1^8\,\varepsilon\,\bar{p}^3\, \|\mathbf{v}\|^2\|\nabla\cdot\mathbf{v}\|^4\\
&\le \varepsilon\,\bar{p}^3\, \|\Delta\mathbf{v}\|^2+ 27\,d_1^8\, M_1\,\bar{p}^2\,(\varepsilon\,\bar{p}\, \|\nabla\cdot\mathbf{v}\|^2)\,\|\nabla\cdot\mathbf{v}\|^2.
\end{aligned}
\end{equation}}Similar to \eqref{n19}, we can show that
{\normalsize
\begin{equation}\label{n23b}
\begin{aligned}
|J_7(t)|&=2\left|\int_{\mathbb{R}^2}p^2\Delta p\nabla\cdot(p\mathbf{v})d\mathbf{x}\right|\\
&\le 2\left(\|p\Delta p\|\|p\|^2_{L^\infty}\|\nabla\cdot\mathbf{v}\|+\|p\Delta p\|\|p\|_{L^\infty}\|\nabla p\|_{L^4}\|\mathbf{v}\|_{L^4}\right)\\
&\le 2\left(d_5^2\,\|p\Delta p\|\|p\|\|\Delta p\|\|\nabla\cdot\mathbf{v}\|+d_1^2d_5\,\|p\Delta p\|\|p\|^{\frac12}\|\Delta p\|\|\nabla p\|^{\frac12}\|\mathbf{v}\|^{\frac12}\|\nabla\cdot\mathbf{v}\|^{\frac12}\right)\\
&\le 2(d_5^2+d_1^2d_5)\,(M_1\delta)^{\frac12}\,\|p\Delta p\|\|\Delta p\|\\
&\le (d_5^2+d_1^2d_5)(M_1\delta)^{\frac12}\left(\|p\Delta p\|^2+\|\Delta p\|^2\right).
\end{aligned}
\end{equation}}By using \eqref{gn2d5}, we can show that
{\normalsize
\begin{equation}\label{n23c}
\begin{aligned}
|J_8(t)|&=2\,\bar{p}\left|\int_{\mathbb{R}^2}p^2\Delta p\nabla\cdot\mathbf{v}\ d\mathbf{x}\right|\\
&\le 2\,\bar{p}\,\|\Delta p\|\|p\|_{L^\infty}^2\|\nabla\cdot\mathbf{v}\|\\
&\le 2\,d_5^2\,\bar{p}\,\|\Delta p\|\|p\|\|\Delta p\|\|\nabla\cdot\mathbf{v}\|\\
&\le 2\,d_5^2\,\bar{p}\,(M_1\delta)^{\frac12}\|\Delta p\|^2.
\end{aligned}
\end{equation}}By using $\eqref{gn2d1}$, we can show that
{\normalsize
\begin{equation}\label{n23d}
\begin{aligned}
|J_9(t)|&=2\left|\int_{\mathbb{R}^2}p|\nabla p|^2\Delta p\ d\mathbf{x}\right|\\
&\le 2\,\|p\Delta p\|\|\nabla p\|_{L^4}^2\\
&\le 2\,d_1^2\,\|p\Delta p\|\|\nabla p\|\|\Delta p\|\\
&\le d_1^2\,\delta^{\frac12}\left(\|p\Delta p\|^2+\|\Delta p\|^2\right).
\end{aligned}
\end{equation}}By using $\eqref{gn2d1}$ and $\eqref{gn2d5}$, we can show that
{\normalsize
\begin{equation}\label{n23e}
\begin{aligned}
&\ |J_{10}(t)|\\
=&\ 2\left|\int_{\mathbb{R}^2}p|\nabla p|^2\nabla\cdot(p\mathbf{v})d\mathbf{x}\right|\\
\le &\ 2\left(\|p\|_{L^\infty}^2\|\nabla p\|_{L^4}^2\|\nabla\cdot\mathbf{v}\|+\|p\|_{L^\infty}\|\nabla p\|_{L^4}^3\|\mathbf{v}\|_{L^4}\right)\\
\le &\ 2\left(d_1^2d_5^2\,\|p\|\|\Delta p\|\|\nabla p\|\|\Delta p\|\|\nabla\cdot\mathbf{v}\|+d_1^4d_5\,\|p\|^{\frac12}\|\Delta p\|^{\frac12}\|\nabla p\|^{\frac32}\|\Delta p\|^{\frac32}\|\mathbf{v}\|^{\frac12}\|\nabla\cdot\mathbf{v}\|^{\frac12}\right)\\
\le &\ 2(d_1^2d_5^2+d_1^4d_5)(M_1\delta^2)^{\frac12}\|\Delta p\|^2.
\end{aligned}
\end{equation}}Again, by using $\eqref{gn2d1}$ and $\eqref{gn2d5}$, we can show that
{\normalsize
\begin{equation}\label{n23f}
\begin{aligned}
|J_{11}(t)|&=2\,\bar{p}\left|\int_{\mathbb{R}^2}p|\nabla p|^2\nabla\cdot\mathbf{v}\ d\mathbf{x}\right|\\
&\le 2\,\bar{p}\,\|p\|_{L^\infty}\|\nabla p\|_{L^4}^2\|\nabla\cdot\mathbf{v}\|\\
&\le 2d_1^2d_5\,\bar{p}\,\|p\|^{\frac12}\|\Delta p\|^{\frac12}\|\nabla p\|\|\Delta p\|\|\nabla\cdot\mathbf{v}\|\\
&\le \frac{\bar{p}^2}{12}\,\|\Delta p\|^2+2916\, d_1^8d_5^4(\bar{p})^{-1}\|p\|^2\|\nabla p\|^4\|\nabla\cdot\mathbf{v}\|^4\\
&\le \frac{\bar{p}^2}{12}\,\|\Delta p\|^2+2916\, d_1^8d_5^4(\bar{p})^{-1}(M_1\delta^2)\|\nabla p\|^2\|\nabla\cdot\mathbf{v}\|^2.
\end{aligned}
\end{equation}}This completes the estimate for the entire RHS of \eqref{n16a}.

{\bf Step 3.} By assembling, \eqref{n18}--\eqref{n23f}, we find
{\normalsize
\begin{equation*}
\begin{aligned}
|N_2(t)| \le &\ \varepsilon\,\bar{p}^3\,\|\Delta\mathbf{v}\|^2+\left[(\mathcal{C}_1+\mathcal{C}_2+\mathcal{D}_{1})\|\nabla p\|^2+\mathcal{C}_3(\varepsilon\,\bar{p}\, \|\nabla\cdot\mathbf{v}\|^2)\right]\|\nabla\cdot\mathbf{v}\|^2\\
&\left(\frac{\bar{p}}{4}+\mathcal{D}_2+\mathcal{D}_3+\mathcal{D}_4+\mathcal{D}_5+\mathcal{D}_6+\mathcal{D}_7+\mathcal{D}_{8}\right)\|\Delta p\|^2+(\mathcal{D}_5+\mathcal{D}_7)\|p\Delta p\|^2,
\end{aligned}
\end{equation*}}where
\begin{equation}\label{D1}
\begin{aligned}
\mathcal{C}_1&=2916\,d_1^8\,\bar{p}^2\,M_1, &\mathcal{C}&_2=3\,d_1^4\,\bar{p}^2, &\mathcal{C}&_3=27\,d_1^8\, M_1\,\bar{p}^2,\\
\mathcal{D}_{1}&=2916\, d_1^8d_5^4(\bar{p})^{-1}(M_1\delta^2), &\mathcal{D}&_2=2(d_5^2+d_1^2d_5)\,\bar{p}\,(M_1\delta)^{\frac12}, &\mathcal{D}&_3=d_1^2\,\bar{p}\,\delta^{\frac12},\\
\mathcal{D}_4&=2\,d_1^2d_5\,\bar{p}\,(M_1\delta)^{\frac12},
&\mathcal{D}&_5=(d_5^2+d_1^2d_5)(M_1\delta)^{\frac12}, &\mathcal{D}&_6=2\,d_5^2\,\bar{p}\,(M_1\delta)^{\frac12}, \\
\mathcal{D}_7&=d_1^2\,\delta^{\frac12}, &\mathcal{D}&_{8}=2(d_1^2d_5^2+d_1^4d_5)(M_1\delta^2)^{\frac12}.
\end{aligned}
\end{equation}
Hence, when
\begin{equation}\label{D2}
\begin{aligned}
M_1\delta^2 &\le 1,\qquad \mathcal{D}_j \le \frac{\bar{p}}{28},\quad j=2,...,8, \qquad \mathrm{and}\\
\mathcal{D}_j &\le \min\left\{\frac{\bar{p}}{28},\ \frac{1}{4}\right\}, \quad j=5,7,
\end{aligned}
\end{equation}
it holds that
{\normalsize
\begin{equation}\label{n25}
\begin{aligned}
|N_2(t)| \le &\ \varepsilon\,\bar{p}^3\,\|\Delta\mathbf{v}\|^2+\frac{\bar{p}}{2}\|\Delta p\|^2+\frac12\|p\Delta p\|^2+\\
&\ \left[(\mathcal{C}_1+\mathcal{C}_2+2916\, d_1^8d_5^4(\bar{p})^{-1})\|\nabla p\|^2+\mathcal{C}_3(\varepsilon\,\bar{p}\, \|\nabla\cdot\mathbf{v}\|^2)\right]\|\nabla\cdot\mathbf{v}\|^2.
\end{aligned}
\end{equation}}By substituting \eqref{n25} into \eqref{n16a}, we obtain
{\normalsize
\begin{equation}\label{n26}
\begin{aligned}
&\ \frac{d}{dt}E_4(t)+D_4(t)\\
\le&\ \left[(\mathcal{C}_1+\mathcal{C}_2+2916\, d_1^8d_5^4(\bar{p})^{-1})\|\nabla p\|^2+\mathcal{C}_3(\varepsilon\,\bar{p}\, \|\nabla\cdot\mathbf{v}\|^2)\right]\|\nabla\cdot\mathbf{v}\|^2,
\end{aligned}
\end{equation}}where
{\normalsize
\begin{equation}\label{e4d4}
\begin{aligned}
E_4(t)&\equiv \bar{p}^2\|\nabla p\|^2+\bar{p}^3\|\nabla\cdot\mathbf{v}\|^2-\bar{p}\int_{\mathbb{R}^2}p|\nabla p|^2\ d\mathbf{x}+\int_{\mathbb{R}^2}p^2|\nabla p|^2\ d\mathbf{x}\\
&=\frac{\bar{p}^2}{2}\,\|\nabla p\|^2+\bar{p}^3\|\nabla\cdot\mathbf{v}\|^2+\frac12\,\|\bar{p}\nabla p-p\nabla p\|^2+\frac12\,\|p\nabla p\|^2,\\
D_4(t)&\equiv \frac{3\,\bar{p}^2}{2}\,\|\Delta p\|^2-2\,\bar{p}\int_{\mathbb{R}^2}p(\Delta p)^2d\mathbf{x}+\varepsilon\,\bar{p}^3\|\Delta\mathbf{v}\|^2+\frac32\int_{\mathbb{R}^2}p^2|\Delta p|^2d\mathbf{x}\\
&=\frac{\bar{p}^2}{2}\,\|\Delta p\|^2+\|\bar{p}\Delta p-p\Delta p\|^2+\varepsilon\,\bar{p}^3\|\Delta\mathbf{v}\|^2+\frac12\,\|p\Delta p\|^2.
\end{aligned}
\end{equation}}We note that $E_4(0)$ is the same as $H_2$ in the statement of Theorem \ref{thm2}. We also note that {\normalsize$\bar{p}^3\|\nabla\cdot\mathbf{v}\|^2\le E_4(t)$}. Hence, we update \eqref{n26} as
{\normalsize
\begin{equation}\label{n26a}
\begin{aligned}
\frac{d}{dt}E_4(t)+D_4(t)\le \frac{1}{\bar{p}}\left[(\mathcal{C}_1+\mathcal{C}_2+2916\, d_1^8d_5^4(\bar{p})^{-1})\|\nabla p\|^2+\mathcal{C}_3(\varepsilon\,\bar{p}\, \|\nabla\cdot\mathbf{v}\|^2)\right]E_4(t).
\end{aligned}
\end{equation}}By applying the Gronwall inequality to \eqref{n26a} and using \eqref{M_1a}, we find
{\normalsize
\begin{equation}\label{n30}
\begin{aligned}
E_4(t) \le E_4(0)\exp\left\{\frac{4}{\bar{p}}\,(\mathcal{C}_1+\mathcal{C}_2+\mathcal{C}_3+2916\, d_1^8d_5^4(\bar{p})^{-1})\,E_3(0)\right\} \equiv \mathcal{C}_4.
\end{aligned}
\end{equation}}Therefore, by the definition of $E_4(t)$, it holds that
{\normalsize
\begin{equation}\label{n31}
\|\nabla p(t)\|^2+\|\nabla \cdot\mathbf{v}(t)\|^2\le \frac{2\,\bar{p}+1}{\bar{p}^3}\,\mathcal{C}_4.
\end{equation}}In view of \eqref{energy1b} and \eqref{n31} we see that
\begin{equation}\label{delta}
\|\nabla p(t)\|^2+\|\nabla \cdot\mathbf{v}(t)\|^2 \le \frac{\delta}{2}.
\end{equation}
In addition, by substituting \eqref{n30} into \eqref{n26a} then integrating with respect to $t$, we have
$$
\begin{aligned}
\int_0^t D_4(\tau)d\tau \le E_4(0)+\frac{4}{\bar{p}}\,(\mathcal{C}_1+\mathcal{C}_2+\mathcal{C}_3+2916\, d_1^8d_5^4(\bar{p})^{-1})\,E_3(0)\,\mathcal{C}_4.
\end{aligned}
$$
According to \eqref{e4d4}, we have
\begin{equation}\label{n32}
\begin{aligned}
&\ \int_0^t\left(\|\Delta p(\tau)(t)\|^2+\varepsilon\,\bar{p}\,\|\Delta\mathbf{v}(\tau)\|^2\right)d\tau \\
\le &\ \frac{3}{\bar{p}^3}\left(E_4(0)+\frac{4}{\bar{p}}\,(\mathcal{C}_1+\mathcal{C}_2+\mathcal{C}_3+2916\, d_1^8d_5^4(\bar{p})^{-1})\,E_3(0)\,\mathcal{C}_4\right)\equiv \mathcal{C}_5,
\end{aligned}
\end{equation}
where the constant $\mathcal{C}_5$ is independent of {\normalsize$t$} and remains bounded as {\normalsize$\varepsilon\to0$}. Next, we estimate the second order spatial derivatives of the solution.

%%%%%%%%%%%%%%%%%%%%%%%%%%%%%%%%%%%%%%%%%%%%%%%%%%%%%%%%%%%%%%%%%%%%%%%%

\subsection{$H^2$-estimate}
By applying {\normalsize$\Delta$} to the equations in \eqref{n1}, then taking the {\normalsize$L^2$} inner products of the resulting equations with {\normalsize$\Delta$} of the targeting functions, we obtain
{\normalsize
\begin{equation}\label{nn90}
\begin{aligned}
&\frac{d}{dt}\left(\|\Delta p\|^2+\bar{p}\,\|\Delta \mathbf{v}\|^2\right)+2\,\|\nabla\Delta p\|^2+2\,\varepsilon\,\bar{p}\,\|\Delta(\nabla\cdot\mathbf{v})\|^2\\
=&-2\int_{\mathbb{R}^2}\nabla\left(\nabla\cdot(p\mathbf{v})\right)\cdot(\nabla\Delta p)\ d\mathbf{x}+2\,\varepsilon\,\bar{p}\int_{\mathbb{R}^2}\Delta(|\mathbf{v}|^2)\Delta(\nabla\cdot\mathbf{v})d\mathbf{x}\\
=&-2\int_{\mathbb{R}^2}\left(\nabla(\nabla p\cdot\mathbf{v})+(\nabla\cdot\mathbf{v})\nabla p\right)\cdot(\nabla\Delta p)\ d\mathbf{x}-2\int_{\mathbb{R}^2}p\nabla(\nabla\cdot\mathbf{v})\cdot(\nabla\Delta p)\ d\mathbf{x}+\\
&4\,\varepsilon\,\bar{p}\int_{\mathbb{R}^2}\left(\mathbf{v}\cdot\Delta\mathbf{v}+|\nabla v_1|^2+|\nabla v_2|^2\right)\Delta(\nabla\cdot\mathbf{v})d\mathbf{x}.
\end{aligned}
\end{equation}}By direct calculations, we can show that
{\normalsize
\begin{equation}\label{nn91}
\begin{aligned}
&\ -\frac{d}{dt}\int_{\mathbb{R}^2}p(\Delta p)^2d\mathbf{x}-2\int_{\mathbb{R}^2}p|\nabla\Delta p|^2d\mathbf{x}\\
=&\ 2\int_{\mathbb{R}^2}\Delta p\nabla p\cdot\nabla(\nabla\cdot(p\mathbf{v}))d\mathbf{x}+2\int_{\mathbb{R}^2}p\nabla\Delta p\cdot\nabla(\nabla\cdot(p\mathbf{v}))d\mathbf{x}+\\
&\ 2\,\bar{p}\int_{\mathbb{R}^2}\Delta p\nabla p\cdot\nabla(\nabla\cdot\mathbf{v})d\mathbf{x}+2\,\bar{p}\int_{\mathbb{R}^2}p\nabla\Delta p\cdot\nabla(\nabla\cdot\mathbf{v})d\mathbf{x}+\\
&\ 2\int_{\mathbb{R}^2}\Delta p\nabla p\cdot\nabla\Delta p\ d\mathbf{x}-\int_{\mathbb{R}^2}(\Delta p)^2\nabla\cdot(p\mathbf{v})d\mathbf{x}-\\
&\ \bar{p}\int_{\mathbb{R}^2}(\Delta p)^2\nabla\cdot\mathbf{v}\ d\mathbf{x}-\int_{\mathbb{R}^2}(\Delta p)^3d\mathbf{x},
\end{aligned}
\end{equation}}and
{\normalsize
\begin{equation}\label{nn92}
\begin{aligned}
&\ \frac{d}{dt}\int_{\mathbb{R}^2}p^2(\Delta p)^2d\mathbf{x}+2\int_{\mathbb{R}^2}p^2|\nabla\Delta p|^2d\mathbf{x}\\
=&\ -4\int_{\mathbb{R}^2}p\Delta p\nabla p\cdot\nabla(\nabla\cdot(p\mathbf{v}))d\mathbf{x}-2\int_{\mathbb{R}^2}p^2\nabla\Delta p\cdot\nabla(\nabla\cdot(p\mathbf{v}))d\mathbf{x}-\\
&\ 4\,\bar{p}\int_{\mathbb{R}^2}p\Delta p\nabla p\cdot\nabla(\nabla\cdot\mathbf{v})d\mathbf{x}-2\,\bar{p}\int_{\mathbb{R}^2}p^2\nabla\Delta p\cdot\nabla(\nabla\cdot\mathbf{v})d\mathbf{x}-\\
&\ 4\int_{\mathbb{R}^2}p\Delta p\nabla p\cdot\nabla\Delta p\ d\mathbf{x}+2\int_{\mathbb{R}^2}p(\Delta p)^2\nabla\cdot(p\mathbf{v})d\mathbf{x}+\\
&\ 2\,\bar{p}\int_{\mathbb{R}^2}p(\Delta p)^2\nabla\cdot\mathbf{v}\ d\mathbf{x}+2\int_{\mathbb{R}^2}p(\Delta p)^3d\mathbf{x}.
\end{aligned}
\end{equation}}The operation: \eqref{nn90}{\normalsize$\times\bar{p}$} + \eqref{nn91} + \eqref{nn92}{\normalsize$\times 1/\bar{p}$} then yields
{\normalsize
\begin{equation}\label{nn9}
\begin{aligned}
&\ \frac{d}{dt}\left(\bar{p}\,\|\Delta p\|^2+\bar{p}^2\|\Delta \mathbf{v}\|^2-\int_{\mathbb{R}^2}p(\Delta p)^2d\mathbf{x}+\frac{1}{\bar{p}}\|p\Delta p\|^2\right)+\\[2mm]
&\quad\quad2\,\bar{p}\,\|\nabla\Delta p\|^2-2\int_{\mathbb{R}^2}p|\nabla\Delta p|^2d\mathbf{x}+\frac{2}{\bar{p}}\|p\nabla\Delta p\|^2+2\,\varepsilon\,\bar{p}^2\|\Delta(\nabla\cdot\mathbf{v})\|^2 \\
=&\ \sum_{k=1}^{19}R_k(t),
\end{aligned}
\end{equation}}where
$$
\begin{aligned}
R_1(t)&=-2\,\bar{p}\int_{\mathbb{R}^2}\left(\nabla(\nabla p\cdot\mathbf{v})\right)\cdot(\nabla\Delta p)\ d\mathbf{x}, &R&_2(t)=-2\,\bar{p}\int_{\mathbb{R}^2}\left((\nabla\cdot\mathbf{v})\nabla p\right)\cdot(\nabla\Delta p)\ d\mathbf{x},\\
R_3(t)&=4\,\varepsilon\,\bar{p}^2\int_{\mathbb{R}^2}\left(\mathbf{v}\cdot\Delta\mathbf{v}\right)\Delta(\nabla\cdot\mathbf{v})d\mathbf{x}, &R&_4(t)=4\,\varepsilon\,\bar{p}^2\int_{\mathbb{R}^2}\left(|\nabla \mathbf{v}|^2\right)\Delta(\nabla\cdot\mathbf{v})d\mathbf{x},\\
R_5(t)&=2\int_{\mathbb{R}^2}\Delta p\nabla p\cdot\nabla(\nabla\cdot(p\mathbf{v}))d\mathbf{x}, &R&_6(t)=2\int_{\mathbb{R}^2}p\nabla\Delta p\cdot\nabla(\nabla\cdot(p\mathbf{v}))d\mathbf{x},\\
R_7(t)&=2\,\bar{p}\int_{\mathbb{R}^2}\Delta p\nabla p\cdot\nabla(\nabla\cdot\mathbf{v})d\mathbf{x}, &R&_8(t)=2\int_{\mathbb{R}^2}\Delta p\nabla p\cdot\nabla\Delta p\ d\mathbf{x},\\
R_9(t)&=-\int_{\mathbb{R}^2}(\Delta p)^2\nabla\cdot(p\mathbf{v})d\mathbf{x}, &R&_{10}(t)=-\bar{p}\int_{\mathbb{R}^2}(\Delta p)^2\nabla\cdot\mathbf{v}\ d\mathbf{x},\\
R_{11}(t)&=-\int_{\mathbb{R}^2}(\Delta p)^3d\mathbf{x}, &R&_{12}(t)=-\frac{4}{\bar{p}}\int_{\mathbb{R}^2}p\Delta p\nabla p\cdot\nabla(\nabla\cdot(p\mathbf{v}))d\mathbf{x},\\
R_{13}(t)&=-\frac{2}{\bar{p}}\int_{\mathbb{R}^2}p^2\nabla\Delta p\cdot\nabla(\nabla\cdot(p\mathbf{v}))d\mathbf{x}, &R&_{14}(t)=-4\int_{\mathbb{R}^2}p\Delta p\nabla p\cdot\nabla(\nabla\cdot\mathbf{v})d\mathbf{x},\\
R_{15}(t)&=-2\int_{\mathbb{R}^2}p^2\nabla\Delta p\cdot\nabla(\nabla\cdot\mathbf{v})d\mathbf{x}, &R&_{16}(t)=-\frac{4}{\bar{p}}\int_{\mathbb{R}^2}p\Delta p\nabla p\cdot\nabla\Delta p\ d\mathbf{x},\\
R_{17}(t)&=\frac{2}{\bar{p}}\int_{\mathbb{R}^2}p(\Delta p)^2\nabla\cdot(p\mathbf{v})d\mathbf{x}, &R&_{18}(t)=2\int_{\mathbb{R}^2}p(\Delta p)^2\nabla\cdot\mathbf{v}\ d\mathbf{x},\\
R_{19}(t)&=\frac{2}{\bar{p}}\int_{\mathbb{R}^2}p(\Delta p)^3d\mathbf{x}.
\end{aligned}
$$
Next, we shall estimate the RHS of \eqref{nn9} term by term. However, due to a large number of terms to be estimated, the detailed analysis is long and not reader-friendly. To simplify the presentation, we only present the resulting estimate, while the detailed arguments are left in the Appendix. Indeed, after a series of energy estimates by using \eqref{gn2d1}--\eqref{gn2d5}, we can show that
{\normalsize
\begin{equation}\label{nn9r}
\begin{aligned}
\sum_{k=1}^{19}|R_k(t)| \le &\frac{\bar{p}}{2}\,\|\nabla\Delta p\|^2+\frac{1}{2\bar{p}}\,\|p\nabla\Delta p\|^2+\varepsilon\,\bar{p}^2\|\Delta(\nabla\cdot\mathbf{v})\|^2+\mathcal{K}_1\left(\|\nabla p\|^2+\|\Delta p\|^2\right)+\\
&\mathcal{K}_2\left(\|\nabla p\|^2+\|\Delta p\|^2+\varepsilon\,\bar{p}\,\|\nabla\cdot\mathbf{v}\|^2\right)\left(\bar{p}\|\Delta p\|^2+\bar{p}^2\|\Delta \mathbf{v}\|^2\right),
\end{aligned}
\end{equation}}for some constants $\mathcal{K}_1$ and $\mathcal{K}_2$ which depend only on $M_1,\delta,\bar{p}$ and the constants in \eqref{gn2d1}--\eqref{gn2d5}, and remain bounded as {\normalsize$\varepsilon\to 0$}. By substituting \eqref{nn9r} into \eqref{nn9}, we have
\begin{equation*}
\begin{aligned}
&\ \frac{d}{dt}\left(\bar{p}\,\|\Delta p\|^2+\bar{p}^2\|\Delta \mathbf{v}\|^2-\int_{\mathbb{R}^2}p(\Delta p)^2d\mathbf{x}+\frac{1}{\bar{p}}\|p\Delta p\|^2\right)+\\[2mm]
&\quad\quad2\,\bar{p}\,\|\nabla\Delta p\|^2-2\int_{\mathbb{R}^2}p|\nabla\Delta p|^2d\mathbf{x}+\frac{2}{\bar{p}}\|p\nabla\Delta p\|^2+2\,\varepsilon\,\bar{p}^2\|\Delta(\nabla\cdot\mathbf{v})\|^2 \\
\le&\ \frac{\bar{p}}{2}\,\|\nabla\Delta p\|^2+\frac{1}{2\bar{p}}\,\|p\nabla\Delta p\|^2+\varepsilon\,\bar{p}^2\|\Delta(\nabla\cdot\mathbf{v})\|^2+\mathcal{K}_1\left(\|\nabla p\|^2+\|\Delta p\|^2\right)+\\
&\mathcal{K}_2\left(\|\nabla p\|^2+\|\Delta p\|^2+\varepsilon\,\bar{p}\,\|\nabla\cdot\mathbf{v}\|^2\right)\left(\bar{p}\|\Delta p\|^2+\bar{p}^2\|\Delta \mathbf{v}\|^2\right),
\end{aligned}
\end{equation*}
which yields
{\normalsize
\begin{equation}\label{nn9s}
\begin{aligned}
&\ \frac{d}{dt}E_5(t)+D_5(t) \\
\le &\ \mathcal{K}_1\left(\|\nabla p\|^2+\|\Delta p\|^2\right)+\mathcal{K}_2\left(\|\nabla p\|^2+\|\Delta p\|^2+\varepsilon\,\bar{p}\,\|\nabla\cdot\mathbf{v}\|^2\right)\left(\bar{p}\|\Delta p\|^2+\bar{p}^2\|\Delta \mathbf{v}\|^2\right),
\end{aligned}
\end{equation}}where
{\normalsize$$
\begin{aligned}
E_5(t)&=\frac{\bar{p}}{2}\,\|\Delta p\|^2+\bar{p}^2\|\Delta \mathbf{v}\|^2+\frac{1}{2\,\bar{p}}\,\|\bar{p}\Delta p-p\Delta p\|^2+\frac{1}{2\,\bar{p}}\,\|p\Delta p\|^2,\\
D_5(t)&=\frac{\bar{p}}{2}\,\|\nabla\Delta p\|^2+\frac{1}{\bar{p}}\,\|\nabla\Delta p-p\nabla\Delta p\|^2+\frac{1}{2\,\bar{p}}\,\|p\nabla\Delta p\|^2+\varepsilon\,\bar{p}^2\|\Delta(\nabla\cdot\mathbf{v})\|^2.
\end{aligned}
$$}Note that {\normalsize$\bar{p}\,\|\Delta p\|^2+\bar{p}^2\|\Delta \mathbf{v}\|^2\le 2\,E_5(t)$}. Hence, we update {\eqref{nn9s} as
{\normalsize
\begin{equation}\label{nn9t}
\begin{aligned}
&\ \frac{d}{dt}E_5(t)+D_5(t) \\
\le &\ \mathcal{K}_1\left(\|\nabla p\|^2+\|\Delta p\|^2\right)+\mathcal{K}_2\left(\|\nabla p\|^2+\|\Delta p\|^2+\varepsilon\,\bar{p}\,\|\nabla\cdot\mathbf{v}\|^2\right)E_5(t).
\end{aligned}
\end{equation}}Applying the Gronwall inequality to \eqref{nn9t}, we find
$$
\begin{aligned}
E_5(t)\le &\ \exp\left\{\mathcal{K}_2\int_0^t \left(\|\nabla p\|^2+\|\Delta p\|^2+\varepsilon\,\bar{p}\,\|\nabla\cdot\mathbf{v}\|^2\right)d\tau\right\}\times\\
&\quad \left(E_5(0)+\mathcal{K}_1\int_0^t\left(\|\nabla p\|^2+\|\Delta p\|^2\right)d\tau\right).
\end{aligned}
$$
In view of the uniform-in-time integrability of {\normalsize$\|\nabla p\|^2$}, {\normalsize$\|\Delta p\|^2$} and {\normalsize$\varepsilon\,\bar{p}\,\|\nabla\cdot\mathbf{v}\|^2$} (cf. \eqref{M_1a}, \eqref{n32}), we see that
$E_5(t)\le \mathcal{K}_3$, for some constant $\mathcal{K}_3$ which is independent of $t$ and remains bounded as $\varepsilon\to0$. This implies that
{\normalsize\begin{equation}\label{nn9u}
\begin{aligned}
\|\Delta p(t)\|^2+\|\Delta\mathbf{v}(t)\|^2 \le 2\left(\frac{1}{\bar{p}}+\frac{1}{\bar{p}^2}\right)\mathcal{K}_3.
\end{aligned}
\end{equation}}By substituting the upper bound of $E_5(t)$ into \eqref{nn9t}, then integrating the result with respect to $t$, we have, in particular,
{\normalsize\begin{equation*}\label{nn9v}
\begin{aligned}
\int_0^t\left(\|\nabla\Delta p(\tau)\|^2+\varepsilon\,\bar{p}\,\|\Delta(\nabla\cdot\mathbf{v}(\tau))\|^2\right)d\tau\le \mathcal{K}_4
\end{aligned}
\end{equation*}}for some constant which is independent of {\normalsize$t$} and remains bounded as {\normalsize$\varepsilon\to0$}.

Next, we shall work on the third order spatial derivatives of the solution, in order to gain the desired energy estimates stated in Theorem \ref{thm2}. The proof is in exactly the same spirit as the $H^1$- and $H^2$-estimates. However, the detailed arguments involve the estimation of more than 50 nonlinear terms, whose presentation is not reader friendly. For the sake of brevity, we omit the technical details here. Moreover, by applying similar arguments as in Sections 2.5--2.8 for the three-dimensional case, we can improve the temporal integrability of {\normalsize$\|\nabla\cdot\mathbf{v}\|_{H^2}^2$}, and establish the global well-posedness, long-time behavior and zero chemical diffusion limit results for the two-dimensional case. The detailed proofs are omitted to simplify the presentation. This completes the proof for Theorem \ref{thm2}.\hfill$\square$

%\section*{Acknowledgments}
%D. Wang  was partially supported  by the National Science Foundation under grants  DMS-1613213 and DMS-1907519. F. Wang was partially supported by the Hunan Provincial Key Laboratory of Mathematical Modeling and Analysis in Engineering (No. 2017TP1017, Changsha University of Science and Technology), the Natural Science Foundation of Hunan Province (No. 2019JJ50659), and the Double First-class International Cooperation Expansion Project (No. 2019IC39). Z. Wang was supported by the the Hong Kong RGC GRF grant No. PolyU 153032/15P. K. Zhao was partially supported by the Simons Foundation Collaboration Grant for Mathematicians (No. 413028).

\section*{Appendix}

\subsection*{1. Derivation of \eqref{nn9r}}

In this appendix, we first provide the detailed derivation of \eqref{nn9r}. The proof involves the estimates of more than 30 nonlinear terms, which is achieved by using \eqref{gn2d1}--\eqref{gn2d5}, Cauchy-Schwarz and Young's inequalities. In the subsequent energy estimates, a generic constant, denoted by $\eta$, which is to be specified at the end of the proof, will appear frequently. To simplify the presentation, we use $d$ to denote a generic constant which has sole dependence on $\eta$ and the constants in \eqref{gn2d1}--\eqref{gn2d5}.
\begin{itemize}
\item For {\normalsize$R_1(t)$}, we have
{\normalsize
\begin{equation*}\label{nn9a1}
\begin{aligned}
\left|R_1(t)\right|&=2\,\bar{p}\left|\int_{\mathbb{R}^2}\left(\nabla(\nabla p\cdot\mathbf{v})\right)\cdot(\nabla\Delta p)\ d\mathbf{x}\right|
=2\,\bar{p}\left|\int_{\mathbb{R}^2}\left[\mathbb{H}(p)\cdot\mathbf{v}+(\nabla\mathbf{v})^{\mathrm{T}}\cdot\nabla p\right]\cdot(\nabla\Delta p)\ d\mathbf{x}\right|\\
&\equiv \left|\tilde{R}_{11}+\tilde{R}_{12}\right|.
\end{aligned}
\end{equation*}}

\begin{itemize}
\item For {\normalsize$\tilde{R}_{11}$}, we have
{\normalsize
\begin{equation*}\label{nn9a11}
\begin{aligned}
\left|\tilde{R}_{11}\right|&=2\,\bar{p}\left|\int_{\mathbb{R}^2}\left[\mathbb{H}(p)\cdot\mathbf{v}\right]\cdot(\nabla\Delta p)\ d\mathbf{x}\right|\\
&\le 2\,\bar{p}\,\|\mathbb{H}(p)\|_{L^4}\|\mathbf{v}\|_{L^4}\|\nabla\Delta p\|
\le d\,\bar{p}\,\|\Delta p\|^{\frac12}\|\nabla\Delta p\|^{\frac32}\|\mathbf{v}\|^{\frac12}\|\nabla\cdot\mathbf{v}\|^{\frac12}\\
&\le \eta\,\|\nabla\Delta p\|^2+d\,\bar{p}^4\,\|\Delta p\|^2\|\mathbf{v}\|^2\|\nabla\cdot\mathbf{v}\|^2
%&\le \eta\,\|\nabla\Delta p\|^2+d\,\bar{p}^4\,\|\Delta p\|^2\|\mathbf{v}\|^2\|\nabla\cdot\mathbf{v}\|^2\\
\le \eta\,\|\nabla\Delta p\|^2+d\,\bar{p}^4\,M_1\delta\,\|\Delta p\|^2,
\end{aligned}
\end{equation*}}for some positive generic constant {\normalsize$\eta$} which will be specified later.

\item For {\normalsize$\tilde{R}_{12}$}, we have
{\normalsize
\begin{equation*}\label{nn9a12}
\begin{aligned}
\left|\tilde{R}_{12}\right|&=2\,\bar{p}\left|\int_{\mathbb{R}^2}\left[(\nabla\mathbf{v})^{\mathrm{T}}\cdot\nabla p\right]\cdot(\nabla\Delta p)\ d\mathbf{x}\right|\\
&\le 2\,\bar{p}\,\|\nabla\mathbf{v}\|_{L^4}\|\nabla p\|_{L^4}\|\nabla\Delta p\|
\le d\,\bar{p}\,\|\nabla\cdot\mathbf{v}\|^{\frac12}\|\Delta \mathbf{v}\|^{\frac12}\|\nabla p\|^{\frac12}\|\Delta p\|^{\frac12}\|\nabla\Delta p\|\\
&\le \eta\,\|\nabla\Delta p\|^2+d\,\bar{p}^2\,\|\nabla\cdot\mathbf{v}\|\|\Delta \mathbf{v}\|\|\nabla p\|\|\Delta p\|
\le \eta\,\|\nabla\Delta p\|^2+d\,\bar{p}^2\,\delta\left(\|\Delta\mathbf{v}\|^2+\|\Delta p\|^2\right).
\end{aligned}
\end{equation*}}

\end{itemize}

\item For {\normalsize$R_2(t)$}, similar to the estimate of $R_{12}$, we have
{\normalsize
\begin{equation*}\label{nn9a2}
\begin{aligned}
\left|R_2(t)\right|&=2\,\bar{p}\left|\int_{\mathbb{R}^2}(\nabla\cdot\mathbf{v})\nabla p\cdot(\nabla\Delta p)\ d\mathbf{x}\right|
\le \eta\,\|\nabla\Delta p\|^2+d\,\bar{p}^2\,\delta\left(\|\Delta\mathbf{v}\|^2+\|\Delta p\|^2\right).
\end{aligned}
\end{equation*}}

\item For {\normalsize$R_3(t)$}, we have
{\normalsize
\begin{equation*}\label{nn9b1}
\begin{aligned}
\left|R_3(t)\right|&=4\,\varepsilon\,\bar{p}^2\left|\int_{\mathbb{R}^2}\left(\mathbf{v}\cdot\Delta\mathbf{v}\right)\Delta(\nabla\cdot\mathbf{v})d\mathbf{x}\right|\\
&\le 4\,\varepsilon\,\bar{p}^2\,\|\mathbf{v}\|_{L^4}\|\Delta\mathbf{v}\|_{L^4}\|\Delta(\nabla\cdot\mathbf{v})\|
\le d\,\varepsilon\,\bar{p}^2\,\|\mathbf{v}\|^{\frac12}\|\nabla\cdot\mathbf{v}\|^{\frac12}\|\Delta \mathbf{v}\|^{\frac12}\|\Delta(\nabla\cdot\mathbf{v})\|^{\frac32}\\
&\le \frac{\varepsilon\,\bar{p}^2}{2}\,\|\Delta(\nabla\cdot\mathbf{v})\|^2+d\,\varepsilon\,\bar{p}^2\|\mathbf{v}\|^2\|\nabla\cdot\mathbf{v}\|^2\|\Delta \mathbf{v}\|^2\\
&\le \frac{\varepsilon\,\bar{p}^2}{2}\,\|\Delta(\nabla\cdot\mathbf{v})\|^2+d\,M_1\,\bar{p}\,(\varepsilon\,\bar{p}\,\|\nabla\cdot\mathbf{v}\|^2)\|\Delta \mathbf{v}\|^2.
\end{aligned}
\end{equation*}}

\item For {\normalsize$R_4(t)$}, we have
{\normalsize
\begin{equation*}\label{nn9b2}
\begin{aligned}
\left|R_4(t)\right|&=4\,\varepsilon\,\bar{p}^2\left|\int_{\mathbb{R}^2}\left(|\nabla \mathbf{v}|^2\right)\Delta(\nabla\cdot\mathbf{v})d\mathbf{x}\right|
\le 4\,\varepsilon\,\bar{p}^2 \|\nabla\mathbf{v}\|_{L^4}^2\|\Delta(\nabla\cdot\mathbf{v})\|\\
&\le d\,\varepsilon\,\bar{p}^2\|\nabla\cdot\mathbf{v}\|\|\Delta\mathbf{v}\|\|\Delta(\nabla\cdot\mathbf{v})\|
\le \frac{\varepsilon\,\bar{p}^2}{2}\,\|\Delta(\nabla\cdot\mathbf{v})\|^2+d\,\bar{p}\,(\varepsilon\,\bar{p}\,\|\nabla\cdot\mathbf{v}\|^2)\|\Delta \mathbf{v}\|^2.
\end{aligned}
\end{equation*}}

\item For {\normalsize$R_5(t)$}, we have
{\normalsize
\begin{equation*}\label{nn9c}
\begin{aligned}
\left|R_5(t)\right|&=2\left|\int_{\mathbb{R}^2}\Delta p\nabla p\cdot\nabla(\nabla\cdot(p\mathbf{v}))d\mathbf{x}\right|\\
&=2\left|\int_{\mathbb{R}^2}\Delta p\nabla p\cdot\left[\mathbb{H}(p)\cdot\mathbf{v}+(\nabla\mathbf{v})^{\mathrm{T}}\cdot\nabla p+(\nabla\cdot\mathbf{v})\nabla p+p\Delta\mathbf{v}\right]d\mathbf{x}\right|\\
&\equiv 2\left|R_{51}+R_{52}+R_{53}+R_{54}\right|.
\end{aligned}
\end{equation*}}

\begin{itemize}
\item For {\normalsize$R_{51}$}, we have
{\normalsize
\begin{equation*}\label{nn9c1}
\begin{aligned}
\left|R_{51}\right|&= 2\left|\int_{\mathbb{R}^2}\Delta p\nabla p\cdot\left[\mathbb{H}(p)\cdot\mathbf{v}\right]d\mathbf{x}\right|
\le 2\,\|\Delta p\|_{L^4}\|\nabla p\|_{L^4}\|\mathbb{H}(p)\|_{L^4}\|\mathbf{v}\|_{L^4}\\
&\le d\,\|\nabla\Delta p\|\|\Delta p\|^{\frac32}\|\nabla p\|^{\frac12}\|\mathbf{v}\|^{\frac12}\|\nabla\cdot\mathbf{v}\|^{\frac12}
\le \eta\,\|\nabla\Delta p\|^2+d\,\|\Delta p\|^2\|\nabla p\|\|\Delta p\|\|\mathbf{v}\|\|\nabla\cdot\mathbf{v}\|\\
&\le \eta\,\|\nabla\Delta p\|^2+d\,(M_1\delta)^{\frac12}\left(\|\nabla p\|^2+\|\Delta p\|^2\right)\|\Delta p\|^2.
\end{aligned}
\end{equation*}}

\item For {\normalsize$R_{52}$}, we have
{\normalsize
\begin{equation*}\label{nn9c2}
\begin{aligned}
\left|R_{52}\right|&=2\left|\int_{\mathbb{R}^2}\Delta p\nabla p\cdot\left[(\nabla\mathbf{v})^{\mathrm{T}}\cdot\nabla p\right]d\mathbf{x}\right|
\le 2\,\|\nabla p\|_{L^4}^2\|\Delta p\|_{L^4}\|\nabla \mathbf{v}\|_{L^4}\\
&\le d\,\|\nabla p\|\|\Delta p\|^{\frac32}\|\nabla\Delta p\|^{\frac12}\|\nabla\cdot\mathbf{v}\|^{\frac12}\|\Delta\mathbf{v}\|^{\frac12}\\
&\le \eta\,\|\nabla\Delta p\|^2+d\,\|\nabla p\|^{\frac43}\|\Delta p\|^2\|\nabla\cdot\mathbf{v}\|^{\frac23}\|\Delta\mathbf{v}\|^{\frac23}\\
&\le \eta\,\|\nabla\Delta p\|^2+d\,\delta^{\frac13}\|\nabla p\|^{\frac43}\|\Delta p\|^{\frac43}\|\Delta p\|^{\frac23}\|\Delta\mathbf{v}\|^{\frac23}\\
&\le \eta\,\|\nabla\Delta p\|^2+d\,\delta^{\frac13}\left(\|\nabla p\|^2\|\Delta p\|^2+\|\Delta p\|^2\|\Delta\mathbf{v}\|^2\right).
\end{aligned}
\end{equation*}}

\item For {\normalsize$R_{53}$}, we have
{\normalsize
\begin{equation*}\label{nn9c3}
\begin{aligned}
\left|R_{53}\right|&=2 \left|\int_{\mathbb{R}^2}\Delta p\nabla p\cdot\left[(\nabla\cdot\mathbf{v})\nabla p\right]d\mathbf{x}\right|
\le 2\,\|\nabla p\|_{L^4}^2\|\Delta p\|_{L^4}\|\nabla\cdot \mathbf{v}\|_{L^4}\\
&\le \eta\,\|\nabla\Delta p\|^2+d\,\delta^{\frac13}\left(\|\nabla p\|^2\|\Delta p\|^2+\|\Delta p\|^2\|\Delta\mathbf{v}\|^2\right).
\end{aligned}
\end{equation*}}

\item For {\normalsize$R_{54}$}, we have
{\normalsize
\begin{equation*}\label{nn9c4}
\begin{aligned}
\left|R_{54}\right|&=2 \left|\int_{\mathbb{R}^2}\Delta p\nabla p\cdot\left[p\Delta\mathbf{v}\right]d\mathbf{x}\right|
\le 2\,\|p\|_{L^\infty}\|\nabla p\|_{L^4}\|\Delta p\|_{L^4}\|\Delta\mathbf{v}\|\\
&\le d\,\|p\|^{\frac12}\|\nabla p\|^{\frac12}\|\Delta p\|^{\frac32}\|\nabla\Delta p\|^{\frac12}\|\Delta\mathbf{v}\|
\le \eta\,\|\nabla\Delta p\|^2+d\,\|p\|^{\frac23}\|\nabla p\|^{\frac23}\|\Delta p\|^2\|\Delta\mathbf{v}\|^{\frac43}\\
&\le \eta\,\|\nabla\Delta p\|^2+d\,(M_1\delta)^{\frac13}\|\Delta p\|^{\frac23}\|\Delta p\|^{\frac43}\|\Delta\mathbf{v}\|^{\frac43}\\
&\le \eta\,\|\nabla\Delta p\|^2+d\,(M_1\delta)^{\frac13}\left(\|\Delta p\|^2+\|\Delta p\|^2\|\Delta\mathbf{v}\|^2\right).
\end{aligned}
\end{equation*}}

\end{itemize}

\item For {\normalsize$R_6(t)$}, we have
{\normalsize
\begin{equation*}\label{nn9d}
\begin{aligned}
\left|R_6(t)\right|&=2\left|\int_{\mathbb{R}^2}p\nabla\Delta p\cdot\left[\mathbb{H}(p)\cdot\mathbf{v}+(\nabla\mathbf{v})^{\mathrm{T}}\cdot\nabla p+(\nabla\cdot\mathbf{v})\nabla p+p\Delta\mathbf{v}\right]d\mathbf{x}\right|\\
&\equiv 2\left|R_{61}+R_{62}+R_{63}+R_{64}\right|.
\end{aligned}
\end{equation*}}

\begin{itemize}
\item For {\normalsize$R_{61}$}, we have
{\normalsize
\begin{equation*}\label{nn9d1}
\begin{aligned}
\left|R_{61}\right|&=2\left|\int_{\mathbb{R}^2}p\nabla\Delta p\cdot\left[\mathbb{H}(p)\cdot\mathbf{v}\right]d\mathbf{x}\right|
\le 2\|p\|_{L^\infty}\|\nabla\Delta p\|\|\mathbb{H}(p)\|_{L^4}\|\mathbf{v}\|_{L^4}\\
&\le d\,\|p\|^{\frac12}\|\Delta p\|\|\nabla\Delta p\|^{\frac32}\|\mathbf{v}\|^{\frac12}\|\nabla\cdot\mathbf{v}\|^{\frac12}
\le \eta\,\|\nabla\Delta p\|^2+d\,\|p\|^2\|\Delta p\|^4\|\mathbf{v}\|^2\|\nabla\cdot\mathbf{v}\|^2\\
&\le \eta\,\|\nabla\Delta p\|^2+d\,M_1^2\delta\,\|\Delta p\|^2\|\Delta p\|^2.
\end{aligned}
\end{equation*}}

\item For {\normalsize$R_{62}$}, we have
{\normalsize
\begin{equation*}\label{nn9d2}
\begin{aligned}
\left|R_{62}\right|&=2\left|\int_{\mathbb{R}^2}p\nabla\Delta p\cdot\left[(\nabla\mathbf{v})^{\mathrm{T}}\cdot\nabla p\right]d\mathbf{x}\right|
\le 2 \|p\|_{L^\infty}\|\nabla\Delta p\|\|\nabla p\|_{L^4}\|\nabla \mathbf{v}\|_{L^4}\\
&\le d\,\|p\|^{\frac12}\|\nabla p\|^{\frac12}\|\nabla\Delta p\|\|\Delta p\|\|\nabla\cdot\mathbf{v}\|^{\frac12}\|\Delta\mathbf{v}\|^{\frac12}
\le \eta\,\|\nabla\Delta p\|^2+c\,\|p\|\|\nabla p\|\|\Delta p\|^2\|\nabla\cdot\mathbf{v}\|\|\Delta\mathbf{v}\|\\
&\le \eta\,\|\nabla\Delta p\|^2+d\,(M_1\delta^2)^{\frac12}\|\Delta p\|^2\|\Delta\mathbf{v}\|
\le \eta\,\|\nabla\Delta p\|^2+d\,(M_1\delta^2)^{\frac12}\left(\|\Delta p\|^2+\|\Delta p\|^2\|\Delta\mathbf{v}\|^2\right).
\end{aligned}
\end{equation*}}

\item For {\normalsize$R_{63}$}, we have
{\normalsize
\begin{equation*}\label{nn9d3}
\begin{aligned}
\left|R_{63}\right|&=2\left|\int_{\mathbb{R}^2}p\nabla\Delta p\cdot\left[(\nabla\cdot\mathbf{v})\nabla p\right]d\mathbf{x}\right|
\le 2\,\|p\|_{L^\infty}\|\nabla\Delta p\|\|\nabla p\|_{L^4}\|\nabla\cdot \mathbf{v}\|_{L^4}\\
&\le \eta\,\|\nabla\Delta p\|^2+d\,(M_1\delta^2)^{\frac12}\left(\|\Delta p\|^2+\|\Delta p\|^2\|\Delta\mathbf{v}\|^2\right).
\end{aligned}
\end{equation*}}

\item For {\normalsize$R_{64}$}, we have
{\normalsize
\begin{equation*}\label{nn9d4}
\begin{aligned}
\left|R_{64}\right|&=2\left|\int_{\mathbb{R}^2}p\nabla\Delta p\cdot\left[p\Delta\mathbf{v}\right]d\mathbf{x}\right|
\le 2\,\|p\|_{L^\infty}^2\|\nabla\Delta p\|\|\Delta\mathbf{v}\|
\le d\,\|p\|\|\Delta p\|\|\nabla\Delta p\|\|\Delta\mathbf{v}\|\\
&\le \eta\,\|\nabla\Delta p\|^2+d\,\|p\|^2\|\Delta p\|^2\|\Delta\mathbf{v}\|^2
\le \eta\,\|\nabla\Delta p\|^2+d\,M_1\,\|\Delta p\|^2\|\Delta\mathbf{v}\|^2.
\end{aligned}
\end{equation*}}

\end{itemize}

\item For {\normalsize$R_7(t)$}, we have
{\normalsize
\begin{equation*}\label{nn9e}
\begin{aligned}
\left|R_7(t)\right|&=2\,\bar{p}\left|\int_{\mathbb{R}^2}\Delta p\nabla p\cdot\nabla(\nabla\cdot\mathbf{v})d\mathbf{x}\right|
\le 2\,\bar{p}\,\|\Delta\mathbf{v}\|\|\nabla p\|_{L^4}\|\Delta p\|_{L^4}
\le d\,\bar{p}\,\|\Delta\mathbf{v}\|\|\nabla p\|^{\frac12}\|\Delta p\|\|\nabla\Delta p\|^{\frac12}\\
&\le \eta\,\|\nabla\Delta p\|^2+d\,\bar{p}^\frac43\,\|\Delta\mathbf{v}\|^{\frac43}\|\nabla p\|^{\frac23}\|\Delta p\|^{\frac43}
\le \eta\,\|\nabla\Delta p\|^2+d\,\bar{p}^\frac43\left(\|\Delta p\|^2\|\Delta\mathbf{v}\|^2+\|\nabla p\|^2\right).
\end{aligned}
\end{equation*}}

\item For {\normalsize$R_8(t)$}, we have
{\normalsize
\begin{equation*}\label{nn9f}
\begin{aligned}
\left|R_8(t)\right|&=2\left|\int_{\mathbb{R}^2}\Delta p\nabla p\cdot\nabla\Delta p\ d\mathbf{x}\right|
\le 2\,\|\Delta p\|_{L^4}\|\nabla p\|_{L^4}\|\nabla\Delta p\|
\le d\,\|\nabla\Delta p\|^{\frac32}\|\nabla p\|^{\frac12}\|\Delta p\|\\
&\le \eta\,\|\nabla\Delta p\|^2+d\,\|\nabla p\|^2\|\Delta p\|^4
\le \eta\,\|\nabla\Delta p\|^2+d\,\delta\,\|\Delta p\|^2\|\Delta p\|^2.
\end{aligned}
\end{equation*}}

\item For {\normalsize$R_9(t)$}, we have
{\normalsize
\begin{equation*}\label{nn9g}
\begin{aligned}
\left|R_9(t)\right|&=\left|\int_{\mathbb{R}^2}(\Delta p)^2\left[p\nabla\cdot\mathbf{v}+\nabla p\cdot\mathbf{v}\right]d\mathbf{x}\right|
\equiv \left|R_{91}+R_{92}\right|.
\end{aligned}
\end{equation*}}

\begin{itemize}
\item For {\normalsize$R_{91}$}, we have
{\normalsize
\begin{equation*}\label{nn9g1}
\begin{aligned}
\left|R_{91}\right|&=\left|\int_{\mathbb{R}^2}(\Delta p)^2\left[p\nabla\cdot\mathbf{v}\right]d\mathbf{x}\right|
\le \|\Delta p\|_{L^4}^2\|p\|_{L^4}\|\nabla\cdot\mathbf{v}\|_{L^4}\\
&\le d\,\|\Delta p\|\|\nabla\Delta p\| \|p\|^{\frac12}\|\nabla p\|^{\frac12}\|\nabla\cdot\mathbf{v}\|^{\frac12}\|\Delta\mathbf{v}\|^{\frac12}
\le \eta\,\|\nabla\Delta p\|^2+d\,\|\Delta p\|^2\|p\|\|\nabla p\|\|\nabla\cdot\mathbf{v}\|\|\Delta\mathbf{v}\|\\
&\le \eta\,\|\nabla\Delta p\|^2+d\,(M_1\delta^2)^{\frac12}\|\Delta p\|^2\|\Delta\mathbf{v}\|
\le \eta\,\|\nabla\Delta p\|^2+d\,(M_1\delta^2)^{\frac12}\left(\|\Delta p\|^2+\|\Delta p\|^2\|\Delta\mathbf{v}\|^2\right).
\end{aligned}
\end{equation*}}

\item For {\normalsize$R_{92}$}, we have
{\normalsize
\begin{equation*}\label{nn9g2}
\begin{aligned}
\left|R_{92}\right|&=\left|\int_{\mathbb{R}^2}(\Delta p)^2\left[\nabla p\cdot\mathbf{v}\right]d\mathbf{x}\right|
\le \|\Delta p\|_{L^4}^2\|\nabla p\|_{L^4}\|\mathbf{v}\|_{L^4}\\
&\le d\,\|\Delta p\|\|\nabla\Delta p\|\|\nabla p\|^{\frac12}\|\Delta p\|^{\frac12}\|\mathbf{v}\|^{\frac12}\|\nabla\cdot\mathbf{v}\|^{\frac12}
\le \eta\,\|\nabla\Delta p\|^2+d\,\|\Delta p\|^3\|\nabla p\|\|\mathbf{v}\|\|\nabla\cdot\mathbf{v}\|\\
&\le \eta\,\|\nabla\Delta p\|^2+d\,(M_1\delta^2)^{\frac12}\|\Delta p\|^3
\le \eta\,\|\nabla\Delta p\|^2+d\,(M_1\delta^2)^{\frac12}\left(\|\Delta p\|^2+\|\Delta p\|^2\|\Delta p\|^2\right).
\end{aligned}
\end{equation*}}

\end{itemize}

\item For {\normalsize$R_{10}(t)$}, we have
{\normalsize
\begin{equation*}\label{nn9h}
\begin{aligned}
\left|R_{10}(t)\right|&=\bar{p}\left|\int_{\mathbb{R}^2}(\Delta p)^2\nabla\cdot\mathbf{v}\ d\mathbf{x}\right|
\le \bar{p}\,\|\Delta p\|_{L^4}^2\|\nabla\cdot\mathbf{v}\|
\le d\,\bar{p}\,\|\Delta p\|\|\nabla\Delta p\|\|\nabla\cdot\mathbf{v}\|\\
&\le \eta\,\|\nabla\Delta p\|^2+d\,\bar{p}^2\|\Delta p\|^2\|\nabla\cdot\mathbf{v}\|^2
\le \eta\,\|\nabla\Delta p\|^2+d\,\bar{p}^2\,\delta\,\|\Delta p\|^2.
\end{aligned}
\end{equation*}}

\item For {\normalsize$R_{11}(t)$}, we have
{\normalsize
\begin{equation*}\label{nn9i}
\begin{aligned}
\left|R_{11}(t)\right|&=\left|\int_{\mathbb{R}^2}(\Delta p)^3d\mathbf{x}\right|
\le \|\Delta p\|_{L^4}^2\|\Delta p\|
\le d\,\|\Delta p\|^2\|\nabla\Delta p\|
\le \eta\,\|\nabla\Delta p\|^2+d\,\|\Delta p\|^2\|\Delta p\|^2.
\end{aligned}
\end{equation*}}

\item For {\normalsize$R_{12}(t)$}, we have
{\normalsize
\begin{equation*}\label{nn9j}
\begin{aligned}
\left|R_{12}(t)\right|&=\frac{4}{\bar{p}}\left|\int_{\mathbb{R}^2}p\Delta p\nabla p\cdot \left[\mathbb{H}(p)\cdot\mathbf{v}+(\nabla\mathbf{v})^{\mathrm{T}}\cdot\nabla p+(\nabla\cdot\mathbf{v})\nabla p+p\Delta\mathbf{v}\right]d\mathbf{x}\right|\\
&\equiv \frac{4}{\bar{p}}\left|R_{121}+R_{122}+R_{123}+R_{124}\right|.
\end{aligned}
\end{equation*}}

\begin{itemize}
\item For {\normalsize$R_{121}$}, we have
{\normalsize
\begin{equation*}\label{nn9j1}
\begin{aligned}
\left|R_{121}\right|&=\frac{4}{\bar{p}}\left|\int_{\mathbb{R}^2}p\Delta p\nabla p\cdot \left[\mathbb{H}(p)\cdot\mathbf{v}\right]d\mathbf{x}\right|\\
&\le \frac{4}{\bar{p}}\,\|p\|_{L^\infty}\|\Delta p\|_{L^4}\|\nabla p\|_{L^4}\|\mathbb{H}(p)\|_{L^4}\|\mathbf{v}\|_{L^4}
\le \frac{4}{\bar{p}}\,\|p\|^{\frac12}\|\nabla\Delta p\|\|\Delta p\|^2\|\nabla p\|^{\frac12}\|\mathbf{v}\|^{\frac12}\|\nabla\cdot\mathbf{v}\|^{\frac12}\\
&\le \eta\,\|\nabla\Delta p\|^2+\frac{d}{\bar{p}^2}\,\|p\|\|\Delta p\|^4\|\nabla p\|\|\mathbf{v}\|\|\nabla\cdot\mathbf{v}\|
\le \eta\,\|\nabla\Delta p\|^2+\frac{d}{\bar{p}^2}\,M_1\delta\,\|\Delta p\|^2\|\Delta p\|^2.
\end{aligned}
\end{equation*}}

\item For {\normalsize$R_{122}$}, we have
{\normalsize
\begin{equation*}\label{nn9j2}
\begin{aligned}
\left|R_{122}\right|&=\frac{4}{\bar{p}}\left|\int_{\mathbb{R}^2}p\Delta p\nabla p\cdot \left[(\nabla\mathbf{v})^{\mathrm{T}}\cdot\nabla p\right]d\mathbf{x}\right|
\le \frac{4}{\bar{p}}\,\|p\|_{L^\infty}\|\nabla p\|_{L^\infty}^2\|\Delta p\|\|\nabla\cdot \mathbf{v}\|\\
&\le \frac{d}{\bar{p}}\,\|p\|^{\frac12}\|\nabla p\|\|\Delta p\|^{\frac32}\|\nabla\Delta p\|\|\nabla\cdot\mathbf{v}\|
\le \eta\,\|\nabla\Delta p\|^2+\frac{d}{\bar{p}^2}\,\|p\|\|\nabla p\|^2\|\Delta p\|^3\|\nabla\cdot\mathbf{v}\|^2\\
&\le \eta\,\|\nabla\Delta p\|^2+\frac{d}{\bar{p}^2}\,(M_1\delta^4)^{\frac12}\left(\|\Delta p\|^2\|\Delta p\|^2+\|\Delta p\|^2\right).
\end{aligned}
\end{equation*}}

\item For {\normalsize$R_{123}$}, we have
{\normalsize
\begin{equation*}\label{nn9j3}
\begin{aligned}
\left|R_{123}\right|&=\frac{4}{\bar{p}}\left|\int_{\mathbb{R}^2}p\Delta p\nabla p\cdot \left[(\nabla\cdot\mathbf{v})\nabla p\right]d\mathbf{x}\right|
\le \frac{4}{\bar{p}}\,\|p\|_{L^\infty}\|\nabla p\|_{L^\infty}^2\|\Delta p\|\|\nabla\cdot \mathbf{v}\|\\
&\le \eta\,\|\nabla\Delta p\|^2+\frac{d}{\bar{p}^2}\,(M_1\delta^4)^{\frac12}\left(\|\Delta p\|^2\|\Delta p\|^2+\|\Delta p\|^2\right).
\end{aligned}
\end{equation*}}

\item For {\normalsize$R_{124}$}, we have
{\normalsize
\begin{equation*}\label{nn9j4}
\begin{aligned}
\left|R_{124}\right|&=\frac{4}{\bar{p}}\left|\int_{\mathbb{R}^2}p\Delta p\nabla p\cdot \left[p\Delta\mathbf{v}\right]d\mathbf{x}\right|
\le \frac{4}{\bar{p}}\,\|p\|_{L^\infty}^2\|\nabla p\|_{L^\infty}\|\Delta p\|\|\Delta\mathbf{v}\|\\
&\le \frac{d}{\bar{p}}\,\|p\|\|\nabla p\|^{\frac12}\|\Delta p\|^2\|\nabla\Delta p\|^{\frac12}\|\Delta\mathbf{v}\|
\le \eta\,\|\nabla\Delta p\|^2+\frac{d}{\bar{p}^\frac43}\,\|p\|^{\frac43}\|\nabla p\|^{\frac23}\|\Delta p\|^{\frac83}\|\Delta\mathbf{v}\|^{\frac43}\\
&\le \eta\,\|\nabla\Delta p\|^2+\frac{d}{\bar{p}^\frac43}\,(M_1^2\delta)^{\frac13}\|\Delta p\|^{\frac43}\|\Delta p\|^{\frac43}\|\Delta\mathbf{v}\|^{\frac43}\\
&\le \eta\,\|\nabla\Delta p\|^2+\frac{d}{\bar{p}^\frac43}\,(M_1^2\delta)^{\frac13}\left(\|\Delta p\|^2\|\Delta p\|^2+\|\Delta p\|^2\|\Delta\mathbf{v}\|^2\right).
\end{aligned}
\end{equation*}}

\end{itemize}

\item For {\normalsize$R_{13}(t)$}, we have
{\normalsize
\begin{equation*}\label{nn9k}
\begin{aligned}
\left|R_{13}(t)\right|&=\frac{2}{\bar{p}}\left|\int_{\mathbb{R}^2}p^2\nabla\Delta p\cdot\left[\mathbb{H}(p)\cdot\mathbf{v}+(\nabla\mathbf{v})^{\mathrm{T}}\cdot\nabla p+(\nabla\cdot\mathbf{v})\nabla p+p\Delta\mathbf{v}\right]d\mathbf{x}\right|\\
&\equiv \frac{2}{\bar{p}}\left|R_{131}+R_{132}+R_{133}+R_{134}\right|.
\end{aligned}
\end{equation*}}

\begin{itemize}
\item For {\normalsize$R_{131}$}, we have
{\normalsize
\begin{equation*}\label{nn9k1}
\begin{aligned}
\left|R_{131}\right|&=\frac{2}{\bar{p}}\left|\int_{\mathbb{R}^2}p^2\nabla\Delta p\cdot\left[\mathbb{H}(p)\cdot\mathbf{v}\right]d\mathbf{x}\right|
\le \frac{2}{\bar{p}}\,\|p\|_{L^\infty}\|p\nabla\Delta p\|\|\nabla^2p\|_{L^4}\|\mathbf{v}\|_{L^4}\\
&\le \frac{d}{\bar{p}}\,\|p\|^{\frac12}\|\Delta p\|\|p\nabla\Delta p\|\nabla\Delta p\|^{\frac12}\|\mathbf{v}\|^{\frac12}\|\nabla\cdot\mathbf{v}\|^{\frac12}\\
&\le 3\,\eta\,\|p\nabla\Delta p\|^{\frac43}\|\nabla\Delta p\|^{\frac23}+\frac{d}{\bar{p}^4}\,\|p\|^2\|\Delta p\|^4\|\mathbf{v}\|^2\|\nabla\cdot\mathbf{v}\|^2\\
&\le \eta\,\|\nabla\Delta p\|^2+2\,\eta\,\|p\nabla\Delta p\|^2+\frac{d}{\bar{p}^4}\,M_1^2\delta\,\|\Delta p\|^2\|\Delta p\|^2.
\end{aligned}
\end{equation*}}

\item For {\normalsize$R_{132}$}, we have
{\normalsize
\begin{equation*}\label{nn9k2}
\begin{aligned}
\left|R_{132}\right|&=\frac{2}{\bar{p}}\left|\int_{\mathbb{R}^2}p^2\nabla\Delta p\cdot\left[(\nabla\mathbf{v})^{\mathrm{T}}\cdot\nabla p\right]d\mathbf{x}\right|\\
&\le \frac{2}{\bar{p}}\,\|p\|_{L^\infty}\|p\nabla\Delta p\|\|\nabla p\|_{L^4}\|\nabla \mathbf{v}\|_{L^4}
\le \frac{d}{\bar{p}}\,\|p\|^{\frac12}\|\nabla p\|^{\frac12}\|p\nabla\Delta p\|\|\Delta p\|\|\nabla\cdot\mathbf{v}\|^{\frac12}\|\Delta\mathbf{v}\|^{\frac12}\\
&\le \eta\,\|p\nabla\Delta p\|^2+\frac{d}{\bar{p}^2}\,\|p\|\|\nabla p\|\|\Delta p\|^2\|\nabla\cdot\mathbf{v}\|\|\Delta\mathbf{v}\|
\le \eta\,\|p\nabla\Delta p\|^2+\frac{d}{\bar{p}^2}\,(M_1\delta^2)^{\frac12}\|\Delta p\|^2\|\Delta\mathbf{v}\|\\
&\le \eta\,\|p\nabla\Delta p\|^2+\frac{d}{\bar{p}^2}\,(M_1\delta^2)^{\frac12}\left(\|\Delta p\|^2+\|\Delta p\|^2\|\Delta\mathbf{v}\|^2\right).
\end{aligned}
\end{equation*}}

\item For {\normalsize$R_{133}$}, we have
{\normalsize
\begin{equation*}\label{nn9k3}
\begin{aligned}
\left|R_{133}\right|&=\frac{2}{\bar{p}}\left|\int_{\mathbb{R}^2}p^2\nabla\Delta p\cdot\left[(\nabla\cdot\mathbf{v})\nabla p\right]d\mathbf{x}\right|
\le \frac{2}{\bar{p}}\,\|p\|_{L^\infty}\|p\nabla\Delta p\|\|\nabla p\|_{L^4}\|\nabla\cdot \mathbf{v}\|_{L^4}\\
&\le \eta\,\|p\nabla\Delta p\|^2+\frac{d}{\bar{p}^2}\,(M_1\delta^2)^{\frac12}\left(\|\Delta p\|^2+\|\Delta p\|^2\|\Delta\mathbf{v}\|^2\right).
\end{aligned}
\end{equation*}}

\item For {\normalsize$R_{134}$}, we have
{\normalsize
\begin{equation*}\label{nn9k4}
\begin{aligned}
\left|R_{134}\right|&=\frac{2}{\bar{p}}\left|\int_{\mathbb{R}^2}p^2\nabla\Delta p\cdot\left[p\Delta\mathbf{v}\right]d\mathbf{x}\right|
\le \frac{2}{\bar{p}}\,\|p\|_{L^\infty}^2\|p\nabla\Delta p\|\|\Delta\mathbf{v}\|
\le \frac{d}{\bar{p}}\,\|p\|\|\Delta p\|\|p\nabla\Delta p\|\|\Delta\mathbf{v}\|\\
&\le \eta\,\|p\nabla\Delta p\|^2+\frac{d}{\bar{p}^2}\,\|p\|^2\|\Delta p\|^2\|\Delta\mathbf{v}\|^2
\le \eta\,\|p\nabla\Delta p\|^2+\frac{d}{\bar{p}^2}\,M_1\,\|\Delta p\|^2\|\Delta\mathbf{v}\|^2.
\end{aligned}
\end{equation*}}
\end{itemize}

\item For {\normalsize$R_{14}(t)$}, we have
{\normalsize
\begin{equation*}\label{nn9l}
\begin{aligned}
\left|R_{14}(t)\right|&=4\left|\int_{\mathbb{R}^2}p\Delta p\nabla p\cdot\nabla(\nabla\cdot\mathbf{v})d\mathbf{x}\right|
\le 4\,\|p\|_{L^\infty}\|\nabla p\|_{L^\infty}\|\Delta p\|\|\Delta \mathbf{v}\|\\
&\le d\,\|p\|^{\frac12}\|\Delta p\|^{\frac12}\|\nabla p\|^{\frac12}\|\nabla\Delta p\|^{\frac12}\|\Delta p\|\|\Delta\mathbf{v}\|
\le \eta\,\|\nabla\Delta p\|^2+d\,\|p\|^{\frac23}\|\Delta p\|^{\frac23}\|\nabla p\|^{\frac23}\|\Delta p\|^{\frac43}\|\Delta\mathbf{v}\|^{\frac43}\\
&\le \eta\,\|\nabla\Delta p\|^2+d\,(M_1\delta)^{\frac13}\|\Delta p\|^{\frac23}\|\Delta p\|^{\frac43}\|\Delta\mathbf{v}\|^{\frac43}
\le \eta\,\|\nabla\Delta p\|^2+d\,(M_1\delta)^{\frac13}\left(\|\Delta p\|^2+\|\Delta p\|^2\|\Delta\mathbf{v}\|^2\right).
\end{aligned}
\end{equation*}}

\item For {\normalsize$R_{15}(t)$}, we have
{\normalsize
\begin{equation*}\label{nn9m}
\begin{aligned}
\left|R_{15}(t)\right|&=2\left|\int_{\mathbb{R}^2}p^2\nabla\Delta p\cdot\nabla(\nabla\cdot\mathbf{v})d\mathbf{x}\right|
\le 2\,\|p\|_{L^\infty}^2\|\nabla\Delta p\|\|\Delta\mathbf{v}\|
\le d\,\|p\|\|\Delta p\|\|\nabla\Delta p\|\|\Delta\mathbf{v}\|\\
&\le \eta\,\|\nabla\Delta p\|^2+d\,\|p\|^2\|\Delta p\|^2\|\Delta\mathbf{v}\|^2
\le \eta\,\|\nabla\Delta p\|^2+d\,M_1\,\|\Delta p\|^2\|\Delta\mathbf{v}\|^2.
\end{aligned}
\end{equation*}}

\item For {\normalsize$R_{16}(t)$}, we have
{\normalsize
\begin{equation*}\label{nn9n}
\begin{aligned}
\left|R_{16}(t)\right|&=\frac{4}{\bar{p}}\left|\int_{\mathbb{R}^2}p\Delta p\nabla p\cdot\nabla\Delta p\ d\mathbf{x}\right|
\le \frac{4}{\bar{p}}\,\|p\nabla\Delta p\|\|\nabla p\|_{L^4}\|\Delta p\|_{L^4}\\
&\le \frac{d}{\bar{p}}\,\|p\nabla\Delta p\|\|\nabla p\|^{\frac12}\|\Delta p\|^{\frac12}\|\Delta p\|^{\frac12}\|\nabla\Delta p\|^{\frac12}
\le 3\,\eta\,\|p\nabla\Delta p\|^{\frac43}\|\nabla\Delta p\|^{\frac32}+\frac{d}{\bar{p}^4}\,\|\nabla p\|^2\|\Delta p\|^2\|\Delta p\|^2\\
&\le \eta\,\|\nabla\Delta p\|^2+2\,\eta\,\|p\nabla\Delta p\|^2+\frac{d}{\bar{p}^4}\,\delta\, \|\Delta p\|^2\|\Delta p\|^2.
\end{aligned}
\end{equation*}}

\item For {\normalsize$R_{17}(t)$}, we have
{\normalsize
\begin{equation*}\label{nn9o}
\begin{aligned}
\left|R_{17}(t)\right|&=\frac{2}{\bar{p}}\left|\int_{\mathbb{R}^2}p(\Delta p)^2 \left[p\nabla\cdot\mathbf{v}+\nabla p\cdot\mathbf{v}\right]d\mathbf{x}\right|
\equiv \frac{2}{\bar{p}}\left|R_{171}+R_{172}\right|.
\end{aligned}
\end{equation*}}

\begin{itemize}
\item For {\normalsize$R_{171}$}, we have
{\normalsize
\begin{equation*}\label{nn9o1}
\begin{aligned}
\left|R_{171}\right|&=\frac{2}{\bar{p}}\left|\int_{\mathbb{R}^2}p(\Delta p)^2\left[p\nabla\cdot\mathbf{v}\right]d\mathbf{x}\right|
\le \frac{2}{\bar{p}}\,\|p\|_{L^\infty}^2\|\Delta p\|_{L^4}^2\|\nabla\cdot\mathbf{v}\|\\
&\le \frac{d}{\bar{p}}\,\|p\|\|\Delta p\|^2\|\nabla\Delta p\|\|\nabla\cdot\mathbf{v}\|
\le \eta\,\|\nabla\Delta p\|^2+\frac{d}{\bar{p}^2}\,M_1\delta\,\|\Delta p\|^2\|\Delta p\|^2.
\end{aligned}
\end{equation*}}

\item For {\normalsize$R_{172}$}, we have
{\normalsize
\begin{equation*}\label{nn9o2}
\begin{aligned}
\left|R_{172}\right|&=\frac{2}{\bar{p}}\left|\int_{\mathbb{R}^2}p(\Delta p)^2\left[\nabla p\cdot\mathbf{v}\right]d\mathbf{x}\right|
\le \frac{2}{\bar{p}}\,\|p\|_{L^\infty}\|\Delta p\|_{L^4}^2\|\nabla p\|_{L^4}\|\mathbf{v}\|_{L^4}\\
&\le \frac{d}{\bar{p}}\,\|p\|^{\frac12}\|\Delta p\|^2\|\nabla\Delta p\|\|\nabla p\|^{\frac12}\|\mathbf{v}\|^{\frac12}\|\nabla\cdot\mathbf{v}\|^{\frac12}
\le \eta\,\|\nabla\Delta p\|^2+\frac{d}{\bar{p}^2}\,M_1\delta\,\|\Delta p\|^2\|\Delta p\|^2.
\end{aligned}
\end{equation*}}
\end{itemize}

\item For {\normalsize$R_{18}(t)$}, we have
{\normalsize
\begin{equation*}\label{nn9p}
\begin{aligned}
\left|R_{18}(t)\right|&=2\left|\int_{\mathbb{R}^2}p(\Delta p)^2\nabla\cdot\mathbf{v}\ d\mathbf{x}\right|
\le 2\,\|p\|_{L^\infty}\|\Delta p\|_{L^4}^2\|\nabla\cdot\mathbf{v}\|
\le d\,\|p\|^{\frac12}\|\Delta p\|^{\frac32}\|\nabla\Delta p\|\|\nabla\cdot\mathbf{v}\|\\
&\le \eta\,\|\nabla\Delta p\|^2+d\,(M_1\delta^2)^{\frac12}\|\Delta p\|^3
\le \eta\,\|\nabla\Delta p\|^2+d\,(M_1\delta^2)^{\frac12}\left(\|\Delta p\|^2+\|\Delta p\|^2\|\Delta p\|^2\right).
\end{aligned}
\end{equation*}}

\item For {\normalsize$R_{19}(t)$}, we have
{\normalsize
\begin{equation*}
\begin{aligned}
\left|R_{19}(t)\right|&=\frac{2}{\bar{p}}\left|\int_{\mathbb{R}^2}p(\Delta p)^3d\mathbf{x}\right|
=\frac{2}{\bar{p}}\left|-\int_{\mathbb{R}^2}(\Delta p)^2|\nabla p|^2d\mathbf{x}-2\int_{\mathbb{R}^2}p\Delta p\nabla p\cdot\nabla\Delta p\ d\mathbf{x}\right|\\
&\le \frac{2}{\bar{p}}\,\|\Delta p\|_{L^4}^2\|\nabla p\|_{L^4}^2+\eta\,\|\nabla\Delta p\|^2+2\,\eta\,\|p\nabla\Delta p\|^2+\frac{d}{\bar{p}^4}\,\delta\,\|\Delta p\|^2\|\Delta p\|^2\\
&\le \frac{d}{\bar{p}}\,\|\nabla p\|\|\Delta p\|^2\|\nabla\Delta p\|+\eta\,\|\nabla\Delta p\|^2+2\,\eta\,\|p\nabla\Delta p\|^2+\frac{d}{\bar{p}^4}\,\delta\,\|\Delta p\|^2\|\Delta p\|^2\\
&\le 2\,\eta\,\|\nabla\Delta p\|^2+2\,\eta\,\|p\nabla\Delta p\|^2+d\left(\frac{1}{\bar{p}^2}+\frac{1}{\bar{p}^4}\right)\delta\, \|\Delta p\|^2\|\Delta p\|^2.
\end{aligned}
\end{equation*}}

\end{itemize}

By assembling the above estimates, we find that the RHS of \eqref{nn9} is controlled by
{\normalsize
\begin{equation*}
\begin{aligned}
& 30\,\eta\,\|\nabla\Delta p\|^2+7\,\eta\,\|p\nabla\Delta p\|^2+
\varepsilon\,\bar{p}^2\|\Delta(\nabla\cdot\mathbf{v})\|^2+\mathcal{K}_1\left(\|\nabla p\|^2+\|\Delta p\|^2\right)\\[2mm]
&+\mathcal{K}_2\left(\|\nabla p\|^2+\|\Delta p\|^2+\varepsilon\,\bar{p}\,\|\nabla\cdot\mathbf{v}\|^2\right)\left(\bar{p}\,\|\Delta p\|^2+\bar{p}^2\,\|\Delta \mathbf{v}\|^2\right),
\end{aligned}
\end{equation*}}where the constants $\mathcal{K}_1$ and $\mathcal{K}_2$ depend only on {\normalsize$\eta,M_1,\delta,\bar{p}$} and the constants in \eqref{gn2d1}--\eqref{gn2d5}, and therefore are independent of {\normalsize$t$} and remain bounded as {\normalsize$\varepsilon\to 0$}. By choosing {\normalsize$\eta=\min\{\bar{p}/60,1/(14\bar{p})\}$}, we arrive at \eqref{nn9r}.

%---------------------------------------------------------------------------------------------------

\subsection*{2. Explicit Examples}

Next, we provide explicit examples of initial data which fulfill the requirements of Theorems \ref{thm1}--\ref{thm2}. In the {\bf three-dimensional} case, let us consider the following functions
{\normalsize\begin{equation*}\label{f1}
\begin{aligned}
&p_0(\mathbf{x})=\left\{
\begin{aligned}
&n^{-\frac54}\left[\sin\left(\frac{r}{n}-\frac{\pi}{2}\right)+1\right]+B,&2&n\pi\le r\le 4n\pi,\\[2mm]
&\quad\quad\quad\quad\quad B,&r&\in(-\infty,2n\pi)\cup (4n\pi,\infty);
\end{aligned}
\right.\\[2mm]
&\mathbf{v}_0(\mathbf{x})=\left\{
\begin{aligned}
&n^{-\frac54}\left[\sin\left(\frac{r}{n}-\frac{\pi}{2}\right)+1\right]\cdot \frac{\mathbf{x}}{r},&2&n\pi\le r\le 4n\pi,\\[2mm]
&\quad\quad\quad\quad\quad \mathbf{0},&r&\in(-\infty,2n\pi)\cup (4n\pi,\infty),
\end{aligned}
\right.
\end{aligned}
\end{equation*}}where {\normalsize$B>0$} is any fixed constant, {\normalsize$n\in \mathbb{N}$} and {\normalsize$r=|\mathbf{x}|$}. Then there hold that {\normalsize$\nabla\times\mathbf{v}_0=\mathbf{0}$}, and
{\normalsize\begin{equation}\label{f2a}
\begin{aligned}
 \|p_0-B\|^2\cong n^{\frac12}, \quad \|\mathbf{v}_0\|^2\cong n^{\frac12}, \quad
 \|\nabla p_0\|^2\cong n^{-\frac32},  \quad\|\nabla\cdot\mathbf{v}_0\|^2\cong n^{-\frac32}.
\end{aligned}
\end{equation}}At the beginning of Section 2, we assumed that
{\normalsize
\begin{equation*}\label{f3}
\begin{aligned}
&\sup_{0\le t\le T}\left(\|p(t)-B\|^2+\|\mathbf{v}(t)\|^2\right)\le N_1, \quad
\sup_{0\le t\le T}\left(\|\nabla p(t)\|^2+\|\nabla\cdot\mathbf{v}(t)\|^2\right)\le \kappa.
\end{aligned}
\end{equation*}}As the proof proceeded, we required that {\normalsize$N_1\kappa$} to be smaller than some absolute constant and obtained the following
{\normalsize\begin{equation*}\label{f5}
\begin{aligned}
&N_1=(1+1/B)\left(\|p_1-B\|^2+B\,\|\mathbf{v}_0\|^2\right)+1, \quad
\kappa=2(1+1/B)\left[\|\nabla p_0\|^2+B\,\|\nabla\cdot\mathbf{v}_0\|^2\right].
\end{aligned}
\end{equation*}}From \eqref{f2a} we see that
%{\normalsize\begin{equation*}\label{f6}
%\begin{aligned}
$N_1\cong n^{\frac12}$, $\kappa\cong n^{-\frac32}$,
%\end{aligned}
%\end{equation*}}
from which we see that when {\normalsize$n\in\mathbb{N}$} is sufficiently large, there holds that
%{\normalsize\begin{equation*}\label{f7'}
%\begin{aligned}
$N_1\kappa\cong n^{-1}$.
%\end{aligned}
%\end{equation*}}
Hence, the smallness of {\normalsize$N_1\kappa$} can be realized by choosing {\normalsize$n\in\mathbb{N}$} to be sufficiently large.

\vspace{.1 in}

In the {\bf two-dimensional} case, let us consider the following functions
{\normalsize\begin{equation}\label{f1}
\begin{aligned}
&p_0(\mathbf{x})=\left\{
\begin{aligned}
&\frac{n}{f(n)}\left[\sin\left(\frac{r}{f(n)}-\frac{\pi}{2}\right)+1\right]+A,\quad&2&f(n)\pi\le r\le 4f(n)\pi,\\[2mm]
&\quad\quad\quad\quad\quad A,&r&\in(-\infty,2f(n)\pi)\cup (4f(n)\pi,\infty);
\end{aligned}
\right.\\[2mm]
&\mathbf{v}_0(\mathbf{x})=\left\{
\begin{aligned}
&\frac{n}{f(n)}\left[\sin\left(\frac{r}{f(n)}-\frac{\pi}{2}\right)+1\right]\cdot \frac{\mathbf{x}}{r},\quad&2&f(n)\pi\le r\le 4f(n)\pi,\\[2mm]
&\quad\quad\quad\quad\quad \mathbf{0},&r&\in(-\infty,2f(n)\pi)\cup (4f(n)\pi,\infty),
\end{aligned}
\right.
\end{aligned}
\end{equation}}where {\normalsize$A>0$} is any fixed constant, {\normalsize$n\in \mathbb{N}$}, {\normalsize$r=|\mathbf{x}|$} and {\normalsize$f(n)>0$} is to be determined. Then there hold that {\normalsize$\nabla\times\mathbf{v}_0=0$}. As the proof in Section 3 proceeded, we obtained the following qualitative relations:
{\normalsize\begin{equation*}\label{f0}
\begin{aligned}
&\bullet\ \ \ M_1=4\,(1+1/A)\ E_3(0)+1,\\[2mm]
&\bullet\ \ \ \|\nabla p(t)\|^2+\|\nabla \cdot\mathbf{v}(t)\|^2\le \frac{2\,(1+1/A)}{A^2}\,E_4(0)\,\exp\left\{\frac{c}{A}\,(2M_1+1)\,E_3(0)\right\},
\end{aligned}
\end{equation*}}where
{\normalsize\begin{equation*}\label{f0a}
\begin{aligned}
\bullet\ \ \ E_3(0)=&\ \frac12\,\|p_0-A\|^2+\frac{A}{2}\,\|\mathbf{v}\|^2-\varepsilon\int_{\mathbb{R}^2}(p_0-A)|\mathbf{v}_0|^2d\mathbf{x}-\frac{2\,\varepsilon+1}{6}\int_{\mathbb{R}^2}(p_0-A)^3d\mathbf{x}+\\[2mm]
&\ k_1\int_{\mathbb{R}^2}(p_0-A)^4d\mathbf{x}+\frac{2\,\varepsilon^2}{A}\int_{\mathbb{R}^2}(p_0-A)^2|\mathbf{v}_0|^2d\mathbf{x}+k_2\int_{\mathbb{R}^2}(v_{01})^4+(v_{02})^4d\mathbf{x}\\[4mm]
\cong&\ \|p_0-A\|^2+\|\mathbf{v}_0\|^2+\|p_0-A\|_{L^4}^4+\varepsilon^2\,\|(p_0-A)\mathbf{v}_0\|^2+\|\mathbf{v}_0\|_{L^4}^4,\\[3mm]
\bullet\ \ \ E_4(0)=&\ \frac{A}{2}\,\|\nabla p_0\|^2+A^3\|\nabla\cdot\mathbf{v}_0\|^2+\frac12\,\|A\nabla p_0-(p_0-A)\nabla p_0\|^2+\frac12\,\|(p_0-A)\nabla p_0\|^2\\[4mm]
\cong&\ \|\nabla p_0\|^2+\|\nabla\cdot\mathbf{v}_0\|^2+\|(p_0-A)\nabla p_0\|^2.
\end{aligned}
\end{equation*}}For the functions in \eqref{f1}, direct calculations show that
{\normalsize\begin{equation}\label{f2}
\begin{aligned}
E_3(0)\cong n^2+\frac{n^4}{[f(n)]^2}, \quad
 E_4(0)\cong \frac{n^2}{[f(n)]^2}+\frac{n^4}{[f(n)]^4},
\end{aligned}
\end{equation}}which imply
{\normalsize\begin{equation}\label{f9}
M_1\cong n^2+\frac{n^4}{[f(n)]^2}.
\end{equation}}Let
{\normalsize$$
\delta =\frac{2\,(1+1/A)}{A^2}\,E_4(0)\,\exp\left\{\frac{c}{A}\,(2M_1+1)\,E_3(0)\right\}.
$$}According to \eqref{f2}, we see that
{\normalsize\begin{equation}\label{f10}
\delta \cong \frac{1}{[f(n)]^2}\left(n^2+\frac{n^4}{[f(n)]^2}\right)\exp\left\{\left(n^2+\frac{n^4}{[f(n)]^2}\right)^2+\left(n^2+\frac{n^4}{[f(n)]^2}\right)\right\}.
\end{equation}}In the proof in Section 3, we required that {\normalsize$\delta$}, {\normalsize$M_1\delta$} and {\normalsize$M_1\delta^2$} are smaller than some absolute constants. From \eqref{f9} and \eqref{f10} we see that
{\normalsize\begin{equation*}\label{f6}
\begin{aligned}
M_1\delta \cong \frac{1}{[f(n)]^2}\left(n^2+\frac{n^4}{[f(n)]^2}\right)^2\exp\left\{\left(n^2+\frac{n^4}{[f(n)]^2}\right)^2+\left(n^2+\frac{n^4}{[f(n)]^2}\right)\right\}.
\end{aligned}
\end{equation*}}Hence, the smallness of {\normalsize$\delta$}, {\normalsize$M_1\delta$} and {\normalsize$M_1\delta^2$} can be realized by choosing {\normalsize$f(n)\cong e^{n^5}$} and {\normalsize$n\in\mathbb{N}$} to be sufficiently large. In addition, by direct calculations, we can show that {\normalsize$\|p_0-A\|_{L^2}^2\cong n^2$} and {\normalsize$\|\mathbf{v}_0\|_{L^2}^2\cong n^2$}. Therefore, the {\normalsize$L^2$}-norm of the zeroth frequency of the initial perturbation can be potentially large.

\section*{Acknowledgments}
D. Wang  was partially supported  by the National Science Foundation under grants  DMS-1613213 and DMS-1907519. F. Wang was partially supported by the Hunan Provincial Key Laboratory of Mathematical Modeling and Analysis in Engineering (No. 2017TP1017, Changsha University of Science and Technology), the Natural Science Foundation of Hunan Province (No. 2019JJ50659), and the Double First-class International Cooperation Expansion Project (No. 2019IC39). Z. Wang was supported by the the Hong Kong RGC GRF grant No. PolyU 153298/16P. K. Zhao was partially supported by the Simons Foundation Collaboration Grant for Mathematicians (No. 413028).

\end{document}